\documentclass[11pt]{amsart}
\def \R{\mathbb{R}}

\usepackage{amsmath,amsthm,amssymb,mathtools,color, tikz}
\usepackage[margin=1in]{geometry}
\def\co{\colon}
\usepackage{enumitem}

\usetikzlibrary{calc}
\usepackage[T1]{fontenc}
\definecolor{one}{RGB}{180,55,60}
\definecolor{two}{RGB}{35,110,180}
\definecolor{three}{RGB}{40,140,85}
\newcommand{\base}{%
 \coordinate (a) at (0,0); \coordinate (b) at (4.8,0);
 \coordinate (c) at (2.4,3.6);
 \coordinate (ab) at ($(a)!.5!(b)$);
 \coordinate (ac) at ($(a)!.5!(c)$);
 \coordinate (bc) at ($(b)!.5!(c)$);
 \coordinate (t) at (2.4,1.2);
}
\newcommand{\mesh}{%
 \draw[gray] (a)--(b)--(c)--cycle;
 \draw[gray] (a)--(t) (b)--(t) (c)--(t)
 (ab)--(t) (ac)--(t) (bc)--(t);
}
\definecolor{redgraph}{RGB}{175,45,55}
\definecolor{bluegraph}{RGB}{30,100,185}
\newcommand{\coords}{%
\coordinate (a) at (0,0);
\coordinate (b) at (2.8,-.85);
\coordinate (c) at (.7,3.3);
\coordinate (d) at (2.8,1.8);
\foreach \u/\v in {a/c,a/d,b/c,b/d}{\coordinate (\u\v) at ($(\u)!.5!(\v)$);}
\coordinate (abc) at (1.166667,.816667);
\coordinate (abd) at (1.866667,.316667);
\coordinate (acd) at (1.166667,1.7);
\coordinate (bcd) at (2.1,1.416667);
\coordinate (o) at (1.575,1.0625);
}
\newcommand{\framepic}{%
\draw[gray!60] (a)--(b)--(c)--cycle;
\draw[gray!60] (a)--(d)--(b) (c)--(d);
\foreach \v in {a,b}{\fill[redgraph] (\v) circle (2pt);}
\foreach \v in {c,d}{\fill[bluegraph] (\v) circle (2pt);}
\node[below left] at (a) {$a$};\node[below right] at (b) {$b$};
\node[above] at (c) {$c$};\node[right] at (d) {$d$};
}
\usepackage{float}
\usetikzlibrary{arrows.meta}
\definecolor{sectorone}{RGB}{190,65,65}
\definecolor{sectortwo}{RGB}{50,110,185}
\definecolor{sectorthree}{RGB}{40,145,95}

\newtheorem{lemma}{Lemma}
\newtheorem{theorem}[lemma]{Theorem}
\newtheorem{proposition}[lemma]{Proposition}
\newtheorem{remark}[lemma]{Remark}
\newtheorem{corollary}[lemma]{Corollary}

\theoremstyle{definition}
\newtheorem{example}[lemma]{Example}
\theoremstyle{plain}

\newcommand{\Int}{\operatorname{Int}}
\newcommand{\F}{\mathcal{F}}

\begin{document}

\title{Universal bounds for decompositions of PL manifolds}
\author{Rustam Sadykov}
\address{Department of Mathematics, Kansas State University}
\email{sadykov@math.ksu.edu}
\date{}

\begin{abstract}
    The geometric dimension of a nonempty compact PL manifold $N$ is the smallest nonnegative integer $d$ for which $N$ admits a handle decomposition with all handles of index at most $d$.  A PL filling of a compact PL manifold is a decomposition into compact codimension $0$ PL submanifolds, called sectors, with pairwise disjoint interiors and intersections described by standard local models. We advance Singhof's program of simplifying the sectors of a filling and their intersections. Specifically, for integers $m, \ell\ge 1$ and nonnegative integers $a_1\le\cdots\le a_\ell$, we prove that every closed PL $m$-manifold admits a filling into $\ell$ sectors whose nonempty $r$-fold intersections have geometric dimension at most $a_r$ for each $1\le r\le \ell$ if and only if
\[
a_{\lambda_1}+\cdots+a_{\lambda_q}\ge m-\ell+1
\]
for every partition $\lambda_1+\cdots+\lambda_q=\ell$ into positive integers. Under the same inequalities, any prescribed $\ell$-sector filling of the boundary of a compact PL $m$-manifold extends to a filling with the corresponding relative handle index bounds. Necessity follows from the cup-product obstructions of Courte, Moussard, Ren, and Zhou. To prove sufficiency, we develop compatible handle trades that simplify intersections in increasing order of depth. For nondecreasing tuples of nonnegative integers, our theorem resolves the PL classification problem 
posed by Courte, Moussard, Ren, and Zhou. In particular, we obtain the bound $(2,2,4)$ for $3$-sector fillings in dimension $6$ and improve the $5$-sector bound $(1,2,3,3,4)$ of Rubinstein and Tillmann to $(1,2,2,3,4)$ in dimension $8$.
\end{abstract}

\maketitle
\section{Introduction}

\subsection{Background and motivation}
Decomposing a manifold into simpler pieces makes it possible to study its
global topology through the pieces themselves and the way they are glued
together.
A classical example is a Heegaard splitting of a closed connected orientable $3$-manifold into two handlebodies~\cite{Heegaard, Scharlemann}. 
In dimension $4$, Gay and Kirby proved~\cite{GayKirby} that every closed connected orientable smooth $4$-manifold admits a trisection into three $1$-handlebodies, whose pairwise intersections are $1$-handlebodies and whose triple intersection is a closed connected surface. 
Bell, Hass, Rubinstein, and Tillmann gave an algorithm that constructs trisections from triangulations of closed orientable PL $4$-manifolds, and obtained explicit linear bounds on the trisection genus in terms of the number of $4$-simplices~\cite{BellHassRubinsteinTillmann}. Islambouli and Naylor generalized trisections of $4$-manifolds by allowing more than three cyclically arranged $1$-handlebody sectors, with adjacent sectors intersecting in $3$-dimensional handlebodies and all sectors meeting along a common closed connected surface~\cite{IslambouliNaylor}. Lambert-Cole and Miller proved that compact oriented smooth $5$-manifolds admit trisections, which can be chosen to extend prescribed trisections of the boundary~\cite{LambertColeMiller}. Their trisections have three $5$-dimensional $1$-handlebody sectors, pairwise intersections built of handles of index at most $2$, and a $3$-dimensional common intersection.

In higher dimensions, there are several ways to extend these constructions. Rubinstein and Tillmann constructed multisections of closed connected PL manifolds of dimension $n$ with $\lfloor n/2\rfloor +1$ sectors, each of which is a $1$-handlebody, such that the intersection of every nonempty proper subcollection of the sectors has a low dimensional spine~\cite{RubinsteinTillmann,RTGeneralized}. The common intersection is a closed connected manifold of dimension $\lceil n/2\rceil$. 
Ben Aribi, Courte, Golla, and Moussard introduced multisections of closed connected orientable smooth or PL $(n+1)$-manifolds into $n$ $1$-handlebodies such that the intersection of every nonempty proper subcollection is a $1$-handlebody, while the common intersection is a closed connected surface~\cite{BACGM}. 
They also established that every closed connected orientable smooth $5$-manifold admits a quadrisection~\cite[Theorem 7.3]{BACGM}. 
Courte, Moussard, Ren, and Zhou subsequently proved that every closed connected smooth or PL $(n+1)$-manifold with $n\ge 3$ admits such a multisection into $n$ $1$-handlebodies~\cite[Theorem~1.1 and \S5.4]{CMRZ}. For manifolds with boundary, Dissler developed a relative version of the Ben Aribi–Courte–Golla–Moussard theory. He proved that a compact connected oriented smooth $5$-manifold with nonempty boundary, equipped with a relative fibration of its boundary over $S^2$, admits a relative quadrisection inducing the prescribed fibration~\cite[Theorem 6.6]{Dissler}.

These constructions suggest a more general question. If the number of pieces is fixed, how simple can the pieces and all their intersections be made simultaneously? In this paper, we measure complexity by geometric dimension. For a nonempty compact PL manifold $W$ and a compact codimension $0$ PL submanifold $E\subseteq \partial W$, possibly empty, the \emph{relative geometric dimension} $\operatorname{gdim}(W, E)$ is the smallest nonnegative integer $d$ for which $W$ admits a handle decomposition relative to a collar of $E$ with all handles of index at most $d$. Taking $E=\emptyset$ gives the absolute \emph{geometric dimension} $\operatorname{gdim} W$ (see \S\ref{ss:dimensions} for conventions). 
For fixed integers $m\ge 1$ and $\ell\ge 1$, we determine all nondecreasing tuples of nonnegative integers $a_1 \le \cdots\le a_\ell$ such that every closed PL $m$-manifold $M$ admits a filling into $\ell$ sectors for which every nonempty $r$-fold intersection of sectors has geometric dimension at most $a_r$ for each $1\le r\le \ell$. We also prove a relative extension theorem with the same bounds (Theorem~\ref{th:intro-main}).

\subsection{Main theorem and consequences}

Following Takens and Singhof~\cite{Ta,Si}, we work with \emph{fillings}, which are decompositions $M=N_1\cup\cdots\cup N_\ell$ of a compact PL $m$-manifold $M$ into compact codimension $0$ PL submanifolds, called \emph{sectors}, with pairwise disjoint interiors and generic intersections described in terms of standard local models in \S\ref{s:fillings}. Empty sectors are permitted. For every nonempty $I\subseteq \{1, \ldots, \ell\}$, put  $N_I=\bigcap_{i\in I}N_i$. Whenever $N_I$ is nonempty, it is a compact PL manifold of dimension $m-|I|+1$. Neither sectors nor their intersections are required to be connected, and no orientability hypothesis is imposed. 

Singhof developed methods for simplifying the sectors of a filling and proposed extending these methods to simplify their intersections \cite[p. 387]{Si}. He observed that such an extension could improve his method for computing the Lusternik–Schnirelmann category, developed in \cite{Si14}. In this paper, we pursue this program in the PL category by adapting handle trading to control the geometric dimensions of the sectors and all their intersections simultaneously.

Let $a=(a_1,\ldots,a_\ell)$ be a nondecreasing tuple of nonnegative integers. We say that $a$ is a \emph{universal relative upper bound} in dimension $m$ if for every compact PL $m$-manifold $M$, every prescribed filling $\partial M=G_1\cup\cdots\cup G_\ell$ extends to a relative filling $M=N_1\cup\cdots\cup N_\ell$ that agrees with the product filling on a collar of the boundary and satisfies $\operatorname{gdim}(N_I,G_I)\le a_{|I|}$ for every nonempty $N_I$, where $G_I=\bigcap_{i\in I}G_i$. If this condition is required only for closed manifolds $M$, we say that $a$ is a \emph{universal upper bound} in dimension $m$. 
Although the bound $a$ is nondecreasing, the maximum geometric dimensions $\max_{|I|=r}\operatorname{gdim} N_I$ of the intersections of a given filling need not be. 
Let $\mathcal A_{m,\ell}$ consist of the nondecreasing tuples
$a=(a_1,\ldots,a_\ell)$ of nonnegative integers satisfying 
\begin{equation}\label{eq:intro-partitions}
a_{\lambda_1}+\cdots+a_{\lambda_q}\ge m-\ell+1
\qquad\text{whenever}\qquad
\lambda_1+\cdots+\lambda_q=\ell,
\end{equation}
where each $\lambda_i$ is a positive integer. Our main result is the following.

\begin{theorem}[Main theorem]\label{th:intro-main}
Let $m, \ell\ge 1$, and let $(a_1,\ldots,a_\ell)$ be a nondecreasing $\ell$-tuple of nonnegative integers. Then the following conditions are equivalent:
\begin{enumerate}
\item[\textup{(i)}] $(a_1,\ldots,a_\ell)\in \mathcal{A}_{m, \ell}$.
\item[\textup{(ii)}] $(a_1,\ldots,a_\ell)$ is a universal upper bound in dimension $m$. 
\item[\textup{(iii)}] $(a_1,\ldots,a_\ell)$ is a universal relative upper bound in dimension $m$.
\item[\textup{(iv)}] $\R P^m$ admits a filling into $\ell$ sectors satisfying $\operatorname{gdim} N_I\le a_{|I|}$ for every nonempty $N_I$. 
\end{enumerate}
\end{theorem}

For fixed $m$ and $\ell$, we call the set of coordinatewise minimal tuples in $\mathcal A_{m,\ell}$ the \emph{cohomological Pareto frontier}. By Theorem~\ref{th:intro-main}, these are precisely the optimal nondecreasing universal upper bounds.

The implications $\textup{(iii)} \Rightarrow \textup{(ii)} \Rightarrow \textup{(iv)}$ are immediate. The implication $\textup{(ii)} \Rightarrow \textup{(i)}$ follows from ~\cite[Corollary~7.17]{CMRZ} after translating their support-function condition into the partition inequalities, and $\textup{(iv)} \Rightarrow \textup{(i)}$ follows by a similar cup-product argument. Our principal contribution is $\textup{(i)} \Rightarrow \textup{(iii)}$.

In the terminology of \cite{CMRZ}, 
a PL \emph{$m$-admissible tuple} is defined in the same way as a universal upper bound, except that its entries need not be nondecreasing and may also equal $-\infty$. Courte, Moussard, Ren, and Zhou posed the problem of classifying all PL $m$-admissible tuples~\cite[Definitions~7.3--7.4 and Question~7.6]{CMRZ}. Theorem~\ref{th:intro-main} resolves this classification problem for nondecreasing PL $m$-admissible tuples of nonnegative integers.

Theorem~\ref{th:intro-main} also settles specific questions about admissible tuples and determines explicit families of optimal bounds. The following examples illustrate these consequences.

\begin{example}[Improved Rubinstein-Tillmann bounds]
By the Rubinstein--Tillmann theorem~\cite{RubinsteinTillmann}, the $k$-tuples $(1, 2, \ldots, k)$ and $(1, 2, \ldots, k-2, k-2, k-1)$ are universal upper bounds in dimensions $m=2k-1$ and $m=2k-2\ge 4$, respectively. Courte, Moussard, Ren, and Zhou already proved that $(1,2,\ldots, k)$ is coordinatewise minimal among nondecreasing universal upper bounds with $k$ entries in dimension $2k-1$~\cite[Proposition 7.18]{CMRZ}. Suppose now that $m=2k-2\ge 8$ and put $r=\lfloor k/2\rfloor$. 
In this case, Courte, Moussard, Ren, and Zhou proposed the smaller tuple 
$(1,2,\ldots,r,r,r+1,\ldots,k-1)$
as a candidate universal upper bound and asked whether it is admissible \cite[Example 7.19 and Question 7.20]{CMRZ}. Theorem~\ref{th:intro-main} proves that this tuple is a universal upper bound in dimension $m$, answering their question affirmatively in the PL category. In dimension $8$, for example, the bound improves from $(1,2,3,3,4)$ to $(1,2,2,3,4)$. To verify the partition inequalities, note that the $j$-th entry $a_j$ of this tuple is $j$ for $j\le r$ and $j-1$ for $j>r$. A partition $\lambda_1+\cdots+\lambda_q=k$ has at most one entry $\lambda_i$ greater than $r$, since $2(r+1)>k$, so the sum of the corresponding entries is at least $k-1=m-k+1$. 
\end{example} 

\begin{example}[Three-sector fillings]
For three sectors, the partition inequalities reduce to $3a_1\ge m-2$, $a_1+a_2\ge m-2$, and $a_3\ge m-2$. It follows that in dimension $m\ge4$, the \emph{cohomological Pareto frontier}, that is, the set of componentwise minimal nondecreasing tuples satisfying
these inequalities, consists of the tuples
\[
a=(b,m-2-b,m-2),\qquad
\left\lceil\frac{m-2}{3}\right\rceil\le b\le
\left\lfloor\frac{m-2}{2}\right\rfloor, 
\]
where $b$ is an integer. 
For $m=4$, the cohomological Pareto frontier consists of the single tuple $(1,1,2)$. These bounds are realized  by Gay-Kirby trisections in the orientable case~\cite{GayKirby}, and, in the nonorientable case, by the multisections of Rubinstein and Tillmann~\cite{RubinsteinTillmann}. For $m=5$, the frontier consists of $(1,2,3)$, which is realized by Rubinstein--Tillmann PL multisections~\cite{RubinsteinTillmann} and, in the orientable case, by smooth trisections of Lambert-Cole and Miller~\cite{LambertColeMiller}. 
For $m=6$, the frontier consists of $(2,2,4)$. Thus, Theorem~\ref{th:intro-main} gives a three sector filling with these bounds for every closed PL $6$-manifold, answering the PL version of~\cite[Question~7.26]{CMRZ}. 
\end{example}

Finally, for $m\ge 3$, the bound $(1,\ldots,1,2)$ with $m-1$ entries recovers the PL handle index bounds of surface centered multisections~\cite{BACGM, CMRZ}.

\subsection{Outline of the paper and relation to earlier work}
In \S\ref{s:PL_preliminaries} we recall PL tools, including regular neighborhoods and collapses. In \S\ref{s:fillings} we review fillings and relative fillings, define absolute and relative geometric dimension, and fix the conventions used throughout the paper.

In \S\ref{s:obstruction} we develop a cohomological obstruction to prescribed bounds on intersections of covering sets. The classical cup-length theorem in Lusternik--Schnirelmann category theory~\cite{CLOT} states that if a topological space $X$ has $r$ cohomology classes of positive degree whose product is nonzero, then $X$ cannot be covered by $r$ open sets each of which can be contracted to a point within $X$. Theorem~\ref{th:cup_product_obstruction} refines this obstruction by incorporating information about the cohomology of multiple intersections of the covering sets. It applies to finite open covers of topological spaces and finite covers of simplicial complexes by subcomplexes. We recover the classical cup-length theorem in Corollary~\ref{c:ls_proof}.
A PL filling admits a triangulation compatible with its sectors, so Theorem~\ref{th:cup_product_obstruction} applies to fillings of arbitrary compact PL manifolds. The resulting obstruction does not require the dimension bounds to be nondecreasing. To obtain the partition inequalities \eqref{eq:intro-partitions}, we specialize to real projective space and use the monotonicity of the upper bounds (Corollary~\ref{c:necessity}). The necessity of these inequalities for universal upper bounds already follows from the obstruction of Courte, Moussard, Ren, and Zhou~\cite[Theorem 7.16 and Corollary~7.17]{CMRZ}. For nondecreasing tuples of nonnegative integers, their support function condition for the torus is equivalent to the partition inequalities.

In \S\ref{s:canonical} we give an explicit construction of canonical fillings from triangulations. This construction is independent of the proof of Theorem~\ref{th:intro-main}.
The barycentric vertex grouping idea underlying our canonical fillings is classical. It appears in Schwarz's proof of the bound for the genus of a fiber space \cite{Schwarz}, where unions of open stars deform onto subcomplexes of controlled dimension.
In \S\ref{s:canonical} we adapt the vertex partition construction of Rubinstein and Tillmann~\cite[Sections 2.5 and 4.1]{RubinsteinTillmann}, using derived neighborhoods in barycentric subdivisions. We describe the intersection spines explicitly in terms of order complexes. For a triangulated closed PL $m$-manifold, and a partition
$\{0,\ldots,m\}=P_1\sqcup\cdots\sqcup P_\ell$ into $\ell$ nonempty blocks with $1\le \ell\le m+1$, our construction gives a canonical filling satisfying the simultaneous bounds
\[
\operatorname{gdim}N_{i_1\cdots i_r}
\le
\sum_{j=1}^{r}\bigl(|P_{i_j}|-1\bigr),
\qquad
1\le i_1<\cdots<i_r\le\ell.
\]

In \S\ref{s:compatible_neighborhoods} we construct compatible neighborhoods for general PL fillings. Ben Aribi, Courte, Golla, and Moussard constructed a compatible transverse product neighborhood for multisections with geometric dimension bounds $(1,1,\ldots,1,2)$ \cite[Lemma 2.4]{BACGM}. We formulate the neighborhood construction for general PL fillings and give coordinates adapted to a terminal handle and the induced filling of its cocore boundary. 
Our proof uses PL block-bundle theory.

In \S\ref{sec:trading} we describe the handle trades used in the sufficiency argument. Handle trading for fillings with $3$ sectors was studied by Takens~\cite[Section~IV]{Ta}. For general fillings handle trading was developed by Singhof~\cite[Section~4]{Si}.  
Singhof's construction removes a terminal handle from a sector while controlling the indices of handles attached to the receiving sectors (but not their intersections). Montejano gave a shorter proof of Singhof's ball covering results using skeleta and their duals~\cite{Montejano}. 
Courte, Moussard, Ren, and Zhou studied related modifications that simplify individual sectors and intersections~\cite[Section~3]{CMRZ}. Their argument relies on bridge position with respect to the multisection induced on the boundary of the chosen intersection. This boundary multisection has geometric dimension bounds $(1,1,\ldots,1,2)$. Our construction applies to general PL fillings without requiring these bounds or bridge position.

The main difficulty in applying handle trading is that removing a handle from one intersection changes other intersections as well.
To remove a terminal handle from an intersection $N_J$ with $1\le j=|J|<\ell$, we use a relative filling of the cocore disc extending the filling already induced on its boundary. Using the compatible product structure from \S\ref{s:compatible_neighborhoods}, we thicken the pieces of this filling and transfer them to the receiving sectors. We express the handle index bounds for the modified intersections in terms of the relative handle decompositions of the cocore filling and its intersections. As a result, the chosen handle disappears from $N_J$, while intersections $N_I$ with $I\neq J$ change only as described in Proposition~\ref{prop:two-receiver-trade}.

In \S\ref{sec:pareto-realization} we  prove the relative theorem for manifolds $M$ of all nonnegative dimensions $m$ simultaneously by induction on the number of labels.  Given $a\in \mathcal{A}_{m, \ell}$, we start with an extension of the prescribed boundary filling, keep the filling fixed on a collar neighborhood of $\partial M$, and simplify intersections in increasing order of depth. The partition inequalities ensure that the relative induction hypothesis applies to the cocore discs, supplying the fillings needed for the handle trades.

With these choices of cocore fillings, every handle introduced in an
intersection of depth $r\le j$ has index at most $a_r$, so the bounds already
achieved are preserved. Intersections containing all donor labels and at least one receiving label have greater depth and are treated later. Consequently, each trade decreases by at least one the total
number of handles of index greater than $a_j$ in the chosen handle
decompositions of the intersections of depth $j$, and no trade creates such
handles at depth at most $j$. Thus the process terminates and produces a filling with the required bounds. The proof of the main theorem is completed at the end of \S\ref{sec:pareto-realization}.

In \S\ref{s:ls_estimate} we bound the Lusternik–Schnirelmann category of a compact connected PL manifold in terms of the geometric dimensions of the sectors of a filling. Finally, in the appendix we verify that the derived neighborhoods used in the canonical construction satisfy the required local models of a filling.

\subsection*{Acknowledgments}
I thank Zachary Marcum for helpful discussions about this project.

\section{PL preliminaries}\label{s:PL_preliminaries}

We use the standard terminology of PL topology \cite{RS}. A \emph{polyhedron} is the underlying space $|K|$ of a locally finite simplicial complex $K$. A triangulation of a PL manifold $M$ is a homeomorphism $h\co |K|\to M$ from a polyhedron. It is \emph{compatible} with the PL structure if $h$ is a PL homeomorphism. All triangulations used in our PL constructions are assumed to be compatible. Once a triangulation is chosen, we identify $|K|$ with $M$. If $M$ is compact, then $K$ is finite. A subset $X \subseteq M$ is a subpolyhedron if it is represented by a subcomplex in some compatible triangulation of $M$.
Every smooth manifold admits a smooth triangulation, and the resulting PL structure is unique up to PL equivalence \cite{Wh}, \cite[Lectures 3-5]{Lu}.

Let $P$ and $Q$ be polyhedra such that $P=Q\cup B^k$ and $Q\cap B^k=B^{k-1}$, where $B^k$ is a PL ball of dimension $k$ and $B^{k-1}$ is a ball contained in $\partial B^k$. Then we say that $P$ \emph{elementarily collapses} onto $Q$. We write $P\searrow Q$ and say that $P$ \emph{collapses} onto $Q$ if $Q$ is obtained from $P$ by a finite sequence of elementary collapses. We will use the following characterization of regular neighborhoods. 

\begin{theorem}[{\cite[Corollary 3.30]{RS}}]\label{th:RS_collapsing}
Let $M$ be a PL manifold, and let $X\subseteq \Int M$ be a compact subpolyhedron. Let $N\subset \Int M$ be a compact  polyhedral neighborhood of $X$. Then $N$ is a regular neighborhood of $X$ in $M$ if and only if $N\searrow X$ and $N$ is a compact PL manifold with boundary.     
\end{theorem}

For a finite poset $\mathcal F$, its \emph{order complex} $\Delta(\mathcal F)$ is the simplicial complex whose vertices are the elements of $\mathcal F$ and whose nonempty simplices are its nonempty chains.
An order preserving map $\varphi\co \F\to\F$ of a poset is a \emph{descending closure operator} if $\varphi^2=\varphi$ and $\varphi(x)\le x$ for all $x\in \F$. Every descending closure operator induces a strong deformation retraction of the order complex $\Delta(\F)$ onto $\Delta(\varphi(\F))$, see, e.g., \cite[Corollary 10.12]{Bj}. We will use a stronger statement for finite posets. 

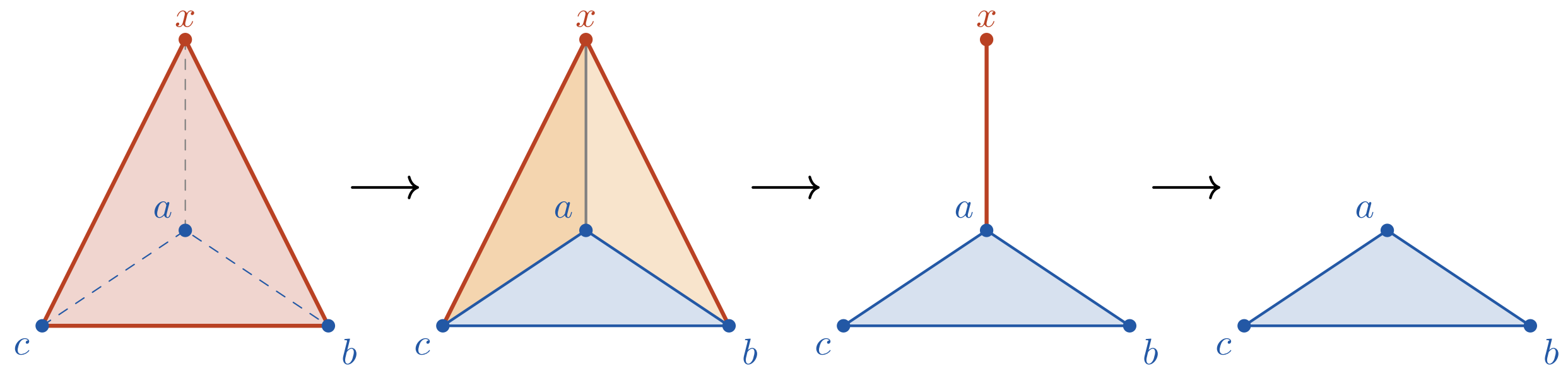
\begin{figure}[h]
    \centering
\begin{tikzpicture}[line cap=round,line join=round]
\definecolor{keep}{RGB}{35,88,165}
\definecolor{remove}{RGB}{185,65,35}
\definecolor{side}{RGB}{235,178,110}

\foreach \stage/\shift in {0/0,1/4.2,2/8.4,3/12.6}{
  \begin{scope}[xshift=\shift cm]
    \coordinate (a) at (1.5,1);
    \coordinate (b) at (3,0);
    \coordinate (c) at (0,0);
    \coordinate (x) at (1.5,3);

    \fill[keep!18] (a)--(b)--(c)--cycle;

    \ifnum\stage=0
      \fill[remove!22] (x)--(b)--(c)--cycle;
      \draw[keep,dashed] (c)--(a)--(b);
      \draw[gray,dashed] (a)--(x);
      \draw[remove,very thick] (x)--(b)--(c)--cycle;
    \fi

    \ifnum\stage=1
      \fill[side!35] (a)--(x)--(b)--cycle;
      \fill[side!55] (a)--(x)--(c)--cycle;
      \draw[keep,thick] (a)--(b)--(c)--cycle;
      \draw[gray,thick] (a)--(x);
      \draw[remove,very thick] (b)--(x)--(c);
    \fi

    \ifnum\stage=2
      \draw[keep,thick] (a)--(b)--(c)--cycle;
      \draw[remove,very thick] (a)--(x);
    \fi

    \ifnum\stage=3
      \draw[keep,thick] (a)--(b)--(c)--cycle;
    \fi

    \foreach \v in {a,b,c}{
      \fill[keep] (\v) circle (2pt);
    }
    \node[above left,keep] at (a) {$a$};
    \node[below right,keep] at (b) {$b$};
    \node[below left,keep] at (c) {$c$};

    \ifnum\stage<3
      \fill[remove] (x) circle (2pt);
      \node[above,remove] at (x) {$x$};
    \fi
  \end{scope}
}

\foreach \shift in {0,4.2,8.4}{
  \draw[->,thick]
    (\shift+3.25,1.45)--(\shift+3.95,1.45);
}
\end{tikzpicture}
 \caption{The order complex of the chain $a<x<b<c$, where $\varphi$ fixes every element except $x$ and $\varphi(x)=a$. The first arrow removes $([x,b,c], [a,x,b,c])$. The middle arrow represents two elementary collapses removing $([x,b], [a,x,b])$ and $([x,c], [a,x,c])$. The final arrow removes $([x], [a,x])$.}
    \label{fig:collapse}
\end{figure}

\begin{proposition}[{\cite[Theorem~2.1]{Ko}}]\label{prop:collapse} 
    Let $\F$ be a finite poset and $\varphi$ a descending closure operator on $\F$. Then $\Delta(\F)\searrow \Delta(\varphi(\F))$.
\end{proposition}
\begin{proof} If every element is fixed, there is nothing to prove. Otherwise, 
let $x$ be a minimal element that is not fixed by $\varphi$. Put $a=\varphi(x)$. Then $a<x$, and $a$ is fixed. Since every $y<x$ is fixed, we have $y=\varphi(y)\le\varphi(x)=a$.
Thus every chain containing $x$ but not $a$ can be enlarged by inserting $a$ in the chain immediately before $x$. Consequently, if $\sigma$ is a simplex in $\Delta(\F)$ containing $x$ but not $a$, then $\sigma$ and $a$ span a simplex $\sigma_+$ in $\Delta(\F)$. Remove all such pairs $(\sigma, \sigma_+)$ in decreasing order of $\dim\sigma$, see Figure~\ref{fig:collapse}. 
Each removal is an elementary collapse. Indeed, when $(\sigma,\ \sigma_+)$ is removed, the only remaining proper coface of $\sigma$ in $\Delta(\F)$ is $\sigma_+$, since all other cofaces belong to pairs removed earlier. These pairs remove exactly the simplices containing $x$, producing 
$
\Delta(\mathcal F)\searrow\Delta(\mathcal F\setminus\{x\}).
$
The restriction of $\varphi$ to $\F\setminus \{x\}$ is again a descending closure operator with the same image. Repeating until only fixed points of $\varphi$ remain proves the proposition. 
\end{proof}

\begin{lemma}\label{l:forNJ}
    Let $N$ be a nonempty compact PL manifold and $X\subseteq N$ a subpolyhedron. If $N\searrow X$, then $N$ admits a handle decomposition with handles of index at most $\dim X$. 
\end{lemma}
\begin{proof} Attaching outward collars to $\partial N$ produces $M=N\cup_{\partial N} (\partial N\times [0,2])$ and $N^+=N\cup_{\partial N}(\partial N\times [0,1])$. Then $N^+$ is a compact polyhedral neighborhood of $X$ in the interior of $M$. Since $N^+\searrow N\searrow X$, the collapsing criterion (Theorem~\ref{th:RS_collapsing}) implies that $N^+$ is a regular neighborhood of $X$ in $M$. Construct another regular neighborhood $N'$ of $X$ in $M$ by thickening the simplices of a triangulation of $X$ in increasing dimension. Under this construction, each $k$-simplex in $X$ contributes a $k$-handle. Hence $\operatorname{gdim}N'\le \dim X$. By uniqueness of regular neighborhoods~\cite[Theorem~3.24]{RS}, $N^+$ is PL homeomorphic to $N'$. Since $N^+$ is PL homeomorphic to $N$, we conclude that $\operatorname{gdim} N\le \dim X$. 
\end{proof}

\section{Smooth and PL fillings}\label{s:fillings}

A filling decomposes a manifold into codimension zero submanifolds whose intersections satisfy a prescribed local model. Fillings were introduced by Takens~\cite{Ta} for three pieces and extended by Singhof~\cite{Si} to arbitrary order. 

Throughout this section, we work either in the smooth category or in the PL category. In the smooth category, sectors of a filling are submanifolds with corners, while in the PL category they are PL submanifolds with boundary. All coordinate charts are understood to belong to the chosen category.  Neither sectors nor their intersections are assumed to be connected. 

\subsection{Fillings}
Let $M$ be a closed smooth or PL manifold of dimension $m$. A \emph{filling} of $M$ is an ordered collection $(N_1, \ldots, N_\ell)$ of compact codimension $0$ submanifolds of $M$, possibly empty, in the chosen category such that $M=N_1\cup \cdots \cup N_\ell$, their interiors are pairwise disjoint, and the manifolds $N_i$, called \emph{sectors}, meet according to the standard local model. Specifically, suppose that $z$ belongs to the sectors $N_{k_1}, \ldots, N_{k_\nu}$, with $k_1<k_2<\cdots<k_{\nu}$, and does not belong to any other sector. Then $\nu\le m+1$. Furthermore, there is a sufficiently small coordinate disc neighborhood $D$ centered at $z$ with coordinates $(x_1, \ldots, x_m)$ such that
    \begin{equation}\label{eq:ts1}
        D\cap N_{k_j} = \{x_1\ge 0, \ldots, x_{j-1}\ge 0, x_j\le 0\}  \quad \text{for} \quad 1\le j<\nu, 
    \end{equation}
    and
    \begin{equation}\label{eq:ts2}
       D\cap N_{k_\nu}=\{x_1\ge 0, \ldots, x_{\nu-1}\ge 0\}. 
    \end{equation}
    When $\nu=1$, (\ref{eq:ts2}) imposes no inequalities, so $D\subset N_{k_1}$.

\begin{figure}[h]
\centering
    \begin{tikzpicture}[font=\small,line cap=round,line join=round]
\definecolor{first}{RGB}{70,112,180}
\definecolor{second}{RGB}{218,142,54}
\definecolor{third}{RGB}{65,150,116}
\fill[first!18] (-2,-1.7) rectangle (0,1.7);
\fill[second!22] (0,-1.7) rectangle (2,0);
\fill[third!22] (0,0) rectangle (2,1.7);
\draw[thick] (0,-1.7)--(0,1.7) (0,0)--(2,0);
\draw[->,gray] (-2.15,0)--(2.35,0) node[right,black] {$x_1$};
\draw[->,gray] (0,-1.85)--(0,2.05) node[above,black] {$x_2$};
\fill (0,0) circle (2pt);\node[below left] at (0,0) {$z$};
\node[first,align=center] at (-1.05,.65) {$N_{1}$\\$x_1\le0$};
\node[third,align=center] at (1.05,.9) {$N_{3}$\\$x_1,x_2\ge0$};
\node[second,align=center] at (1.05,-.9) {$N_{2}$\\$x_1\ge0,\ x_2\le0$};
\begin{scope}[xshift=6.4cm]
\fill[first!18] (-2,-1.7) rectangle (0,1.7);
\fill[second!22] (0,-1.7) rectangle (2,0);
\fill[third!22] (0,0) rectangle (2,1.7);
\draw[very thick] (0,-1.7)--(0,1.7);
\draw[thick] (0,0)--(2,0);
\draw[->,gray] (-2.15,0)--(2.35,0) node[right,black] {$t$};
\draw[->,gray] (0,-1.85)--(0,2.05) node[above,black] {$y$};
\fill (0,0) circle (2pt);\node[below left] at (0,0) {$z$};
\node[first,align=center] at (-1.05,.6) {$N'_1$\\$t\le0$};
\node[third,align=center] at (1.05,.9) {$N_{2}$\\$t\ge0,\ y\ge0$};
\node[second,align=center] at (1.05,-.9) {$N_{1}$\\$t\ge0,\ y\le0$};
\end{scope}
\end{tikzpicture}
\caption{Left: the standard local model for three ordered sectors meeting at an interior point $z\in M$. Right: the standard local model for two ordered sectors meeting at a boundary point $z\in \partial M$ of the manifold $M=\{t\ge 0\}$.}
\end{figure}
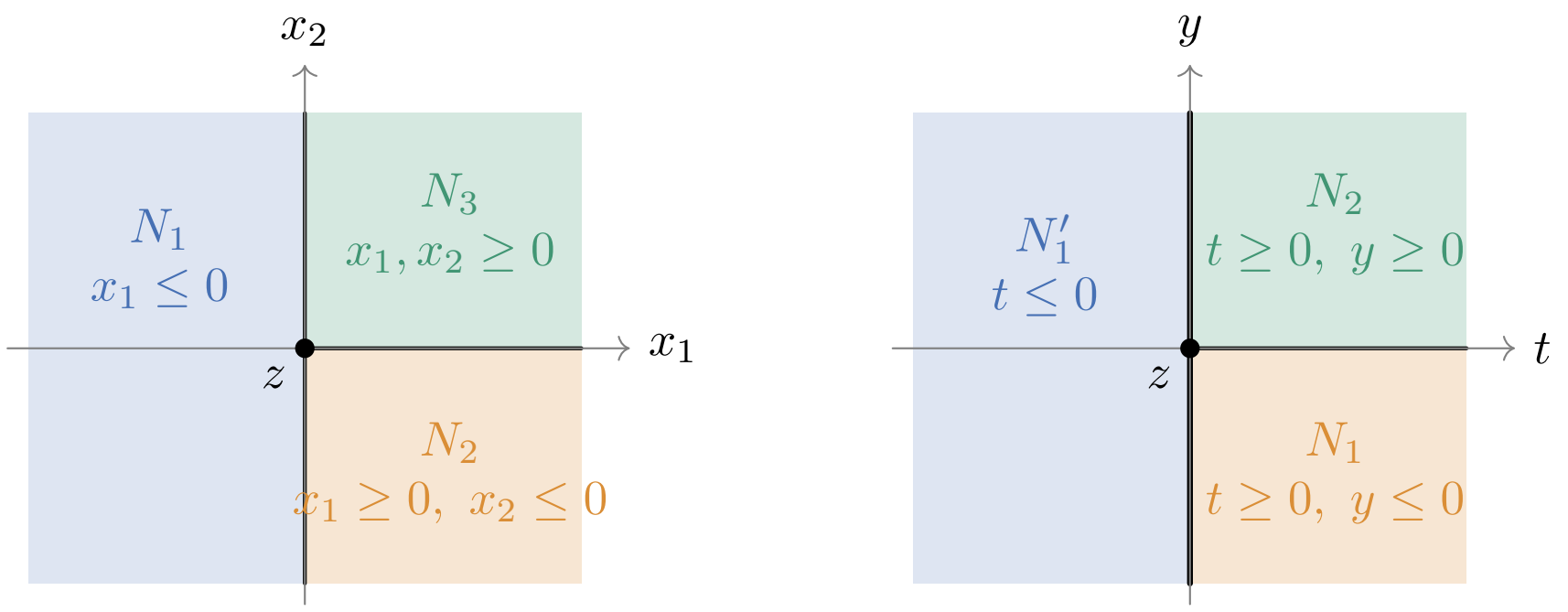

    We say that $\ell$ is the \emph{order} of the filling $(N_i)$, even if some of the sectors are empty. Every nonempty intersection of $r$ distinct sectors is a submanifold of $M$ of dimension $m-r+1$ with corners in the smooth category and with boundary in the PL category. The ordering of the sectors matters in the smooth category since a permutation of the sectors need not define another filling. In the PL category, every permutation defines another filling, since permutations of the standard local sectors can be realized by PL homeomorphisms.    

\subsection{Relative fillings}
   Relative smooth fillings were introduced in \cite{ST}. We use the analogous notion in both categories. Let $M$ be a compact manifold with possibly empty boundary. Attaching an outward collar $\partial M\times [0,1]$ to $M$
   results in a manifold
   \[
      \tilde M = M\cup_{\partial M\sim \partial M\times \{0\}} (\partial M\times [0,1]). 
   \]
   An ordered family $(N_1, \ldots, N_\ell)$ of compact codimension zero submanifolds of $M$, possibly empty, in the chosen category is a \emph{relative filling} of $M$ if the ordered family $(N_1', \ldots, N_{\ell+1}')$, where $N_1'=\partial M\times [0,1]$ and $N_i'=N_{i-1}$ for $i=2, \ldots, \ell+1$, covers $\tilde{M}$, has pairwise disjoint interiors, and satisfies the local model \eqref{eq:ts1}--\eqref{eq:ts2} at every point of $M\subset \tilde{M}$. 

   Put $G_i=N_i\cap \partial M$. These boundary pieces may be empty or disconnected. With the original labels retained, they form the induced filling of $\partial M$. For a nonempty set of labels $I$, we write $N_I=\bigcap_{i\in I}N_i$ and $G_I=N_I\cap \partial M$. We say that $|I|$ is the \emph{depth} of the intersection $N_I$. 
   Every nonempty intersection $N_I$ of sectors is a compact manifold of dimension $m-|I|+1$, with corners in the smooth category and with boundary in the PL category, and $\partial N_I=G_I\cup\bigcup_{i\notin I}N_{I\cup\{i\}}$.

\subsection{Dimensions}\label{ss:dimensions}

Fix a field $\mathbb{F}$. Given a topological space $X$, the \emph{cohomological dimension} $\operatorname{cd}_{\Bbb F}(X)$ over $\Bbb F$ is the supremum of degrees $q$ for which $H^q(X;\Bbb F)\neq 0$, where cohomology means ordinary singular cohomology. The \emph{homotopy dimension} $\operatorname{hdim}(X)$ is the least nonnegative integer $d$ such that $X$ is homotopy equivalent to a CW complex of dimension at most $d$, and is $+\infty$ if no such integer exists. For a compact manifold $W$ and a compact codimension zero submanifold $E\subseteq \partial W$, possibly empty or disconnected, the \emph{relative geometric dimension} $\operatorname{gdim}(W, E)$ is the least nonnegative integer $d$ such that $W$ admits a handle decomposition relative to a collar of $E$, with all handles of index at most $d$. Additional $0$-handles in the relative handle decomposition are allowed, so $W$ may have components disjoint from $E$. 
We set $\operatorname{gdim}(\emptyset)=\operatorname{gdim}(\emptyset, \emptyset)=\operatorname{hdim}(\emptyset)=-\infty$. 
We put $\operatorname{gdim}(W)=\operatorname{gdim}(W, \emptyset)$. In the smooth category, in the definition of relative handle decompositions corners are rounded before taking relative handle decompositions.  
For every nonempty compact manifold $W$, we have $\operatorname{cd}_{\Bbb F}(W)\le \operatorname{hdim}(W)\le \operatorname{gdim} (W)\le \dim (W)$. 

Let $\mathcal{N}=(N_1, \ldots, N_\ell)$ be a filling of a compact possibly nonorientable manifold of dimension $m$. The \emph{cohomological dimension} $\operatorname{cd}_{\mathbb{F}}(\mathcal{N})$ of a filling $\mathcal{N}=\{N_i\}$ over $\mathbb{F}$ is the maximum of the numbers $\operatorname{cd}_{\Bbb F}(N_J)+|J|-1$ over all nonempty subsets $J\subseteq \{1, \ldots, \ell\}$ with $N_J\neq \emptyset$. 

\subsection{Dimension profiles}\label{ss:dimension_profiles}

The \emph{geometric dimension profile} of $\mathcal{N}$ is the tuple $(d_1, \ldots, d_\ell)$, where $d_r=\max_{|I|=r}\operatorname{gdim}N_I$. Similarly, the \emph{homotopy dimension profile} is a tuple $(h_1, \ldots, h_\ell)$, where $h_r$ is the maximum of $\operatorname{hdim}N_I$ over all intersections $N_I$ of depth $r$. The maxima include empty intersections, so $d_r=h_r=-\infty$ if every intersection of depth $r$ is empty. In general, we have $h_r\le d_r$ for $r=1, \ldots, \ell$. For a relative filling, its relative geometric dimension profile is $(d_1^{\mathrm{rel}},\ldots,d_\ell^{\mathrm{rel}})$, where $d_r^{\mathrm{rel}}=\max_{|I|=r}\operatorname{gdim}(N_I,G_I)$. In particular, $d_r^{\mathrm{rel}}=-\infty$ if every such intersection is empty.  

The profile of a filling need not be nondecreasing. We consider nondecreasing tuples $(a_1,\ldots,a_\ell)$ of nonnegative integers and fillings $\mathcal{N}$ with $d_r(\mathcal{N})\le a_r$ for every $r$. We say that a filling $\mathcal{N}$ \emph{satisfies} the bounds $(a_1, \ldots, a_\ell)$ if $d_r(\mathcal{N})\le a_r$ for every $r$. For relative fillings, we require $d_r^{\mathrm{rel}}(\mathcal{N})\le a_r$ for every $r$.

For tuples $a,b\in\mathbb Z_{\ge0}^\ell$, we write $a\le b$ if $a_r\le b_r$ for every $r$. Let $\mathcal A\subseteq\mathbb Z_{\ge0}^\ell$ be a set of tuples. A tuple $a\in\mathcal A$ is \emph{Pareto minimal} if there is no $b\in\mathcal A$ distinct from $a$ such that $b\le a$. The set of Pareto minimal tuples in $\mathcal A$ is its \emph{Pareto frontier}.

\section{Cohomological obstructions}\label{s:obstruction}

Let $X$ be a nonempty topological space and let
$\mathcal U=(U_1,\ldots,U_\ell)$ be a finite open cover of $X$. Alternatively, let $X$ be a simplicial complex covered by finitely many subcomplexes $U_i$.
For $\emptyset\ne I\subseteq\{1,\ldots,\ell\}$, put
$U_I=\bigcap_{i\in I}U_i$, and let $s$ denote the largest $|I|$ for which $U_I\neq\emptyset$. Fix a field $\mathbb{F}$. Unless otherwise stated, cohomology means singular cohomology with coefficients in
$\mathbb F$.

 \subsection{Mayer-Vietoris spectral sequence}

We suppress the coefficients $\mathbb F$ from the notation.
Define the Mayer--Vietoris bicomplex by $C^{p,q}=\bigoplus_{|I|=p+1}C^q(U_I)$,  with $p,q\ge 0$.
Its horizontal differential $\delta$ is the alternating restriction map,
and its vertical differential $d$ is the cochain differential. Its total
complex $\operatorname{Tot}^n C=\bigoplus_{p+q=n}C^{p,q}$ has differential
 $D=\delta+(-1)^p d\quad\text{on }C^{p,q}$. The filtration $F^p\operatorname{Tot}^n C=\bigoplus_{i\ge p}C^{i,n-i}$
induces a filtration
$H^n(X)=F^0H^n(X)\supseteq F^1H^n(X)\supseteq\cdots\supseteq F^sH^n(X)=0$
and gives the Mayer--Vietoris spectral sequence, concentrated in the
columns $0\le p\le s-1$, with
\begin{equation}\label{eq:spec}
  E_1^{p,q}=\bigoplus_{|I|=p+1}H^q(U_I)
 \Longrightarrow H^{p+q}(X),
 \qquad
 E_\infty^{p,n-p}=F^pH^n(X)/F^{p+1}H^n(X).
\end{equation}

The product on the total complex respects the filtration and induces the
cup product, so
$F^aH^u(X)\cdot F^bH^v(X)\subseteq F^{a+b}H^{u+v}(X)$. Finally, 
\begin{equation}\label{eq:MV-edge-kernel}
 F^1H^n(X)=\ker\Bigl(H^n(X)\longrightarrow\bigoplus_{i=1}^{\ell}H^n(U_i)\Bigr),
\end{equation}
since $\operatorname{Tot} C/F^1\operatorname{Tot} C$ is the $0$-th column and the
induced map on cohomology is the restriction map.

\subsection{Cup product obstruction}
We now use the multiplicative structure of the Mayer--Vietoris spectral
sequence. Recall that the filtration on $H^*(X)$ is multiplicative and that
$F^sH^*(X)=0$. The following theorem gives a cup product obstruction for
covers and contains the classical cup-length theorem as a special case.

\begin{theorem}[Cup product obstruction]\label{th:cup_product_obstruction} Let $s$ denote the largest $|I|$ for which $U_I\neq\emptyset$. Let $u_1, \ldots, u_t$ be cohomology classes of $X$ of positive degrees $n_1, \ldots, n_t$, respectively. Suppose that each $u_j$ restricts to zero on every $U_i$, and that there are integers $1\le p_1,\ldots,p_t\le\ell$ such that $H^{n_j-r+1}(U_I)=0$ whenever $2\le |I|=r\le p_j$. If
$p_1+\cdots+p_t\ge s$, then $u_1\smile\cdots\smile u_t=0$.
\end{theorem}
\begin{proof}
Since $u_j$ restricts to zero on each $U_i$, we have $u_j\in F^1H^{n_j}(X)$ by \eqref{eq:MV-edge-kernel}.
The condition on intersections implies
\[
E_1^{r-1,n_j-r+1}
=
\bigoplus_{|I|=r}H^{n_j-r+1}(U_I)
=0, 
\qquad 2\le r\le p_j.
\]
These entries remain zero on every subsequent page. By \eqref{eq:spec}, $F^1H^{n_j}(X)
=
F^{p_j}H^{n_j}(X)$. Therefore, $u_1\smile\cdots\smile u_t
\in
F^{p_1+\cdots+p_t}H^{*}(X)=0$, since $p_1+\cdots+p_t\ge s$.
\end{proof}

Recall that the cup-length $\operatorname{cl}_{\mathbb{F}}X$ of a topological space is the supremum of the numbers of positive degree factors in a nonzero cup product in $H^*(X; \mathbb{F})$. It is zero if there are no nonzero positive degree classes, and it may be infinite. We define $\operatorname{cat}(X)$ as the least number of open subsets covering $X$ whose inclusions in $X$ are nullhomotopic. If no finite such cover exists, then $\operatorname{cat}(X)=\infty$.

\begin{corollary}[Cup-length bound, {\cite[Proposition~1.5]{CLOT}}]\label{c:ls_proof}
 Let $X$ be a nonempty topological space. Then $\operatorname{cat}(X)\ge \operatorname{cl}_{\mathbb F}(X)+1$.
\end{corollary}
\begin{proof} We may assume that $\ell=\operatorname{cat}(X)$ is finite.
    Let $X=U_1\cup\cdots\cup U_\ell$ be a cover by open sets whose inclusions into
$X$ are nullhomotopic. Every positive degree cohomology class restricts to zero on each $U_i$. Apply Theorem~\ref{th:cup_product_obstruction} to any $\ell$ classes of positive degree, with  $p_1=\ldots=p_\ell=1$. The condition on intersections in Theorem~\ref{th:cup_product_obstruction} is vacuous, and $s\le \ell$, so every product of $\ell$ positive degree classes vanishes. Hence $\operatorname{cl}_{\mathbb F}(X)\le\ell-1$. 
\end{proof}

\subsection{Applications to fillings}

Let $M$ be a compact, connected, possibly nonorientable PL manifold of dimension $m$, and let $\mathcal N=(N_1,\ldots,N_\ell)$ be a filling of $M$. Choose a finite triangulation compatible with the sectors. The preceding results apply to the cover by the subcomplexes $N_i$. 

\begin{proposition}\label{prop:rpm}
     Let $m\ge 1$, and let $N_1,\ldots,N_\ell$ be a filling of $\mathbb{RP}^m$. Let $a=(a_1,\ldots,a_\ell)$ be a nondecreasing tuple of nonnegative integers such that $\operatorname{gdim}N_I\le a_{|I|}$ for every $I$ with $N_I\neq \emptyset$. Then $a\in \mathcal{A}_{m, \ell}$. 
\end{proposition}
\begin{proof} Let $x\in H^1(\mathbb{RP}^m;\mathbb Z_2)$ be the generator. Consider a partition $\lambda_1+\cdots+\lambda_t=\ell$, and put $n_j=a_{\lambda_j}+\lambda_j$. For $1\le r\le\lambda_j$, we have
$n_j-r+1 = a_{\lambda_j}+\lambda_j-r+1 > a_r$, since $a$ is nondecreasing. 
Consequently, every intersection of depth $r$ has trivial cohomology in degree $n_j-r+1$, since its geometric dimension is at most $a_r$. In particular, $x^{n_j}$ restricts to zero on every sector. Theorem~\ref{th:cup_product_obstruction} therefore applies to the classes $x^{n_j}$, with $p_j=\lambda_j$, and gives $x^{n_1+\cdots+n_t}=0$.
Since \(x^m\ne0\), we obtain
\[
m+1
\le
\sum_{j=1}^t n_j
=
\ell+\sum_{j=1}^t a_{\lambda_j}.
\]
Thus $\sum_j a_{\lambda_j}\ge m-\ell +1$, and therefore $a\in \mathcal{A}_{m, \ell}$.
\end{proof}

Similarly, if every sector of a filling of $\mathbb{RP}^m$ has geometric dimension at most $n$, then
 $\ell(n+1)\ge m+1$. Indeed, for $m\ge 1$, the class $x^{n+1}$ restricts to zero on every sector. Applying Theorem~\ref{th:cup_product_obstruction} to $\ell$ copies of this class, with every $p_j=1$, gives $x^{\ell(n+1)}=0$. Hence $m+1\le \ell(n+1)$. For $n=0$, this recovers the classical bound $\operatorname{cat}(\R P^m)\ge m+1$. Proposition~\ref{prop:rpm} immediately implies the following.

\begin{corollary}[Necessary dimension bounds]\label{c:necessity} Every nondecreasing universal upper bound $a=(a_1, \ldots, a_\ell)$ in dimension $m\ge 1$ belongs to $\mathcal{A}_{m, \ell}$. 
\end{corollary}

\section{Canonical fillings of triangulated manifolds}\label{s:canonical}

\subsection{Construction of canonical fillings}

Let $M$ be a compact PL manifold of dimension $m$, and let $K$ be a finite triangulation of $M$ compatible with its PL structure. We write $K'$ and $K''$ for the first and second barycentric subdivisions of $K$. We use $\sigma, \tau$ and $\eta$ to denote simplices in $K, K'$ and $K''$, respectively. 
The barycenter of a simplex $\alpha$ is denoted by $\hat\alpha$. 

 Recall that the \emph{order complex} $\Delta(\mathcal{F})$ of a poset $\mathcal{F}$ is the simplicial complex whose vertices are the elements of $\mathcal{F}$ and whose $k$-simplices are strict chains $p_0<\cdots < p_k$ in $\mathcal{F}$. If $\mathcal F$ is a set of nonempty
simplices of $K$ ordered by inclusion, then $\Delta(\mathcal F)$ is a
subcomplex of $K'$, where the simplex corresponding to a chain
$\sigma_0\subsetneq\cdots\subsetneq\sigma_k$ is identified with
$[\hat\sigma_0,\ldots,\hat\sigma_k]$.

Choose a partition $\mathcal{P}=(P_1, \ldots, P_\ell)$
of $\{0,\dots,m\}$ into nonempty disjoint blocks. A vertex $\hat{\sigma}$ of $K'$ has \emph{color} $i$ if $\dim \sigma\in P_i$. 
Let $Y_i$ be the full subcomplex of $K'$ spanned by the vertices of color $i$, and let $Y_i'$ be its first barycentric
subdivision.  
 Let
$N_i$ denote the closed star of $Y_i'$ in $K''$. We will show that $\{N_i\}$ is a filling of $M$ (a relative filling if
$\partial M\neq\emptyset$), and call it the \emph{canonical filling associated with $K$ and $\mathcal{P}$}.

For a nonempty subset $J\subseteq\{1,\dots,\ell\}$, let $\mathcal{F}_J$ denote the poset of nonempty simplices of $K'$ containing at least one vertex of each color in $J$, ordered by inclusion.  Vertices of other colors are also allowed.  
Let $\mathcal{Q}_J\subseteq \mathcal{F}_J$ denote the subposet of simplices of $K'$ whose vertices have exactly the colors in $J$, and put $Z_J=\Delta(\mathcal{Q}_J)$. Finally, write $N_J=\cap_{i\in J}N_i$ and $P_J=\cup_{i\in J}P_i$.

\begin{figure}[h!]
\centering
\begin{tikzpicture}[scale=.8]
\begin{scope}
\base
\fill[gray!12] (a)--(ab)--(t)--cycle;
\mesh
\foreach \p in {a,b,c} \fill[one] (\p) circle (3pt);
\foreach \p in {ab,ac,bc} \fill[two] (\p) circle (3pt);
\fill[three] (t) circle (3pt);
\node[below left] at (a) {$\widehat v$};
\node[below] at (ab) {$\widehat e$};
\node[right] at (t) {$\widehat \sigma$};
\node[above left=8pt, font=\Large] at (ac) {$\sigma$};
\end{scope}
\begin{scope}[xshift=6cm]
\base
\fill[one!20] (a)--(ab)--(t)--(ac)--cycle;
\fill[two!20] (b)--(bc)--(t)--(ab)--cycle;
\fill[three!20] (c)--(ac)--(t)--(bc)--cycle;
\mesh
\draw[black,line width=1pt] (ab)--(t)--(ac) (t)--(bc);
\node[below left,text=one] at (a) {$a$};
\node[below right,text=two] at (b) {$b$};
\node[above,text=three] at (c) {$c$};
\node at (.85,.45) {$N_1$};
\node at (3.95,.45) {$N_2$};
\node at (2.4,2.6) {$N_3$};
\node[above left=8pt, font=\Large] at (ac) {$\tau$};
\end{scope}
\begin{scope}[xshift=12cm]
\base\mesh
\draw[two,line width=2pt] (ab)--(t);
\fill[two] (ab) circle (3pt);
\fill[black] (t) circle (3pt);
\node[below left] at (a) {$a$};
\node[below right] at (b) {$b$};
\node[above] at (c) {$c$};
\node[below,align=center] at (ab) {$\widehat{ab}$};
\node[right] at (t) {$\widehat\tau$};
\node[right,text=two] at (2.45,.65) {$N_{\{1,2\}}\cap\tau$};
\node[above left=8pt, font=\Large] at (ac) {$\tau$};
\end{scope}
\end{tikzpicture}
\caption{Left: a simplex $\sigma$ in $K$ with its subdivision in $K'$, colored according to the partition $P_i=\{i-1\}$, $i=1, 2, 3$. Middle and right: the simplex $\tau=[a,b,c]=[\hat{v}, \hat{e}, \hat{\sigma}]$ in $K'$ with its barycentric subdivision in $K''$.}
\end{figure}
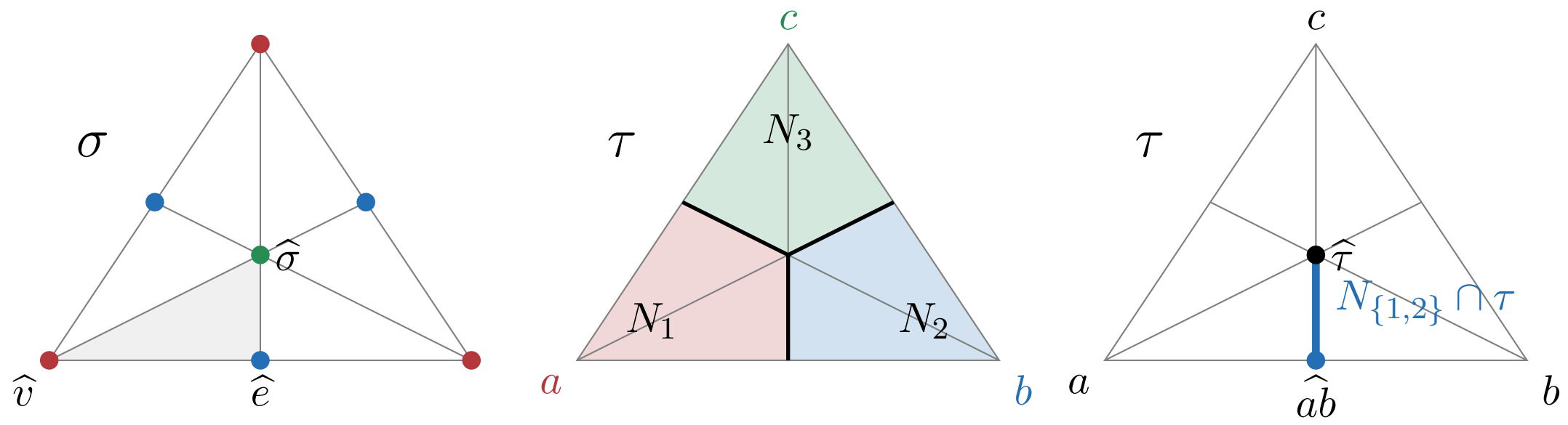

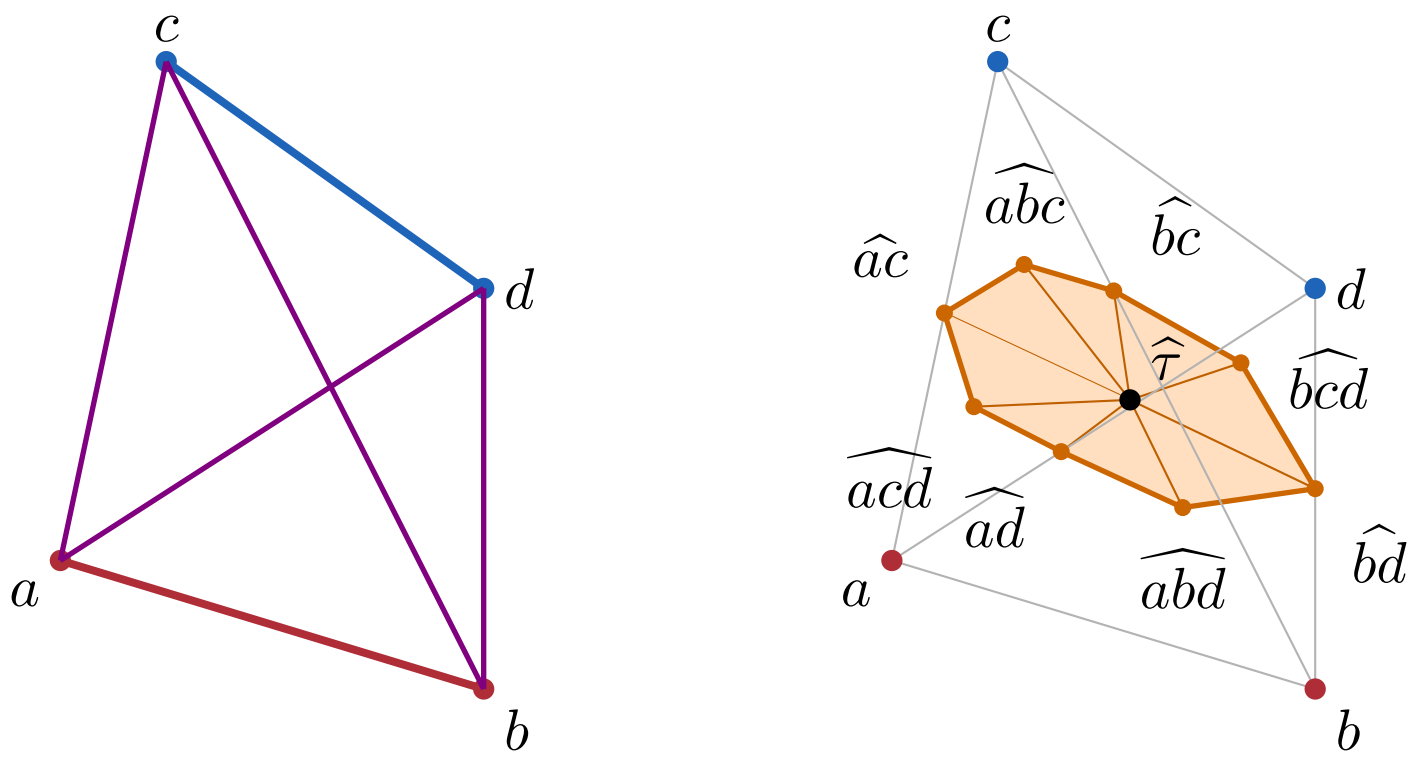
\begin{figure}[h!]\centering
\begin{tikzpicture}[scale=1]
\coords\framepic
\draw[redgraph,line width=1.5pt] (a)--(b);
\draw[bluegraph,line width=1.5pt] (c)--(d);
\foreach \u/\v in {a/c,a/d,b/c,b/d}{\draw[violet,line width=1pt] (\u)--(\v);}
\begin{scope}[xshift=5.5cm]
\coords

\coordinate (ac) at ($(c)!1.7cm!(a)$);
\coordinate (bc) at ($(c)!1.7cm!(b)$);
\coordinate (ad) at ($(a)!0.4!(d)$);
\coordinate (bd) at ($(b)!0.5!(d)$);

\coordinate (q) at ($(ac)!0.5!(bc)$);
\coordinate (abc) at ($(q)!0.25cm!90:(bc)$);

\coordinate (q) at ($(bc)!0.5!(bd)$);
\coordinate (bcd) at ($(q)!0.25cm!90:(bd)$);

\coordinate (q) at ($(bd)!0.5!(ad)$);
\coordinate (abd) at ($(q)!0.25cm!90:(ad)$);

\coordinate (q) at ($(ad)!0.5!(ac)$);
\coordinate (acd) at ($(q)!0.25cm!90:(ac)$);

\foreach \u/\v in {ac/abc,abc/bc,bc/bcd,bcd/bd,bd/abd,abd/ad,ad/acd,acd/ac}{
  \fill[orange!25] (o)--(\u)--(\v)--cycle;
  \draw[orange!75!black,thin] (o)--(\u);
  \draw[orange!80!black,line width=1pt] (\u)--(\v);
}

\framepic
\foreach \p in {ac,ad,bc,bd,abc,abd,acd,bcd}{
  \fill[orange!80!black] (\p) circle (1.6pt);
}
\fill[black] (o) circle (2pt);
\node[above right] at (o) {$\widehat{\tau}$};
\node[above=4pt]       at (abc) {$\widehat{abc}$};
\node[right=5pt,yshift=-3pt] at (bcd) {$\widehat{bcd}$};
\node[below=4pt]       at (abd) {$\widehat{abd}$};
\node[below left=4pt]  at (acd) {$\widehat{acd}$};

\node[above left=3pt]  at (ac) {$\widehat{ac}$};
\node[above right=3pt] at (bc) {$\widehat{bc}$};
\node[below left=3pt]  at (ad) {$\widehat{ad}$};
\node[below right=3pt] at (bd) {$\widehat{bd}$};


\end{scope}

\end{tikzpicture}
\caption{Left: the red and blue edges are $Y_1\cap\tau$ and $Y_2\cap\tau$, respectively. Right: a schematic drawing of the eight-triangle disc $N_{\{1,2\}}\cap \tau$, with barycenters repositioned for clarity. }
\end{figure}

Lemma~\ref{prop:canonical} and Proposition~\ref{prop:filling} show that every compatible triangulation together with a partition $\mathcal{P}$ determines a PL filling of $M$ when $M$ is closed, and a relative PL filling when $M$ has boundary.

\begin{lemma}\label{prop:canonical} The complexes $N_1,\ldots,N_\ell$ cover $M$ and have pairwise disjoint interiors in $M$. For every nonempty $J\subseteq\{1,\ldots,\ell\}$, the finite subcomplexes $N_J$ and $Z_J$ of $K''$ satisfy
\[
N_J=\Delta(\mathcal{F}_J), \qquad
\dim N_J=m-|J|+1,
\qquad
N_J\searrow Z_J,
\qquad
\dim Z_J\le |P_J|-|J|.
\]
\end{lemma}
\begin{proof}
A simplex of $K''$ determined by a chain $\tau_0\subsetneq \cdots \subsetneq \tau_k$ of nonempty simplices of $K'$ belongs to $N_i$ precisely when $\tau_0$ has a vertex of color $i$. Consequently, the sectors cover $M$, and $N_J=\Delta(\mathcal{F}_J)$. Every top dimensional simplex of $K''$ belongs to exactly one sector, since the smallest member $\tau_0$ of its defining chain is a vertex. Thus the sector interiors are pairwise disjoint.

 Every simplex in $\mathcal F_J$ has at least $|J|$ vertices and at most $m+1$ vertices. Consequently, a strict chain in $\mathcal F_J$ determines a simplex of dimension at most $m-|J|+1$.
 To show equality, we note that every maximal simplex $\tau$ of $K'$ has a vertex of each color. Choose a face consisting of one vertex of each color in $J$, and successively add the remaining vertices of $\tau$. The resulting chain lies in $\mathcal{F}_J$ and determines a simplex of $N_J$ of dimension $m-|J|+1$. 
Similarly, every simplex in $\mathcal Q_J$ has at least $|J|$ vertices and at most $|P_J|$ vertices. Indeed, its vertices correspond to a chain of simplices of $K$ whose distinct dimensions belong to $P_J$. The same counting argument gives $\dim Z_J\le |P_J|-|J|$.

Define 
 $\rho_J\co \mathcal{F}_J\to \mathcal{F}_J$ by deleting from a simplex all vertices whose colors are not in $J$. Since every simplex in $\mathcal{F}_J$ contains a vertex of every color in $J$, the resulting face is nonempty and belongs to $\mathcal{Q}_J$. 
Furthermore, the map $\rho_J$ is order preserving and satisfies $\rho_J^2=\rho_J$ and $\rho_J(\tau)\subseteq \tau$ for each $\tau$. Thus, $\rho_J$ is a  
descending closure operator. By Proposition~\ref{prop:collapse}, the map $\rho_J$ produces $N_J\searrow Z_J$. 
\end{proof}

We defer to the appendix the proof of the following proposition. 

\begin{proposition}\label{prop:filling}
Under the assumptions of Lemma~\ref{prop:canonical}, the ordered collection $(N_1, \ldots, N_\ell)$ is a PL filling if $\partial M=\emptyset$ and a relative PL filling otherwise.    
\end{proposition}

\begin{theorem}\label{th:NJ}
    Under the assumptions of Lemma~\ref{prop:canonical}, we have $\operatorname{gdim}N_J\le |P_J|-|J|$ for every nonempty $J\subseteq \{1, \ldots, \ell\}$. 
\end{theorem}
\begin{proof}
    By Lemma~\ref{prop:canonical} and Proposition~\ref{prop:filling}, $N_J$ is a compact PL manifold with possibly empty boundary such that $N_J\searrow Z_J$ with $\dim Z_J\le |P_J|-|J|$. Therefore, Theorem~\ref{th:NJ} follows from Lemma~\ref{l:forNJ}. 
\end{proof}

\subsection{Consequences of the canonical construction}

The following results describe bounds obtained directly from the canonical construction and Corollary~\ref{c:necessity}.

\begin{corollary}[Prescribed sector bound]\label{prop:conseq}
Let $M$ be a compact, possibly nonorientable, PL manifold of dimension $m$ equipped with a compatible triangulation. Let $1\le \ell\le m+1$, and let $n_1\ge \cdots \ge n_\ell$ be nonnegative integers with $\sum_{i=1}^{\ell}(n_i+1)\ge m+1$. 
Then $M$ admits a canonical PL filling (relative if $\partial M\neq \emptyset$) 
such that  $\operatorname{gdim}N_J\le\sum_{i\in J}n_i$ for every nonempty $J\subseteq\{1,\ldots,\ell\}$. In particular, $\operatorname{gdim}N_i\le n_i$ for every $i$.
\end{corollary}
\begin{proof}
Fix a partition
$\{0,\ldots,m\}=P_1\sqcup\cdots\sqcup P_\ell$ into nonempty blocks with $|P_i|\le n_i+1$. Let $(N_1,\ldots,N_\ell)$ be the associated canonical filling. Then, by Theorem~\ref{th:NJ},
$\operatorname{gdim}N_J  \le |P_J|-|J| =\sum_{i\in J}(|P_i|-1)
\le\sum_{i\in J}n_i$. 
\end{proof}

\begin{corollary}[Sharp universal sector bound] Let $1\le\ell\le m+1$, and let $n\ge 0$ be an integer. Then every nonempty closed PL $m$-manifold admits a filling into $\ell$ sectors, each of geometric dimension at most $n$, if and only if $\ell(n+1)\ge m+1$. Consequently, the least universal common bound for the geometric dimension of the sectors of closed PL $m$-manifolds is $n\ge \left\lceil\frac{m+1}{\ell}\right\rceil-1$.
\end{corollary}
\begin{proof}
    Sufficiency follows from Corollary~\ref{prop:conseq} by taking $n_1=\cdots=n_\ell=n$. Necessity follows from Corollary~\ref{c:necessity}.
\end{proof}

\begin{corollary}[Balanced intersection bounds]\label{th:canonical-depth-r}
Let $M$ be a compact PL manifold of dimension $m$, and let
$1\le \ell\le m+1$. Write $m+1-\ell=q\ell+e$, where 
$0\le e<\ell$. Then $M$ admits a canonical PL filling (relative if $\partial M\neq\emptyset$)
such that 
$\operatorname{gdim}N_J
\le
jq+\min\{j,e\}
$ for every nonempty $J\subseteq\{1,\ldots,\ell\}$, where $j=|J|$. 
\end{corollary}
\begin{proof}
Choose a partition
$\mathcal P=(P_1,\ldots,P_\ell)$
of $\{0,\ldots,m\}$ into nonempty blocks such that 
\[
|P_i|-1=
\begin{cases}
q+1,& 1\le i\le e,\\
q,& e<i\le\ell.
\end{cases}
\]
Choose a compatible triangulation of $M$, and let
$M=N_1\cup \cdots\cup N_\ell$
be the canonical filling associated with the partition $\mathcal P$.
For every nonempty
$J\subseteq\{1,\ldots,\ell\}$, put $j=|J|$. Then, by Theorem~\ref{th:NJ}, 
\[
\operatorname{gdim}N_J
\le |P_J|-j
=\sum_{i\in J}(|P_i|-1)
=jq+|J\cap\{1,\ldots,e\}|
\le jq+\min\{j,e\}.
\]
\end{proof}

\section{Compatible neighborhoods and collars}\label{s:compatible_neighborhoods}

In this section we construct product coordinates near a terminal handle of an intersection $N_J$. These coordinates will be used to describe how a handle trade changes the sectors and their intersections. 

We first recall two results of Rourke and Sanderson concerning the existence and uniqueness of block bundle structures on regular neighborhoods. Throughout this section, all maps are PL and all embeddings of manifolds are locally flat. An embedded manifold $M\subset V$ is called \emph{proper} if $M\cap \partial V=\partial M$. We write $D^k$ for the standard PL ball $[-1,1]^k$ with $D^0$ a point. Cells in a PL cell complex are understood to be closed PL balls. A \emph{complex} is a finite collection $K$ of cells which cover a polyhedron $X$ such that the boundary of each cell, and the intersection of two cells are unions of cells, and the interiors of different cells are disjoint. We then write $|K|=X$. Every finite PL complex $K$ admits a subdivision given by a triangulation of $|K|$.

A \emph{block bundle} $\xi$ of rank $q$ over a PL cell complex $K$ consists of a polyhedron $E=E(\xi)$ with $|K|\subset E$, together with a family of balls $\{\beta_\sigma\}_{\sigma\in K}$, called \emph{blocks}, satisfying the following conditions: 
\begin{itemize}
\item For each $k$-cell $\sigma\in K$, the block $\beta_\sigma\subset E$ is a PL ball of dimension $k+q$, and $(\beta_\sigma,\sigma)$ is PL homeomorphic to $(D^k\times D^q, D^k\times \{0\})$. 
\item The blocks $\beta_\sigma$ cover $E$ and have pairwise disjoint interiors. 
\item For any $\sigma, \tau\in K$,  the intersection $\beta_\sigma\cap \beta_\tau$ is the union of the blocks over cells in $\sigma\cap \tau$. 
\end{itemize}
 An isomorphism $\xi\cong \eta$ of block bundles over $K$ is a PL homeomorphism $E(\xi)\to E(\eta)$ that fixes $|K|$ pointwise and takes each block $\beta_\sigma(\xi)$ to the corresponding block $\beta_\sigma(\eta)$. 
For a subcomplex $L\subset K$, the restriction $\xi|L$ is the block bundle over $L$ that consists of the blocks $\beta_\sigma$ with $\sigma\in L$. The bundle $\xi$ is called \emph{trivial} if it admits an isomorphism to the product block bundle with total space $|K|\times D^q$ and blocks $\sigma\times D^q$. No projection onto $|K|$ is part of the definition of a block bundle. 

Rourke and Sanderson established the existence and uniqueness results for block bundles.

\begin{theorem}[{\cite[Theorem 4.3]{RSBlock}}] 
Let $Q$ be a compact PL manifold of dimension $d+q$ and let $B\subset Q$ be a properly embedded compact PL submanifold of dimension $d$ such that $Q\searrow B$.
Let $K$ be a PL cell complex and $L\subset K$ a subcomplex such that $(B, \partial B)=(|K|, |L|)$. Then there is a block bundle $\xi$ of rank $q$ over $K$ with $E(\xi)=Q$. Furthermore, given a block bundle $\eta$ of rank $q$ over $L$ whose total space lies in $\partial Q$, the block bundle $\xi$ may be chosen so that $\xi|L=\eta$.   
\end{theorem}

\begin{theorem}[{\cite[Theorem 4.4]{RSBlock}}]\label{th:RS4.4} 
\begin{enumerate}
   \item[(a)] Let $V$ be a PL manifold of dimension $d+q$, and let $B\subset V$ be a properly embedded compact PL submanifold of dimension $d$. Let $K$ be a PL cell complex and $L\subset K$ be a subcomplexes such that $(B, \partial B)=(|K|, |L|)$. Let $\xi$ and $\eta$ be block bundles of rank $q$ over $K$ whose total spaces are embedded in $V$, with both embeddings restricting to the given inclusion $B\subset V$. Then there is a PL isotopy of the inclusion $E(\xi)\subset V$ fixed on $B$  whose terminal map is an isomorphism from $\xi$ to $\eta$.  
   \item[(b)]
   If, in addition, $E(\xi)\cap \partial V=E(\xi|L)$ and $E(\eta)\cap \partial V=E(\eta|L)$,  then the isotopy of (a) can be chosen to be ambient. 
   \item[(c)] Under the hypotheses of (b), suppose that $\xi|L=\eta|L$. Then the isotopy of (b) can be chosen to fix $\partial V\cup B$ pointwise.
   \end{enumerate}
\end{theorem}

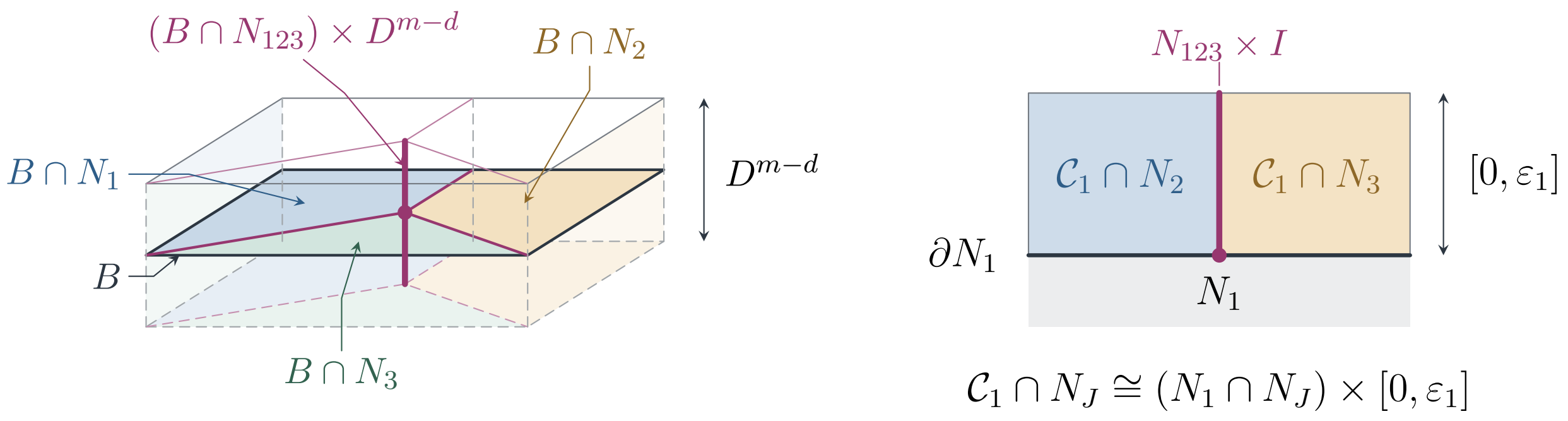
\begin{figure}
    \centering
\begin{tikzpicture}[line join=round,line cap=round,font=\large,>=stealth]
\definecolor{sectorblue}{RGB}{58,116,170}
\definecolor{sectorgold}{RGB}{218,160,61}
\definecolor{sectorgreen}{RGB}{77,151,121}
\definecolor{seam}{RGB}{151,55,111}
\definecolor{ink}{RGB}{44,54,65}
\begin{scope}[x={(1cm,0cm)},y={(0.51cm,0.32cm)},z={(0cm,1cm)}]
  \fill[sectorblue!7] (0,0,-.75)--(0,2.8,-.75)--(0,2.8,.75)--(0,0,.75)--cycle;
  \fill[sectorgold!7] (4,0,-.75)--(4,2.8,-.75)--(4,2.8,.75)--(4,0,.75)--cycle;
  \fill[sectorgreen!7] (0,0,-.75)--(4,0,-.75)--(4,0,.75)--(0,0,.75)--cycle;
  \fill[sectorblue!12] (0,0,-.75)--(2,1.4,-.75)--(2,2.8,-.75)--(0,2.8,-.75)--cycle;
  \fill[sectorgold!14] (4,0,-.75)--(4,2.8,-.75)--(2,2.8,-.75)--(2,1.4,-.75)--cycle;
  \fill[sectorgreen!12] (0,0,-.75)--(4,0,-.75)--(2,1.4,-.75)--cycle;
  \draw[ink!45,densely dashed] (0,0,-.75)--(4,0,-.75)--(4,2.8,-.75)--(0,2.8,-.75)--cycle;
  \draw[seam!55,densely dashed] (0,0,-.75)--(2,1.4,-.75)--(4,0,-.75);
  \draw[seam!55,densely dashed] (2,1.4,-.75)--(2,2.8,-.75);
  \fill[sectorblue!28] (0,0,0)--(2,1.4,0)--(2,2.8,0)--(0,2.8,0)--cycle;
  \fill[sectorgold!32] (4,0,0)--(4,2.8,0)--(2,2.8,0)--(2,1.4,0)--cycle;
  \fill[sectorgreen!25] (0,0,0)--(4,0,0)--(2,1.4,0)--cycle;
  \draw[ink,thick] (0,0,0)--(4,0,0)--(4,2.8,0)--(0,2.8,0)--cycle;
  \draw[seam,thick] (0,0,0)--(2,1.4,0)--(4,0,0);
  \draw[seam,thick] (2,1.4,0)--(2,2.8,0);
  \draw[ink!65] (0,0,.75)--(4,0,.75)--(4,2.8,.75)--(0,2.8,.75)--cycle;
  \draw[seam!65] (0,0,.75)--(2,1.4,.75)--(4,0,.75);
  \draw[seam!65] (2,1.4,.75)--(2,2.8,.75);
  \foreach \x/\y in {0/0,4/0,0/2.8,4/2.8,2/2.8}{
    \draw[ink!45,densely dashed] (\x,\y,-.75)--(\x,\y,.75);
  }
  \draw[seam,line width=1.8pt] (2,1.4,-.75)--(2,1.4,.75);
  \fill[seam] (2,1.4,0) circle[radius=2.2pt];
  \coordinate (targetQone)   at (.8,1.7,0);
  \coordinate (targetQtwo)   at (3.1,1.7,0);
  \coordinate (targetQthree) at (2,.45,0);
  \coordinate (targetB)      at (.35,0,0);
  \coordinate (targetSeam)   at (2,1.4,.48);

  \draw[<->,ink,thin]
    (4.42,2.8,-.75)--(4.42,2.8,.75);
  \node[anchor=west,font=\normalsize,inner sep=2pt]
    at (4.56,2.8,0) {$D^{m-d}$};

  \begin{scope}[
    x={(1cm,0cm)}, y={(0cm,1cm)},
    every node/.style={font=\normalsize,inner sep=2pt},
    every path/.style={line width=.4pt}
  ]
    \node[anchor=east,text=sectorblue!80!black]
      (labelQone) at (-.18,.85) {$B\cap N_1$};
    \draw[->,sectorblue!80!black]
      (labelQone.east)--(.45,.85)--(targetQone);

    \node[anchor=south,text=sectorgold!65!black]
      (labelQtwo) at (4.65,1.98) {$B\cap N_2$};
    \draw[->,sectorgold!65!black]
      (labelQtwo.south)--(4.65,1.73)--(targetQtwo);

    \node[anchor=north,text=sectorgreen!65!black]
      (labelQthree) at (2.05,-1) {$B\cap N_3$};
    \draw[->,sectorgreen!65!black]
      (labelQthree.north)--(2.05,-.45)--(targetQthree);

    \node[anchor=east,text=ink]
      (labelB) at (-.18,-.22) {$B$};
    \draw[->,ink]
      (labelB.east)--(.08,-.22)--(targetB);

    \node[anchor=south,text=seam]
      (labelSeam) at (1.65,2.05)
      {$(B\cap N_{123})\times D^{m-d}$};
    \draw[->,seam]
      (labelSeam.south)--(2.08,1.70)--(targetSeam);
  \end{scope}
\end{scope}
\begin{scope}[shift={(9.25cm,0cm)}]
   \fill[ink!9] (0,-.75) rectangle (4,0);
  \fill[sectorblue!25] (0,0) rectangle (2,1.7);
  \fill[sectorgold!30] (2,0) rectangle (4,1.7);
  \draw[ink!65] (0,0)--(0,1.7)--(4,1.7)--(4,0);
  \draw[ink,line width=1.2pt] (0,0)--(4,0);
  \draw[seam,line width=1.8pt] (2,0)--(2,1.7);
  \fill[seam] (2,0) circle[radius=2.2pt];
  \node at (2,-.4) {$N_1$};
  \node[anchor=east] at (-.18,0) {$\partial N_1$};
  \node[text=sectorblue!80!black] at (.97,.83) {$\mathcal{C}_1\cap N_2$};
  \node[text=sectorgold!65!black] at (3.03,.83) {$\mathcal{C}_1\cap N_3$};
   \node[seam,font=\normalsize] at (2,2.2) {$N_{123}\times I$};
  \draw[seam] (2,2.02)--(2,1.78);
  \draw[<->,ink] (4.35,0)--(4.35,1.7);
  \node[anchor=west] at (4.47,.85) {$[0, \varepsilon_1]$};
  \node at (2,-1.42) {$\mathcal{C}_1\cap N_J\cong(N_1\cap N_J)\times [0, \varepsilon_1]$};
\end{scope}
\end{tikzpicture}
\caption{Left: a compatible product neighborhood of $B$ in which the filling of $B$ extends along the normal discs. Right: a compatible exterior collar of $N_1$, in which the induced filling of $\partial N_1$ extends along the collar. }
    \label{fig:compatibility}
\end{figure}

We now apply these ideas to neighborhoods compatible with a filling. Let $M$ be a closed connected PL manifold of dimension $m$ with a filling $M=N_1\cup \cdots \cup N_\ell$, where $\ell\ge 2$. Let $B$ be a closed PL submanifold of $M$ of dimension $d$ transverse to each stratum of the filling $\{N_i\}$. Suppose that $B$ admits a regular neighborhood that is PL homeomorphic to $B\times D^{m-d}$ under a homeomorphism that takes $B$ to $B\times \{0\}$. Then there is a regular neighborhood $V$ of $B$ and a PL homeomorphism $\psi\co V\cong B\times D^{m-d}$ such that $\psi(b)=(b, 0)$ for each $b\in B$ and 
\[
\psi(V\cap N_J) = (B\cap N_J)\times D^{m-d}
\]
for every nonempty $J\subseteq \{1, \ldots, \ell\}$, see \cite[Lemma 3.6]{Si}. We say that such a product neighborhood is \emph{compatible} with the filling. 

Let $N_\beta$ be one of the sectors of the filling. Then an \emph{exterior collar} of $N_\beta$ is an embedding
\[
c_\beta:\partial N_\beta\times[0,\varepsilon_\beta]
\longrightarrow \overline{M\setminus N_\beta},
\qquad c_\beta(z,0)=z,
\]
such that $c_\beta(z,t)\notin N_\beta$ for $t>0$.
We say that an exterior collar $c_\beta$ of $N_\beta$ is \emph{compatible} with the filling provided that for every $z\in \partial N_\beta$, every $t\in [0, \varepsilon_\beta]$, and every $i\neq \beta$, we have
$c_\beta(z, t)\in N_i$ if and only if $z\in N_i$. 
Such collars exist by the construction in the proof of \cite[Lemma 4.2]{Si} for $\ell\ge 3$, and by the PL collaring theorem for $\ell=2$. In particular, if $\mathcal{C}_\beta$ is the image of $c_\beta$, then 
\[
  c_\beta^{-1}(\mathcal{C}_\beta\cap N_J)=(\partial N_\beta\cap N_J)\times [0, \varepsilon_\beta]
\]
for every nonempty $J\subseteq \{1, \ldots, \ell\}\setminus \{\beta\}$.

Fix a nonempty proper subset $J\subsetneq \{1, \ldots, \ell\}$ with $N_J\neq\emptyset$, and put $j=|J|$.
After relabeling, we may assume that $J=\{1,\ldots,j\}$. The sectors indexed by $J$ are called
\emph{donor sectors}; the others are called \emph{receiving sectors}.
Choose a handle decomposition of the manifold $N_J$ with a terminal handle $H$ so that 
\[
 N_J=A\cup H,\qquad H=D^p\times D^h,\qquad p+h=m-j+1,
\]
where $A$ is the union of the other handles in $N_J$, the handle $H$ is attached along $D^p\times\partial D^h$, 
$\partial D^p\times D^h$ is the free boundary of $H$, and $h, p\ge 0$. 
Since $H$ is a terminal handle, the receiving sectors meet $H$ only along its free boundary and cover it.  Let
$\Delta=D^p\times\{0\}\cong D^p$ be the cocore, and $S=\partial\Delta\cong \partial D^p$. After a small isotopy, we may assume that $S$ is transverse to the strata of the
induced filling of $\partial N_J$. For $i>j$, put $E_i=S\cap N_i$.
Then the sets $E_i$, including empty labeled pieces, form the induced filling of 
$S$. If $p=0$, all these sets are empty.

 Let $\widetilde D^p\supset D^p$ and $\widetilde D^h\supset D^h$ be open PL balls, and let $S\times [0, \varepsilon)\subset \widetilde D^p\setminus \Int D^p$ be an outward collar of $S$ for a small $\varepsilon>0$. 
Put 
\[C=[0,1]^{j-1},\qquad\widetilde C=(-\varepsilon,1+\varepsilon)^{j-1},\qquad \widetilde E_i=E_i\times[0,\varepsilon)\subset\widetilde D^p, \qquad i>j.\]
We write points of $\widetilde D^p\times \widetilde D^h\times \widetilde C$ as $(x, y, t)$ with $t=(t_2, \ldots, t_j)$. 
All zero-dimensional factors are understood to be points. We write $t\ge 0$ when every coordinate of $t=(t_2, \ldots, t_j)$ is nonnegative.

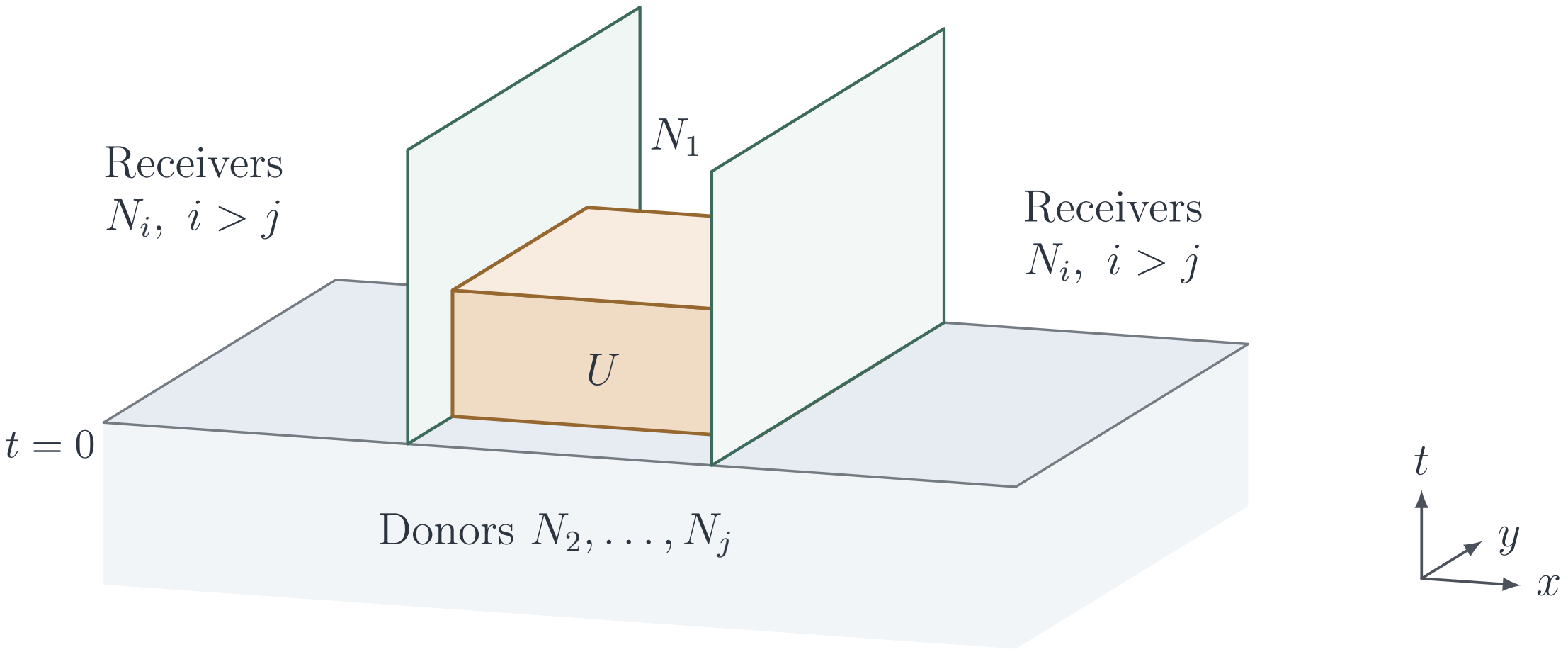
\begin{figure}
    \centering
    \begin{tikzpicture}[
  font=\large, line cap=round, line join=round,
  x={(1.42cm,-.10cm)},
  y={(.70cm,.43cm)},
  z={(0cm,1.12cm)}]
  \definecolor{ink}{RGB}{43,52,63}
  \definecolor{donors}{RGB}{105,143,181}
  \definecolor{receivers}{RGB}{86,151,132}
  \definecolor{transfer}{RGB}{208,143,65}
  \tikzset{
    edge/.style={draw=ink!65,line width=.65pt},
    hidden/.style={draw=ink!48,line width=.55pt,densely dashed},
    wall/.style={draw=receivers!70!black,line width=.8pt},
    boxedge/.style={draw=transfer!72!black,line width=.9pt},
    label/.style={text=ink,inner sep=2pt,align=center},
    axis/.style={->,draw=ink!85,line width=.7pt,>=latex}
  }

  \fill[donors!9]
    (-3,-1.5,-1.35)--(3,-1.5,-1.35)--(3,-1.5,0)--(-3,-1.5,0)--cycle;
  \fill[donors!9]
    (3,-1.5,-1.35)--(3,1.6,-1.35)--(3,1.6,0)--(3,-1.5,0)--cycle;

  \filldraw[edge,fill=donors!17]
    (-3,-1.5,0)--(3,-1.5,0)--(3,1.6,0)--(-3,1.6,0)--cycle;

  \fill[receivers!9]
    (-1,-1.5,0)--(-1,1.6,0)--(-1,1.6,2.45)--(-1,-1.5,2.45)--cycle;
  \draw[wall]
    (-1,-1.5,0)--(-1,1.6,0)--(-1,1.6,2.45)--(-1,-1.5,2.45)--cycle;

  \draw[hidden] (-1,-.9,0)--(-1,.9,0)--(1,.9,0);
  \draw[hidden] (-1,.9,0)--(-1,.9,1.05);
  \filldraw[boxedge,fill=transfer!25]
    (-1,-.9,0)--(-1,.9,0)--(-1,.9,1.05)--(-1,-.9,1.05)--cycle;
  \filldraw[boxedge,fill=transfer!31]
    (-1,-.9,0)--(1,-.9,0)--(1,-.9,1.05)--(-1,-.9,1.05)--cycle;
  \filldraw[boxedge,fill=transfer!19]
    (1,-.9,0)--(1,.9,0)--(1,.9,1.05)--(1,-.9,1.05)--cycle;
  \filldraw[boxedge,fill=transfer!17]
    (-1,-.9,1.05)--(1,-.9,1.05)--(1,.9,1.05)--(-1,.9,1.05)--cycle;

  \fill[receivers!7]
    (1,-1.5,0)--(1,1.6,0)--(1,1.6,2.45)--(1,-1.5,2.45)--cycle;
  \draw[wall]
    (1,-1.5,0)--(1,1.6,0)--(1,1.6,2.45)--(1,-1.5,2.45)--cycle;

  \node[label] at (0,-.92,.49) {$U$};
  \node[label] at (0,.05,2.05) {$N_1$};
  \node[label] at (-2.95,-.4,1.5) {Receivers\\$N_i,\ i>j$};
  \node[label] at (2.9,0,1.5) {Receivers\\$N_i,\ i>j$};
  \node[label] at (0,-1.55,-.65) {Donors $N_2,\ldots,N_j$};
  \node[label,anchor=north east,font=\normalsize]
    at (-3,-1.5,0) {$t=0$};

  \begin{scope}[x={(1cm,0cm)},y={(0cm,1cm)}]
    \coordinate (O) at (7.0,-1.8);
    \draw[axis] (O)--++(.93,-.065) node[right,text=ink] {$x$};
    \draw[axis] (O)--++(.57,.35) node[right,text=ink] {$y$};
    \draw[axis] (O)--++(0,.83) node[above,text=ink] {$t$};
  \end{scope}
\end{tikzpicture}
    \caption{Sectors near the handle $H$.}
    \label{fig:sectors_near_the_handle}
\end{figure}

\begin{proposition}\label{prop:compatible-handle-neighborhood}
After an ambient isotopy of $N_J$ fixing $\Delta$ pointwise, the factors $\widetilde D^p, \widetilde D^h$ and $\widetilde C$ can be chosen so that there is a PL homeomorphism 
\[
 \widetilde\Phi:\widetilde D^p\times\widetilde D^h\times\widetilde C
       \longrightarrow\widetilde U
\]
onto an open subset $\tilde U$ of $M$ with $\tilde\Phi(D^p\times D^h\times \{0\})=H$ and
\begin{align*}
 \widetilde\Phi^{-1}(N_1)
   &=\{(x,y,t):x\in D^p,\ t_k\ge0\text{ for all }2\le k\le j\},\\
 \widetilde\Phi^{-1}(N_\beta)
   &=\{(x,y,t):t_\beta\le0,\ t_k\ge 0\text{ for all }\beta<k\le j\},
       &&2\le\beta\le j,\\
 \widetilde\Phi^{-1}(N_i)
   &=\{(x,y,t):x\in\widetilde E_i,\ t_k\ge0\text{ for all }2\le k\le j\},
       &&i>j.
\end{align*}
Unmentioned coordinates in these formulas are unrestricted within their
respective factors. The same conclusion holds for a relative filling of $M$ that is a product on a fixed collar of $\partial M$, provided that the terminal handle is disjoint from that collar.
\end{proposition}

It follows that the restriction
$\Phi=\widetilde\Phi|D^p\times D^h\times C$
is a PL embedding into $N_1$ with image $U$, and 
\[
 \Phi^{-1}(N_i)=
 \begin{cases}
 D^p\times D^h\times\{t\in C:t_i=0\},&2\le i\le j,\\
 E_i\times D^h\times C,&i>j.
 \end{cases}
\]

\begin{figure}[p]
\centering
\begin{tikzpicture}[scale=0.9,transform shape,font=\large,line cap=round,line join=round,
 x={(1.6cm,0cm)},y={(.72cm,.50cm)},z={(0cm,1.15cm)},
 edge/.style={line width=.55pt,draw=ink!45},
 hidden/.style={densely dashed,line width=.65pt,draw=ink!55},
 boxedge/.style={line width=1pt,draw=ink!80},
 leader/.style={-{Stealth[length=1.7mm]},line width=.6pt},
 lab/.style={inner sep=2pt,align=center,text=ink}]
\definecolor{ink}{RGB}{44,55,70}
\definecolor{donorTwo}{RGB}{196,138,36}
\definecolor{donorThree}{RGB}{45,104,172}
\definecolor{receiverA}{RGB}{34,135,142}
\definecolor{receiverB}{RGB}{115,83,161}
\definecolor{core}{RGB}{190,42,65}
\definecolor{volume}{RGB}{101,151,185}
\filldraw[draw=none,fill=donorTwo!6] (-2.8,-1.65,-0.65)--(2.8,-1.65,-0.65)--(2.8,-1.65,0)--(-2.8,-1.65,0)--cycle;
\filldraw[draw=none,fill=donorTwo!8] (2.8,-1.65,-0.65)--(2.8,2.35,-0.65)--(2.8,2.35,0)--(2.8,-1.65,0)--cycle;
\filldraw[edge,fill=donorTwo,fill opacity=.13] (-2.8,-1.65,0)--(2.8,-1.65,0)--(2.8,2.35,0)--(-2.8,2.35,0)--cycle;
\filldraw[draw=receiverA!75!black,line width=.8pt,fill=receiverA,fill opacity=.06] (-1.3,-1.65,0)--(-1.3,2.35,0)--(-1.3,2.35,2.35)--(-1.3,-1.65,2.35)--cycle;
\filldraw[draw=receiverB!75!black,line width=.8pt,fill=receiverB,fill opacity=.06] (1.3,-1.65,0)--(1.3,2.35,0)--(1.3,2.35,2.35)--(1.3,-1.65,2.35)--cycle;

\filldraw[draw=none,fill=donorTwo,fill opacity=.22] (-1.3,-0.9,0)--(1.3,-0.9,0)--(1.3,0.9,0)--(-1.3,0.9,0)--cycle;
\filldraw[draw=none,fill=receiverA,fill opacity=.14] (-1.3,-0.9,0)--(-1.3,0.9,0)--(-1.3,0.9,1.3)--(-1.3,-0.9,1.3)--cycle;
\filldraw[draw=none,fill=receiverB,fill opacity=.14] (1.3,-0.9,0)--(1.3,0.9,0)--(1.3,0.9,1.3)--(1.3,-0.9,1.3)--cycle;
\filldraw[draw=none,fill=volume,fill opacity=.05] (-1.3,-0.9,1.3)--(1.3,-0.9,1.3)--(1.3,0.9,1.3)--(-1.3,0.9,1.3)--cycle;
\draw[hidden] (-1.3,-0.9,0)--(-1.3,0.9,0)--(1.3,0.9,0);
\draw[hidden] (-1.3,0.9,0)--(-1.3,0.9,1.3);
\draw[boxedge] (-1.3,-0.9,0)--(1.3,-0.9,0)--(1.3,0.9,0)--(1.3,0.9,1.3)--(-1.3,0.9,1.3)--(-1.3,-0.9,1.3)--cycle;
\draw[boxedge] (-1.3,-0.9,1.3)--(1.3,-0.9,1.3)--(1.3,0.9,1.3);
\draw[boxedge] (1.3,-0.9,0)--(1.3,-0.9,1.3);
\draw[core,line width=2.3pt] (-1.3,0,0)--(1.3,0,0);
\fill[core] (-1.3,0,0) circle[radius=2pt];
\fill[core] (1.3,0,0) circle[radius=2pt];
\coordinate (lefttarget) at (-1.3,0.0,.7);
\coordinate (righttarget) at (1.3,0.0,.7);
\node[lab] at (.2,0,2.5) {$N_1$};
\node[lab,text=receiverA!85!black] at (-2.3,0,2) {$N_3$};
\node[lab,text=receiverB!85!black] at (2.9,0,2.35) {$N_4$};
\node[lab,text=donorTwo!80!black] at (-2,-1.5,-.65) {$N_2$};
\coordinate (bottomtarget) at (.7,-.5,0);
\begin{scope}[x={(1cm,0cm)},y={(0cm,1cm)}]
\node[lab,anchor=east,text=receiverA!85!black] (leftlabel) at (-3.4,1.55) {$U\cap N_3$};
\draw[leader,receiverA!85!black] (leftlabel.east)--(lefttarget);
\node[lab,anchor=west,text=receiverB!85!black] (rightlabel) at (4.05,1.55) {$U\cap N_4$};
\draw[leader,receiverB!85!black] (rightlabel.west)--(righttarget);
\node[lab,text=donorTwo!80!black] at (-4.35,-.6) {$t_2=0$};
\node[lab,text=donorTwo!80!black] (hlabel) at (0,-1.95) {$H=U\cap N_2$};
\draw[leader,donorTwo!80!black] (hlabel.north)--(bottomtarget);
\node[lab,text=core] at (0,.23) {$\Delta$};
\coordinate (O) at (6.45,-1.9);
\draw[leader,ink!80] (O)--++(.8,0) node[right] {$x$};
\draw[leader,ink!80] (O)--++(.52,.36) node[right] {$y$};
\draw[leader,ink!80] (O)--++(0,.85) node[above] {$t_2$};
\end{scope}
\end{tikzpicture}
\caption{Here $C=[0,1]$, $J=\{1,2\}$, $p=1$, and $h=1$. The whole box $U$ lies in $N_1$. The sector $N_2$ meets it in the bottom face, $t_2=0$. The sectors $N_3$ and $N_4$ are the receiving sectors.}

\medskip
\medskip
\medskip

\begin{tikzpicture}[scale=.9,transform shape,font=\large,line cap=round,line join=round,
 x={(1.6cm,0cm)},y={(.72cm,.50cm)},z={(0cm,1.15cm)},
 edge/.style={line width=.55pt,draw=ink!45},
 hidden/.style={densely dashed,line width=.65pt,draw=ink!55},
 boxedge/.style={line width=1pt,draw=ink!80},
 leader/.style={-{Stealth[length=1.7mm]},line width=.6pt},
 lab/.style={inner sep=2pt,align=center,text=ink}]
\definecolor{ink}{RGB}{44,55,70}
\definecolor{donorTwo}{RGB}{196,138,36}
\definecolor{donorThree}{RGB}{45,104,172}
\definecolor{receiverA}{RGB}{34,135,142}
\definecolor{receiverB}{RGB}{115,83,161}
\definecolor{core}{RGB}{190,42,65}
\definecolor{volume}{RGB}{101,151,185}
\filldraw[draw=none,fill=donorThree!6] (-2.8,-1.65,-0.65)--(2.8,-1.65,-0.65)--(2.8,-1.65,0)--(-2.8,-1.65,0)--cycle;
\filldraw[draw=none,fill=donorThree!8] (2.8,-1.65,-0.65)--(2.8,2.35,-0.65)--(2.8,2.35,0)--(2.8,-1.65,0)--cycle;
\filldraw[edge,fill=donorThree,fill opacity=.13] (-2.8,-1.65,0)--(2.8,-1.65,0)--(2.8,2.35,0)--(-2.8,2.35,0)--cycle;
\filldraw[draw=donorTwo!60!black,line width=.7pt,fill=donorTwo,fill opacity=.09] (-2.8,0,0)--(2.8,0,0)--(2.8,0,2.35)--(-2.8,0,2.35)--cycle;
\filldraw[draw=receiverA!75!black,line width=.8pt,fill=receiverA,fill opacity=.06] (-1.3,0,0)--(-1.3,2.35,0)--(-1.3,2.35,2.35)--(-1.3,0,2.35)--cycle;
\filldraw[draw=receiverB!75!black,line width=.8pt,fill=receiverB,fill opacity=.06] (1.3,0,0)--(1.3,2.35,0)--(1.3,2.35,2.35)--(1.3,0,2.35)--cycle;

\filldraw[draw=none,fill=donorThree,fill opacity=.22] (-1.3,0,0)--(1.3,0,0)--(1.3,1.5,0)--(-1.3,1.5,0)--cycle;
\filldraw[draw=none,fill=donorTwo,fill opacity=.16] (-1.3,0,0)--(1.3,0,0)--(1.3,0,1.3)--(-1.3,0,1.3)--cycle;
\filldraw[draw=none,fill=receiverA,fill opacity=.14] (-1.3,0,0)--(-1.3,1.5,0)--(-1.3,1.5,1.3)--(-1.3,0,1.3)--cycle;
\filldraw[draw=none,fill=receiverB,fill opacity=.14] (1.3,0,0)--(1.3,1.5,0)--(1.3,1.5,1.3)--(1.3,0,1.3)--cycle;
\filldraw[draw=none,fill=volume,fill opacity=.05] (-1.3,0,1.3)--(1.3,0,1.3)--(1.3,1.5,1.3)--(-1.3,1.5,1.3)--cycle;
\draw[hidden] (-1.3,0,0)--(-1.3,1.5,0)--(1.3,1.5,0);
\draw[hidden] (-1.3,1.5,0)--(-1.3,1.5,1.3);
\draw[boxedge] (-1.3,0,0)--(1.3,0,0)--(1.3,1.5,0)--(1.3,1.5,1.3)--(-1.3,1.5,1.3)--(-1.3,0,1.3)--cycle;
\draw[boxedge] (-1.3,0,1.3)--(1.3,0,1.3)--(1.3,1.5,1.3);
\draw[boxedge] (1.3,0,0)--(1.3,0,1.3);
\draw[core,line width=2.3pt] (-1.3,0,0)--(1.3,0,0);
\fill[core] (-1.3,0,0) circle[radius=2pt];
\fill[core] (1.3,0,0) circle[radius=2pt];
\coordinate (lefttarget) at (-1.3,0.75,.7);
\coordinate (righttarget) at (1.3,0.75,.7);
\coordinate (bottomtarget) at (0,0.7,0);
\coordinate (fronttarget) at (0,0,.7);
\begin{scope}[x={(1cm,0cm)},y={(0cm,1cm)}]
\node[lab,anchor=east,text=receiverA!85!black] (leftlabel) at (-4.75,1.55) {$U\cap N_4$};
\draw[leader,receiverA!85!black] (leftlabel.east)--(lefttarget);
\node[lab,anchor=west,text=receiverB!85!black] (rightlabel) at (4.85,1.55) {$U\cap N_5$};
\draw[leader,receiverB!85!black] (rightlabel.west)--(righttarget);
\node[lab,text=donorTwo!80!black] at (-3.25,2.3) {$t_2=0$};
\node[lab,anchor=east,text=donorTwo!80!black] (n2label) at (-3.0,-1.5) {$U\cap N_2$};
\draw[leader,donorTwo!80!black] (n2label.east)--(-2.5,-1.2)--(fronttarget);
\fill[donorTwo!80!black] (fronttarget) circle[radius=1.5pt];
\node[lab,text=donorThree!85!black] at (-3.75,-.6) {$t_3=0$};
\node[lab,anchor=west,text=donorThree!85!black] (n3label) at (1,-1.1) {$U\cap N_3$};
\draw[leader,donorThree!85!black] (n3label.west)--(0.5,-1.1)--(bottomtarget);
\fill[donorThree!85!black] (bottomtarget) circle[radius=1.5pt];
\node[lab,text=core,anchor=north] (hlabel) at (-.5,-1.5) {$H=\Delta=U\cap N_{123}$};
\draw[leader,core] (hlabel.north)--(-.5,-.4)--(-.5,0);
\coordinate (O) at (6.45,-1.9);
\draw[leader,ink!80] (O)--++(.8,0) node[right] {$x$};
\draw[leader,ink!80] (O)--++(.52,.36) node[right] {$t_2$};
\draw[leader,ink!80] (O)--++(0,.85) node[above] {$t_3$};
\end{scope}
\end{tikzpicture}
\caption{Here $C=[0,1]^2$, $J=\{1,2,3\}$, $p=1$, and $h=0$. The box $U$ lies in $N_1$. The sector $N_2$ meets it in the front face $\{t_2=0\}$, and $N_3$ meets it in the bottom face $\{t_3=0\}$. Their common red edge is $H=\Delta\subset N_{123}$. The receiving sectors $N_4$ and $N_5$ meet $U$ along two copies of $C$.}
    \label{fig:placeholder}

\medskip
\medskip
\medskip

\begin{tikzpicture}[scale=.9,transform shape,font=\large,line cap=round,line join=round,
 x={(1.6cm,0cm)},y={(.72cm,.50cm)},z={(0cm,1.15cm)},
 edge/.style={line width=.55pt,draw=ink!45},
 hidden/.style={densely dashed,line width=.65pt,draw=ink!55},
 boxedge/.style={line width=1pt,draw=ink!80},
 leader/.style={-{Stealth[length=1.7mm]},line width=.6pt},
 lab/.style={inner sep=2pt,align=center,text=ink}]
\definecolor{ink}{RGB}{44,55,70}
\definecolor{donorTwo}{RGB}{196,138,36}
\definecolor{donorThree}{RGB}{45,104,172}
\definecolor{receiverA}{RGB}{34,135,142}
\definecolor{receiverB}{RGB}{115,83,161}
\definecolor{core}{RGB}{190,42,65}
\definecolor{volume}{RGB}{101,151,185}
\filldraw[draw=none,fill=donorThree!6] (-2.8,-1.65,-0.65)--(2.8,-1.65,-0.65)--(2.8,-1.65,0)--(-2.8,-1.65,0)--cycle;
\filldraw[draw=none,fill=donorThree!8] (2.8,-1.65,-0.65)--(2.8,2.35,-0.65)--(2.8,2.35,0)--(2.8,-1.65,0)--cycle;
\filldraw[edge,fill=donorThree,fill opacity=.13] (-2.8,-1.65,0)--(2.8,-1.65,0)--(2.8,2.35,0)--(-2.8,2.35,0)--cycle;
\filldraw[draw=donorTwo!60!black,line width=.7pt,fill=donorTwo,fill opacity=.09] (-2.8,0,0)--(2.8,0,0)--(2.8,0,2.35)--(-2.8,0,2.35)--cycle;

\filldraw[draw=none,fill=donorThree,fill opacity=.22] (-1.3,0,0)--(1.3,0,0)--(1.3,1.5,0)--(-1.3,1.5,0)--cycle;
\filldraw[draw=none,fill=donorTwo,fill opacity=.16] (-1.3,0,0)--(1.3,0,0)--(1.3,0,1.3)--(-1.3,0,1.3)--cycle;
\filldraw[draw=none,fill=ink,fill opacity=.14] (-1.3,0,0)--(-1.3,1.5,0)--(-1.3,1.5,1.3)--(-1.3,0,1.3)--cycle;
\filldraw[draw=none,fill=ink,fill opacity=.14] (1.3,0,0)--(1.3,1.5,0)--(1.3,1.5,1.3)--(1.3,0,1.3)--cycle;
\filldraw[draw=none,fill=volume,fill opacity=.05] (-1.3,0,1.3)--(1.3,0,1.3)--(1.3,1.5,1.3)--(-1.3,1.5,1.3)--cycle;
\draw[hidden] (-1.3,0,0)--(-1.3,1.5,0)--(1.3,1.5,0);
\draw[hidden] (-1.3,1.5,0)--(-1.3,1.5,1.3);
\draw[boxedge] (-1.3,0,0)--(1.3,0,0)--(1.3,1.5,0)--(1.3,1.5,1.3)--(-1.3,1.5,1.3)--(-1.3,0,1.3)--cycle;
\draw[boxedge] (-1.3,0,1.3)--(1.3,0,1.3)--(1.3,1.5,1.3);
\draw[boxedge] (1.3,0,0)--(1.3,0,1.3);
\draw[core,line width=2.3pt] (-1.3,0,0)--(1.3,0,0);
\fill[core] (-1.3,0,0) circle[radius=2pt];
\fill[core] (1.3,0,0) circle[radius=2pt];
\fill[white] (0,0,0) circle[radius=4pt];\fill[core] (0,0,0) circle[radius=2.8pt];
\coordinate (lefttarget) at (-1.3,0.75,.7);
\coordinate (righttarget) at (1.3,0.75,.7);
\coordinate (bottomtarget) at (.9,.65,0);\coordinate (fronttarget) at (-.6,0,.75);
\begin{scope}[x={(1cm,0cm)},y={(0cm,1cm)}]
\node[lab,anchor=east,text=ink!85!black] (leftlabel) at (-4.75,1.55) {$\{-1\}\times C$};
\draw[leader,ink!85!black] (leftlabel.east)--(lefttarget);
\node[lab,anchor=west,text=ink!85!black] (rightlabel) at (4.85,1.55) {$\{1\}\times C$};
\draw[leader,ink!85!black] (rightlabel.west)--(righttarget);
\node[lab,text=donorTwo!80!black] at (-3.25,2.3) {$t_2=0$};
\node[lab,text=donorTwo!80!black,font=\normalsize] at (-.2,.93) {$U\cap N_2$};
\node[lab,text=donorThree!85!black] at (-3.75,-.6) {$t_3=0$};
\node[lab,anchor=west,text=donorThree!85!black] (n3label) at (3.35,-.5) {$U\cap N_3$};
\draw[leader,donorThree!85!black] (n3label.west)--(bottomtarget);
\node[lab,text=core,anchor=north] (hlabel) at (-2,-1.5) {$H=U\cap N_{123}$};
\draw[leader,core] (hlabel.north)--(-2,-.4)--(-1.5,0);
\node[lab,text=core,anchor=north] (dlabel) at (.9,-1.5) {$\Delta=\{0\}$};
\draw[leader,core] (dlabel.north)--(.9,-.5)--(.08,-.05);
\coordinate (O) at (6.45,-1.9);
\draw[leader,ink!80] (O)--++(.8,0) node[right] {$y$};
\draw[leader,ink!80] (O)--++(.52,.36) node[right] {$t_2$};
\draw[leader,ink!80] (O)--++(0,.85) node[above] {$t_3$};
\end{scope}
\end{tikzpicture}
\caption{Here $C=[0,1]^2$, $J=\{1,2,3\}$, $p=0$, and $h=1$. The whole $U=D^0\times D^1\times C$ lies in $N_1$. The sectors $N_2$ and $N_3$ meet it in the front face $\{t_2=0\}$ and the bottom face $\{t_3=0\}$, respectively. Their common red edge is the handle $H\subset N_{123}$; its midpoint is the cocore $\Delta=D^0\times\{0\}$. The side faces are $D^0\times\partial D^1\times C$. The top and back faces of $U$ are given by the collar cutoffs $\{t_3=1\}$ and $\{t_2=1\}$. Every receiving sector is disjoint from $U$.}
\label{fig:compatible-handle-p0-h1}
\end{figure}

\begin{proof} We first arrange the free boundary. Let $e\co D^p\times D^h\to H$ be a parametrization of the terminal handle such that $e|D^p\times \{0\}$ is a parametrization of $\Delta$. Suppose that $p, h>0$. The free boundary $R_0=e(S\times D^h)$ of $H$ is a trivial regular neighborhood of $S$ in $\partial N_J$.  If there is only one receiving sector, take $R=R_0$ and $r=e|S\times D^h$. Otherwise, by \cite[Lemma 3.6]{Si} recalled above, there is a  regular neighborhood $R$ of $S$ in $\partial N_J$ and a PL homeomorphism $r\co S\times D^h\to R$ such that $r(x, 0)=x$ for each $x\in S$, and $r(E_i\times D^h)=R\cap N_i$ for every $i>j$.    

Choose a triangulation $T$ of $S$ with respect to which every $E_i$ is a subcomplex. The parametrizations $e|S\times D^h$ and $r$ give $R_0$ and $R$ trivial block bundle structures over $T$. By Theorem~\ref{th:RS4.4}, applied to the closed manifold $\partial N_J$, there is an ambient isotopy $G_t$, with $t\in [0,1]$, of $\partial N_J$ fixed on $S$ that takes corresponding blocks onto one another. In particular, $G_1(e(E_i\times D^h))=R\cap N_i$. The isotopy $G_t$ extends to an ambient isotopy $\bar{G}_t$ of $N_J$ that is identity outside a collar neighborhood of $\partial N_J$ in $N_J$ and that fixes $\Delta$ pointwise. Replacing $H$, $A$, and $e$ by $\bar{G}_1(H)$, $\bar{G}_1(A)$, and $\bar{G}_1\circ e$, we may assume that the free boundary of $H$ is a compatible regular neighborhood of $S$ in $\partial N_J$. Since the receiving sectors $N_i$, where $i>j$, meet $H$ only in its free boundary, $e^{-1}(N_i)=E_i\times D^h$ for $i>j$. 

If $p=0$, the free boundary is empty, so $H$ misses every receiving sector. If $h =0$, then $H = \Delta$ and the required equality follows directly from the definition of $E_i$. In both cases no adjustment is needed. Thus, in all cases there is a parametrization $e$ of $H$ such that $e^{-1}(N_i) = E_i\times D^h$ for every $i > j$.

We next extend $e$ beyond the attaching face. If $h>0$, choose $\lambda\in (0,1)$ and put $H'=e(D^p\times \lambda D^h)$. Both $H$ and $H'$ are trivial block bundle neighborhoods of $\Delta$ in $N_J$, and their restrictions over $S$ are their intersections with $\partial N_J$. By Theorem~\ref{th:RS4.4}(b), there is an ambient isotopy of $N_J$ fixed on $\Delta$ that takes $H$ to $H'$. Transporting the handle decomposition of $N_J$ by this isotopy makes $H'$ a terminal handle. The restriction $e|D^p\times \Int(D^h)$ extends the parametrization $e|D^p\times \lambda D^h$ of the handle $H'$. We now denote $H'$ by $H$ and relabel the pair $(\Int D^h, \lambda D^h)$ as $(\widetilde D^h, D^h)$. With this notation we have $e(D^p\times \widetilde D^h)\subset N_J$ and $e(D^p\times D^h)=H$. For $h=0$, set $\widetilde D^h=D^h$ and keep $e$ unchanged.

We now add the collar coordinates $t_2, \ldots, t_j$, one for each donor sector other than $N_1$. Put
\[
\begin{aligned}
Q_k&=N_{\{k,\ldots,j\}}, && 1\le k\le j,\\
P_k&=\partial Q_{k+1}, && 1\le k<j, \\
 P_j&=M. &&
\end{aligned}
\]
For $1\le k<j$, the boundary $P_k$ of $Q_{k+1}=N_{\{k+1,\ldots,j\}}$ is covered by its intersections with the sectors $N_i$ with $i\notin \{k+1, \ldots, j\}$. Hence the induced filling of $P_k$ consists of $Q_k=P_k\cap N_k$, and the sectors $P_k\cap N_i$ with $i<k$ or $i>j$. For $k=j$, we have $P_j=M$ and $Q_j=N_j$. Note that $Q_1\subset\cdots\subset Q_j$ and $Q_k\subset P_k$.  

Use a compatible exterior collar of $Q_1=N_J$ in $P_1$ to extend $e$ to an embedding $f_1\co \widetilde D^p\times \widetilde D^h\to P_1$ such that $f_1^{-1}(Q_1)=D^p\times \widetilde D^h$, and $f_1^{-1}(N_i)=\widetilde E_i\times \widetilde D^h$ for $i>j$. The preimages of $Q_1$ and $N_i$ meet precisely along $E_i\times \widetilde D^h$. When $p=0$, no exterior collar is needed, and we take $f_1=e$. 

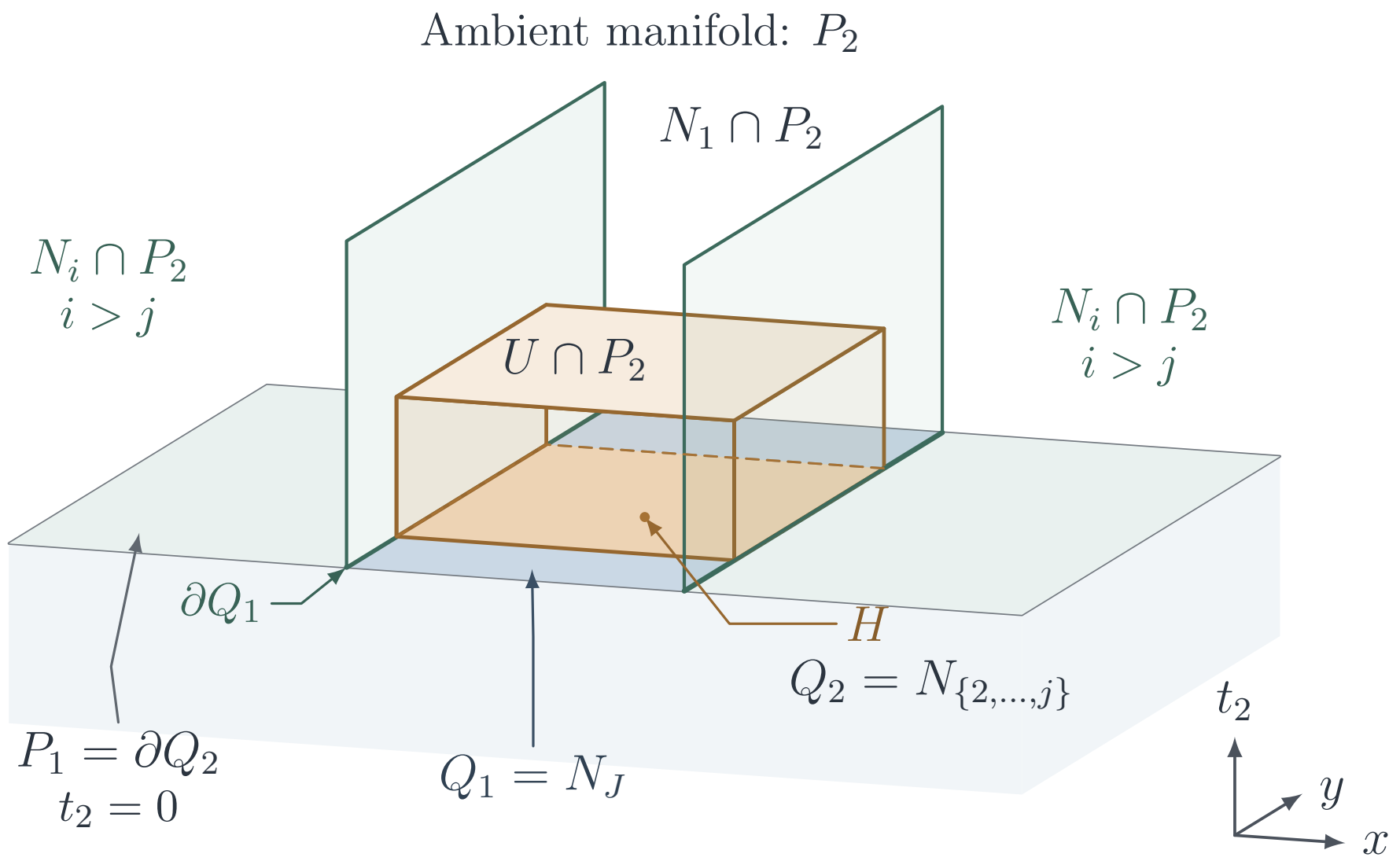
\begin{figure}[!htbp]
\centering
\begin{tikzpicture}[
  font=\large, line cap=round, line join=round,
  x={(1.42cm,-.10cm)},
  y={(.70cm,.43cm)},
  z={(0cm,1.12cm)}]
  \definecolor{ink}{RGB}{43,52,63}
  \definecolor{donors}{RGB}{105,143,181}
  \definecolor{receivers}{RGB}{86,151,132}
  \definecolor{transfer}{RGB}{208,143,65}
  \tikzset{
    edge/.style={draw=ink!65,line width=.65pt},
    hidden/.style={draw=ink!48,line width=.55pt,densely dashed},
    wall/.style={draw=receivers!70!black,line width=.8pt},
    boxedge/.style={draw=transfer!72!black,line width=.9pt},
    label/.style={text=ink,inner sep=2pt,align=center},
    leader/.style={->,draw=ink!75,line width=.6pt,>=latex},
    axis/.style={->,draw=ink!85,line width=.7pt,>=latex}
  }

  \fill[donors!9]
    (-3,-1.5,-1.35)--(3,-1.5,-1.35)--(3,-1.5,0)--(-3,-1.5,0)--cycle;
  \fill[donors!9]
    (3,-1.5,-1.35)--(3,1.6,-1.35)--(3,1.6,0)--(3,-1.5,0)--cycle;

  \filldraw[edge,fill=donors!12]
    (-3,-1.5,0)--(3,-1.5,0)--(3,1.6,0)--(-3,1.6,0)--cycle;
  \fill[receivers!13]
    (-3,-1.5,0)--(-1,-1.5,0)--(-1,1.6,0)--(-3,1.6,0)--cycle;
  \fill[receivers!13]
    (1,-1.5,0)--(3,-1.5,0)--(3,1.6,0)--(1,1.6,0)--cycle;
  \fill[donors!35]
    (-1,-1.5,0)--(1,-1.5,0)--(1,1.6,0)--(-1,1.6,0)--cycle;

  \draw[wall,line width=1.25pt] (-1,-1.5,0)--(-1,1.6,0);
  \draw[wall,line width=1.25pt] (1,-1.5,0)--(1,1.6,0);

  \fill[receivers!9]
    (-1,-1.5,0)--(-1,1.6,0)--(-1,1.6,2.45)--(-1,-1.5,2.45)--cycle;
  \draw[wall]
    (-1,-1.5,0)--(-1,1.6,0)--(-1,1.6,2.45)--(-1,-1.5,2.45)--cycle;

  \fill[transfer!35]
    (-1,-.9,0)--(1,-.9,0)--(1,.9,0)--(-1,.9,0)--cycle;
  \draw[hidden] (-1,-.9,0)--(-1,.9,0)--(1,.9,0);
  \draw[hidden] (-1,.9,0)--(-1,.9,1.05);
  \filldraw[boxedge,fill=transfer,fill opacity=.08,draw opacity=1]
    (-1,-.9,0)--(-1,.9,0)--(-1,.9,1.05)--(-1,-.9,1.05)--cycle;
  \filldraw[boxedge,fill=transfer,fill opacity=.06,draw opacity=1]
    (-1,-.9,0)--(1,-.9,0)--(1,-.9,1.05)--(-1,-.9,1.05)--cycle;
  \filldraw[boxedge,fill=transfer,fill opacity=.06,draw opacity=1]
    (1,-.9,0)--(1,.9,0)--(1,.9,1.05)--(1,-.9,1.05)--cycle;
  \filldraw[boxedge,fill=transfer!17]
    (-1,-.9,1.05)--(1,-.9,1.05)--(1,.9,1.05)--(-1,.9,1.05)--cycle;

  \fill[receivers,opacity=.07]
    (1,-1.5,0)--(1,1.6,0)--(1,1.6,2.45)--(1,-1.5,2.45)--cycle;
  \draw[wall]
    (1,-1.5,0)--(1,1.6,0)--(1,1.6,2.45)--(1,-1.5,2.45)--cycle;

  \draw[hidden,draw=transfer!80!black]
    (-1,-.9,0)--(-1,.9,0)--(1,.9,0);

  \node[label,text=receivers!65!black] at (-2.95,-.4,1.5)
    {$N_i\cap P_2$\\[-1pt]\normalsize $i>j$};
  \node[label,text=receivers!65!black] at (2.9,0,1.5)
    {$N_i\cap P_2$\\[-1pt]\normalsize $i>j$};

  \coordinate (targetPone) at (-2.4,-1.15,0);
  \coordinate (targetQone) at (0,-1.30,0);
  \coordinate (targetBoundary) at (-1,-1.5,0);
  \coordinate (targetHandle) at (.15,-.25,0);

  \begin{scope}[x={(1cm,0cm)},y={(0cm,1cm)}]
    \node[label,font=\normalsize] at (0,3.95)
      {Ambient manifold: $P_2$};
    \node[label] at (.85,3.15) {$N_1\cap P_2$};
    \node[label] at (-.55,1.20) {$U\cap P_2$};
    \node[label,anchor=north east] (labelPone) at (-3.45,-1.85)
      {$P_1=\partial Q_2$\\[-1pt]\normalsize $t_2=0$};
    \draw[leader] (labelPone.north)--(-4.45,-1.38)--(targetPone);
    \node[label,anchor=north,text=donors!45!black] (labelQone) at (-.90,-2.05)
      {$Q_1=N_J$};
    \draw[leader,draw=donors!55!black]
      (labelQone.north)--(-.90,-1.13)--(targetQone);
    \node[label,text=ink] at (2.45,-1.54)
      {$Q_2=N_{\{2,\ldots,j\}}$};
    \node[label,anchor=east,font=\normalsize,text=receivers!65!black]
      (labelBoundary) at (-3.10,-.85) {$\partial Q_1$};
    \draw[leader,draw=receivers!65!black]
      (labelBoundary.east)--(-2.85,-.85)--(targetBoundary);
    \node[label,anchor=west,text=transfer!65!black] (labelHandle) at (1.65,-1.02)
      {$H$};
    \draw[leader,draw=transfer!72!black]
      (labelHandle.west)--(.75,-1.02)--(targetHandle);
    \fill[transfer!80!black] (targetHandle) circle[radius=1.2pt];

    \coordinate (O) at (5.0,-2.8);
    \draw[axis] (O)--++(.93,-.065) node[right,text=ink] {$x$};
    \draw[axis] (O)--++(.57,.35) node[right,text=ink] {$y$};
    \draw[axis] (O)--++(0,.83) node[above,text=ink] {$t_2$};
  \end{scope}
\end{tikzpicture}
\caption{The first collar step, viewed locally in $P_2$ (here $p=h=1$).
The horizontal plane is $P_1=\partial Q_2$, and its central strip is
$Q_1=N_J$. 
The brown block is $U\cap P_2$, whose entire bottom face is
$H=U\cap P_1$. When $j=2$, we have $P_2=M$, $Q_2=N_2$, and the brown
block is $U$.}
\label{fig:sectors_near_U}
\end{figure}

\begin{figure}[!htbp]
\centering
\begin{tikzpicture}[
  font=\large, line cap=round, line join=round,
  x={(1.42cm,0cm)},
  y={(.70cm,.60cm)},
  z={(0cm,1.12cm)}]
  \definecolor{ink}{RGB}{43,52,63}
  \definecolor{donors}{RGB}{105,143,181}
  \definecolor{receivers}{RGB}{86,151,132}
  \definecolor{transfer}{RGB}{208,143,65}
  \definecolor{core}{RGB}{177,49,65}
  \tikzset{
    edge/.style={draw=ink!60,line width=.65pt},
    hidden/.style={draw=ink!40,line width=.55pt,densely dashed},
    seam/.style={draw=receivers!70!black,line width=1pt},
    boxedge/.style={draw=transfer!75!black,line width=1pt},
    label/.style={text=ink,inner sep=2pt,align=center},
    leader/.style={->,draw=ink!75,line width=.6pt,>=latex},
    axis/.style={->,draw=ink!85,line width=.7pt,>=latex}
  }

  \fill[ink!7]
    (-3,-1.6,-1.25)--(3,-1.6,-1.25)--(3,-1.6,0)--(-3,-1.6,0)--cycle;
  \fill[ink!10]
    (3,-1.6,-1.25)--(3,2.3,-1.25)--(3,2.3,0)--(3,-1.6,0)--cycle;
  \draw[edge] (-3,-1.6,-1.25)--(-3,-1.6,0);
  \draw[edge] (3,-1.6,-1.25)--(3,-1.6,0);

  \filldraw[edge,fill=donors!24]
    (-3,-1.6,0)--(3,-1.6,0)--(3,0,0)--(-3,0,0)--cycle;
  \filldraw[edge,fill=receivers!12]
    (-3,0,0)--(3,0,0)--(3,2.3,0)--(-3,2.3,0)--cycle;
  \fill[donors!7]
    (-1,0,0)--(1,0,0)--(1,2.3,0)--(-1,2.3,0)--cycle;
  \draw[edge] (-3,2.3,0)--(3,2.3,0);

  \draw[seam] (-1,0,0)--(-1,2.3,0);
  \draw[seam] (1,0,0)--(1,2.3,0);

  \filldraw[boxedge,fill=transfer!31]
    (-1,0,0)--(1,0,0)--(1,1.25,0)--(-1,1.25,0)--cycle;

  \draw[draw=ink!65,line width=.95pt] (-3,0,0)--(3,0,0);
  \draw[draw=core,line width=2.15pt] (-1,0,0)--(1,0,0);
  \fill[core] (-1,0,0) circle[radius=2.2pt];
  \fill[core] (1,0,0) circle[radius=2.2pt];

  \draw[->,>=latex,draw=transfer!65!black,line width=1.15pt]
    (.68,.10,0)--(.68,1.06,0);
  \draw[->,>=latex,draw=donors!40!black,line width=1.05pt]
    (-.62,-.08,0)--(-.62,-.94,0);

  \node[label] at (0,1.94,0) {$N_1\cap P_2$};
  \node[label,text=receivers!65!black] at (-2.1,1.40,0) {$N_4\cap P_2$};
  \node[label,text=receivers!65!black] at (2.1,1.40,0) {$N_4\cap P_2$};
  \node[label,text=transfer!60!black] at (-.27,.83,0) {$U\cap P_2$};
  \node[label,text=donors!40!black] at (-1.45,-1.05,0) {$Q_2=N_{23}$};
  \node[label] at (1.45,-1.6,-.67) {$Q_3=N_3$};

  \coordinate (targetPtwo) at (2.50,2.02,0);
  \coordinate (targetPone) at (-2.68,0,0);
  \coordinate (targetHandle) at (-.06,0,0);

  \begin{scope}[x={(1cm,0cm)},y={(0cm,1cm)}]
    \node[label] at (0,2.08) {$k=2$};
    \node[label,anchor=south] (labelPtwo) at (4.17,1.94)
      {$P_2=\partial Q_3$};
    \draw[leader] (labelPtwo.south)--(4.65,1.70)--(targetPtwo);

    \node[label,anchor=north east] (labelPone) at (-3.83,-2.55)
      {$P_1=\partial Q_2$};
    \draw[leader] (labelPone.north)--(-4.55,-1.12)--(targetPone);

    \node[label,anchor=north,text=core] (labelHandle) at (-.75,-2.64)
      {$H=\Delta\subset Q_1=N_{123}$};
    \draw[leader,draw=core]
      (labelHandle.north)--(-.75,-1.18)--(targetHandle);

    \node[label,font=\normalsize,text=transfer!65!black]
      at (2.62,.25) {$t_2>0$};
    \node[label,font=\normalsize,text=donors!40!black]
      at (-1.92,-.23) {$t_2<0$};

    \coordinate (O) at (5.10,-2.57);
    \draw[axis] (O)--++(.88,0) node[right,text=ink] {$x$};
    \draw[axis] (O)--++(.58,.50) node[right,text=ink] {$t_2$};
    \draw[axis] (O)--++(0,.88) node[above,text=ink] {$t_3$};
  \end{scope}
\end{tikzpicture}
\caption{The collar construction at $k=2$, with $p=1$, $h=0$, and $j=3$.
The region $Q_3=N_3$ lies below the plane $P_2=\partial Q_3$.
Within that plane, $Q_2=N_{23}$ occupies the side $t_2\leq0$, and
$P_1=\partial Q_2$ is the dividing line. The red interval is
$H=\Delta$, a component of $Q_1=N_{123}$.
The positive collar thickens $H$ to the brown rectangle
$U\cap P_2\subset N_1\cap P_2$. The negative collar enters $Q_2$.
Both collar directions at this stage lie in $P_2$. Above $P_2=\partial Q_3$ is the exterior of $Q_3$ occupied by the other sectors $N_1, N_2$ and $N_4$. At $k=2$, the construction stays within $P_2$.}
\label{fig:collar_h0_a2}
\end{figure}
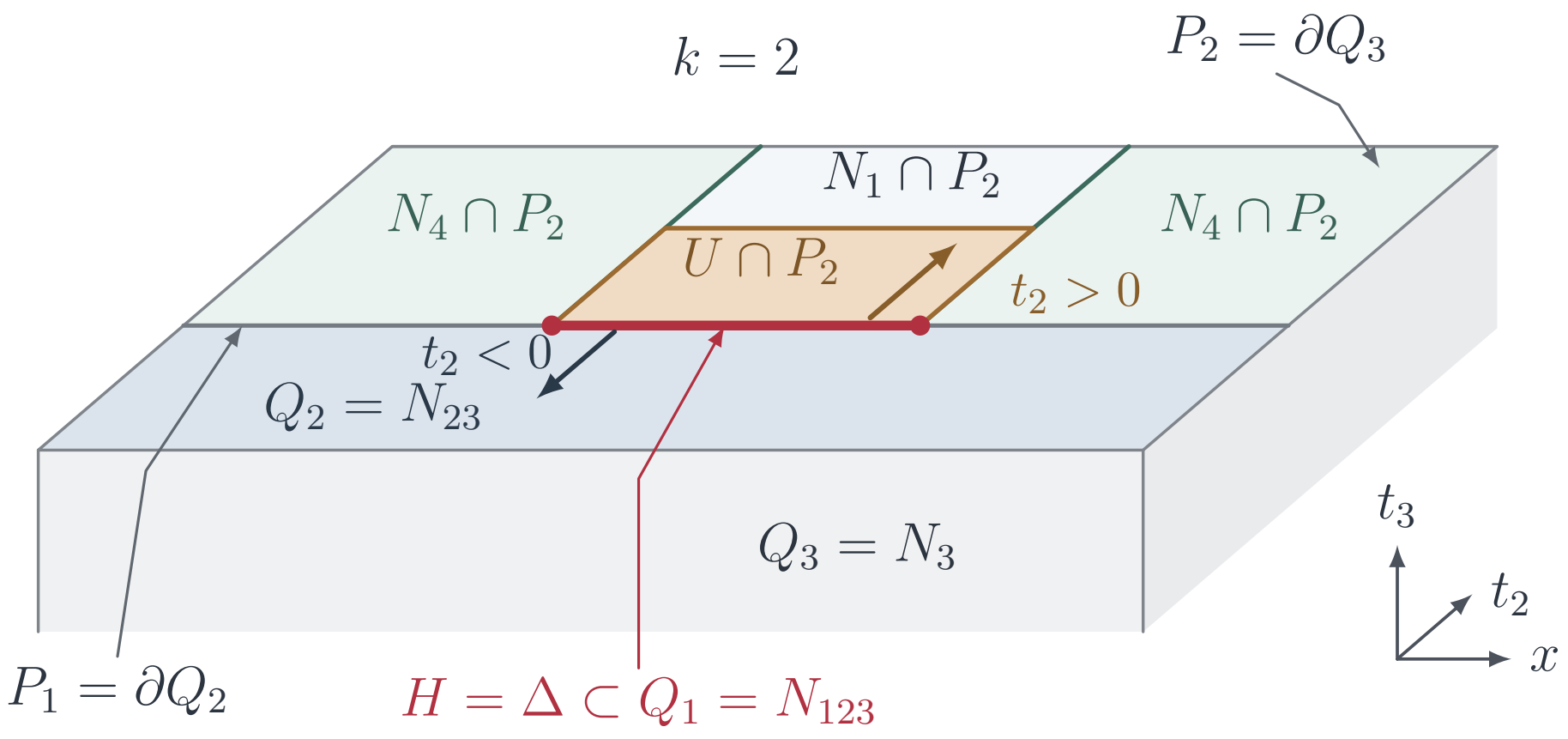

For each $k=2,\ldots,j$, regard $Q_k$ as a sector of the induced filling of $P_k$. Choose a bicollar of $\partial Q_k$ in $P_k$, with its negative half an interior collar in $Q_k$, and its positive half a compatible exterior collar of the sector $Q_k$ in $P_k$:
\[
b_k:P_{k-1}\times(-\varepsilon,1+\varepsilon)
\longrightarrow P_k,
\qquad b_k(z,0)=z.
\]
Compatibility of the exterior collar implies that for each $i<k$ or $i>j$, each $z\in P_{k-1}$, and each $t\in [0, 1+\varepsilon)$, a point $b_k(z, t)$ is in $P_k\cap N_i$ if and only if $z$ is in $P_k\cap N_i$. Starting with $f_1\co \widetilde D^p\times\widetilde D^h\to P_1$ define successively
\[
f_k:\widetilde D^p\times\widetilde D^h
\times(-\varepsilon,1+\varepsilon)^{k-1}
\longrightarrow P_k,
\]
\[
f_k(x,y,t_2,\ldots,t_k)
=
b_k\bigl(f_{k-1}(x,y,t_2,\ldots,t_{k-1}),t_k\bigr),
\]
for $k=2, \ldots, j$.
Each $f_k$ is well defined since $f_{k-1}$ takes values in $P_{k-1}$, which is the base of the bicollar $b_k$. It is a PL embedding since it is a composition of $f_{k-1}\times \operatorname{id}$ with $b_k$. Set $\widetilde\Phi=f_j$. Its domain is a manifold of dimension $m$ without boundary. Hence, its image $\widetilde U$ is open in $M$. Since every bicollar fixes its zero slice,
$\widetilde\Phi(D^p\times D^h\times\{0\})=H$. 

It remains to verify the sector formulas. Throughout this verification, sector labels refer to the original filling $M=N_1\cup\ldots\cup N_\ell$. Fix a point $\widetilde\Phi(x, y, t_2, \ldots, t_j)$. At stage $k$, write $z=f_{k-1}(x, y, t_2, \ldots, t_{k-1})$.
If $t_k<0$, then the point $b_k(z, t_k)=f_k(x, y, t_2, \ldots, t_k)$ lies in the interior of the sector $Q_k$ of $P_k$. The other sectors of the induced filling of $P_k$ are $P_k\cap N_i$ with $i<k$ or $i>j$. These meet $Q_k$ only along its boundary. Consequently, the point $b_k(z, t_k)$ is not in $N_i$ for $i<k$ or $i>j$. Since $Q_k=N_k\cap \ldots\cap N_j$, the point $b_k(z, t_k)$ belongs to each $N_i$ with $k\le i\le j$. Thus, immediately after a negative step at stage $k$, its sector labels are $\{k, k+1, \ldots, j\}$. 

If $t_k=0$, then $b_k(z, t_k)=z$ so its membership in every original sector $N_i$ is unchanged. If $t_k> 0$, then $b_k(z, t_k)$ retains the labels $\{k+1, \ldots, j\}$ since $P_k\subset N_i$ for every $k<i\le j$. Compatibility of the exterior collar also preserves its membership in every $N_i$ with $i<k$ or $i>j$. It does not lie to $N_k$ since it lies outside $Q_k=P_k\cap N_k$. Thus a positive step removes only the label $k$. 

Thus, membership of the initial point $f_1(x, y)$ in sector $N_1$, or in a receiving sector $N_i$ with $i>j$ persists precisely when all collar parameters $t_2, \ldots, t_j$ are nonnegative. Moreover, a label $1$ or $i>j$ that is absent initially cannot appear at a later stage. Since $f_1(x,y)\in N_1$ if and only if $x\in D^p$, and $f_1(x, y)\in N_i$ if and only if $x\in \widetilde E_i$ for $i>j$, this gives the asserted formulas for $N_1$ and the receiving sectors. 

Membership of the final point in a donor sector $N_\beta$ with $j\ge \beta\ge 2$ occurs precisely when
$t_\beta\le0$, and $t_k\ge 0$ for every $k>\beta$, with no restriction on earlier parameters. Indeed, immediately before stage $\beta$, the point lies in $P_{\beta-1}=\partial Q_\beta\subset N_\beta$, regardless of the earlier parameters. At stage $\beta$, membership in $N_\beta$ is retained exactly when $t_\beta\le0$. Every subsequent negative step removes this membership, whereas every subsequent zero or positive step preserves its presence or absence.
 Restricting to $D^p\times D^h\times C$ gives the stated formulas for the closed neighborhood $U\subset N_1$. Indeed, on this domain all collar parameters are nonnegative, so the donor condition reduces to $t_\beta=0$, while the receiving condition reduces to $x\in E_i$, since $D^p\cap\widetilde E_i=E_i$.

For a relative filling, choose a neighborhood of $H$ disjoint from the fixed collar of $\partial M$. The handle adjustments and collar constructions can be carried out within this neighborhood, with the isotopies extended by the identity outside it. Thus the same argument applies and preserves the prescribed boundary filling.
\end{proof}

\section{Trading}\label{sec:trading}

In this section we describe a modification that removes a terminal handle from an
intersection of sectors. The handle is thickened into one of these
sectors, and the resulting neighborhood is distributed among the
remaining sectors. This modification will allow us to
simplify intersections in increasing order of depth.

Let $M=N_1\cup\cdots\cup N_\ell$ be a filling of a compact PL
manifold of dimension $m$. Choose a nonempty proper set $J$ of donor labels and relabel
the sectors so that $J=\{1,\ldots,j\}$. Then the labels in $R=\{j+1,\ldots,\ell\}$ 
are the receiving labels. Choose a handle decomposition of the compact manifold $N_J$, and let $H$ be a terminal handle of index $h\ge 0$. Then, for some compact manifold $A$ of dimension $h+p$,  
\[
N_J=A\cup H,\qquad H=D^p\times D^h,\qquad p+h=m-j+1,
\]
\[
A\cap H=D^p\times\partial D^h,
\qquad
H\cap\partial N_J=\partial D^p\times D^h.
\]
The set $A\cap H$ is the attaching region of the handle, and $H\cap \partial N_J$ is the free
boundary of the handle. The disk $\Delta=D^p\times\{0\}$ is its cocore.
We allow $p=0$, $h=0$, or both, with $\partial D^0=\emptyset$. When $j=1$, we assume that $A\neq \emptyset$ so that every originally nonempty sector remains nonempty after the trade. 


If $M$ has boundary, we work with a relative filling that is a product on a prescribed collar of $\partial M$. We assume that the handle adjustment and the enlarged neighborhood used below are disjoint from this collar. After a small isotopy, assume that $S=\partial\Delta$ is transverse to the strata of the induced filling of $\partial N_J$. For $i\in R$, put $E_i=S\cap N_i$. The sets $E_i$, including empty labeled pieces, form the induced filling of $\partial D^p$.

\medskip
\noindent\emph{The neighborhood to be transferred.}
By Proposition~\ref{prop:compatible-handle-neighborhood}, we may adjust the terminal handle
by an ambient isotopy of $N_J$ so that there is an open neighborhood
\[
\widetilde\Phi:
\widetilde D^p\times\widetilde D^h\times\widetilde C
\longrightarrow\widetilde U,
\qquad
\widetilde C=(-\varepsilon,1+\varepsilon)^{j-1},
\]
of the traded region
\[
U=\widetilde\Phi(D^p\times D^h\times C),
\qquad
C=[0,1]^{j-1}.
\]
We identify the source of $\widetilde\Phi$ with its image $\widetilde U$ under $\widetilde\Phi$.
The coordinates of $C$ are indexed by $2,\ldots,j$, and the handle $H$
corresponds to $D^p\times D^h\times \{0\}$ in $U$. When $j=1$, $C$ is a point.
In these coordinates, 
\[
\begin{aligned}
N_1\cap\widetilde U
&=\{(x,y,t):x\in D^p,\ t_k\ge0\text{ for }2\le k\le j\},\\
N_\beta\cap\widetilde U
&=\{(x,y,t):t_\beta\le0,\ t_k\ge0\text{ for }\beta<k\le j\},
&&2\le\beta\le j,\\
N_i\cap\widetilde U
&=\{(x,y,t):x\in\widetilde E_i,\ t_k\ge0\text{ for }2\le k\le j\},
&&i\in R.
\end{aligned}
\]
Unmentioned coordinates are unrestricted within their respective factors. In particular, 
\begin{equation}\label{eq:N_i}
N_i\cap U=
\begin{cases}
D^p\times D^h\times\{t\in C:t_i=0\},&2\le i\le j,\\
E_i\times D^h\times C,&i\in R,
\end{cases}
\end{equation}
where the sets $E_i$ with $i\in R$ form the induced filling of $S$. All subsequent calculations near $U$ use these formulas on $\widetilde U$.

Suppose that the filling of $\partial D^p$ by the sectors $E_i$, with $i\in R$,
extends to a relative filling of $(D^p, \partial D^p)$ by the relative 
sectors $(B_i, E_i)$, with $i\in R$, and $B_i\cap\partial D^p=E_i$.
In a collar of $\partial D^p$, this extension is the product of the
given filling $\{E_i\}$ by an interval. Empty labeled pieces are permitted.
For every nonempty $T\subseteq R$, put
\[
B_T=\bigcap_{i\in T}B_i,\qquad E_T=\bigcap_{i\in T}E_i.
\]

Define a new decomposition of $M$ by redistributing $U$ among the receiving sectors:
\begin{equation}\label{eq:N_i'}
\begin{aligned}
N'_1&=\overline{N_1\setminus U},\\
N'_\beta&=N_\beta, &&2\le\beta\le j,\\
N'_i&=N_i\cup(B_i\times D^h\times C), &&i\in R.
\end{aligned}
\end{equation}
Thus only the distinguished donor sector $N_1$ is reduced.
The other donor sectors remain fixed, and the receiving sectors acquire
the pieces according to the decomposition of $D^p$ by $\{B_i\}$. 
Let $\partial_+C$ denote the union of the upper faces $\{t_\beta=1\}$ of $\partial C$, and put
\[
Q=(\partial D^h\times C)\cup(D^h\times\partial_+C).
\]
If $C$ is a point, put $\partial_+C=\emptyset$. Then the remaining part $N_1'$ of $N_1$ meets $U$ along
$N'_1\cap U=D^p\times Q$.
Every lower face with $t_\beta=0$ belongs to unchanged donor
sectors, and its intersection with $D^p\times Q$ also belongs to $N_1'$.

\begin{lemma}\label{l:prodcuct_filling}
Suppose that a decomposition $X=\bigcup P_i$ of a PL manifold satisfies the local PL filling condition.
Then the decomposition of $X\times\mathbb R$ given by the pieces $P_i\times[0,\infty)$
and one additional piece $X\times(-\infty,0]$ satisfies the local filling condition.
\end{lemma}

\begin{proof}
At a point with positive last coordinate, use the original local model
with one additional free coordinate. At a point with negative last
coordinate, only the additional piece occurs. At a point with last
coordinate zero, put the last coordinate first in a local filling chart.
Its negative half-space is the additional piece, and its nonnegative
half-space carries the original filling.
\end{proof}

\begin{proposition}\label{prop:trade-is-filling}
The modified sectors $N_1', \ldots, N_\ell'$ satisfy the local PL filling condition,
cover $M$, and have pairwise disjoint interiors.
Suppose that the original sectors are nonempty, and in addition that
$A\ne\emptyset$ when $j=1$. Then the new sectors are nonempty as well,
and hence form a filling with the same number $\ell$ of sectors.
The construction fixes the prescribed external boundary collar.
\end{proposition}
\begin{proof}
The pieces $B_i$ cover $D^p$ and have pairwise disjoint interiors.
Consequently, the new sectors cover $M$ and have pairwise disjoint
interiors. Away from $U$ the filling is unchanged.
It remains to check the local models along $U$. Throughout the proof we use the product coordinates on the open neighborhood $\tilde U$ given by Proposition~\ref{prop:compatible-handle-neighborhood}. 

\medskip
\noindent\emph{Step 1. The local filling condition when $C$ is a point.}
In this case $j=1$, and therefore there are no unchanged donors. We have
$U=H=D^p\times D^h$ and $Q=\partial D^h$. Write a point of $U$ as $(x,y)$.
If $y\in\operatorname{Int}D^h$ and $x\in\operatorname{Int}D^p$,
the new filling is the product of $\{B_i\}$ with free coordinates.
If $y\in\operatorname{Int}D^h$ and $x\in\partial D^p$, the product
collar of each $B_i$ agrees with the old receiving sector $N_i$ on the
other side of $\partial D^p$. Thus in this case the local model is the filling
$\{E_i\}$ with additional free coordinates.

Now suppose that $y\in\partial D^h$ and $x\in\operatorname{Int}D^p$. 
 Choose a signed normal coordinate
$s$ for this boundary, with $D^h$ on the side $s\ge0$.
The side $s\ge 0$ carries the
product filling $\{B_i\times D^h\}$, and the side $s\le 0$ is the remaining
donor $N'_1$. The pieces $B_i\times \partial D^h$ define a filling on $D^p\times \partial D^h$, and the decomposition near $(x,y)$ is obtained by extending this filling to the side $s\ge 0$ and assigning the entire side $s\le 0$ to $N_1'$ as in Lemma~\ref{l:prodcuct_filling}.

If $y\in\partial D^h$ and $x\in\partial D^p$, let $\rho$ be the inward collar coordinate
of $D^p$, and let $a$ denote coordinates along its boundary. In particular, $\rho> 0$ corresponds to $\Int D^p$ and $\rho=0$ corresponds to $\partial D^p$.
In the two normal coordinates $(\rho,s)$, the modified donor $N_1'$ occupies
the quadrant $\{\rho\ge 0,\ s\le 0\}$, the transferred region $U$ covers the region $\{\rho\ge 0, s\ge 0\}$, while the rest $\{\rho\le 0\}$ is covered by the receiving pieces. 
Their labels depend only on which $E_i$ contains $a$. 
A planar PL homeomorphism sends the donor quadrant $N_1'$ to a half-space, 
keeping all coordinates except $\rho$ and $s$ fixed.
We obtain the local model of Lemma~\ref{l:prodcuct_filling} applied to $\{E_i\}$, see Figure~\ref{fig:11}. 

If $h=0$, then $Q=\emptyset$, so only the cases with
$y\in\operatorname{Int}D^h$ occur. If $p=0$, the point $D^0$ belongs
to exactly one nonempty $B_i$, and the cases involving $\partial D^p$
are absent. This also includes $p=h=0$.
Thus the local filling condition holds when $C$ is a point.

\medskip
\noindent\emph{Step 2. Points at which all collar coordinates are positive.}
Let $z=(x, y, t)$ be a point in $D^p\times D^h\times C$.
Suppose that $t_\beta(z)>0$ for every $\beta=2,\ldots,j$.
Then none of the unchanged donor sectors $N_\beta$ meets a sufficiently small
neighborhood of the point. Thus only $N'_1$ and the receiving sectors
need to be considered. Write $w=(y, t)$ and regard $D^h\times C$ as a single PL ball. 

If $w\in\operatorname{Int}(D^h\times C)$, the first two cases of Step 1 apply. If $w\in\partial(D^h\times C)$, then $w\in Q$, because all lower faces are excluded by the positivity of the collar coordinates. Choose a local signed PL collar coordinate for the boundary of the whole product ball $D^h\times C$, with that ball on the nonnegative side. The last two cases of Step 1 now apply with $D^h\times C$ in place of $D^h$. 


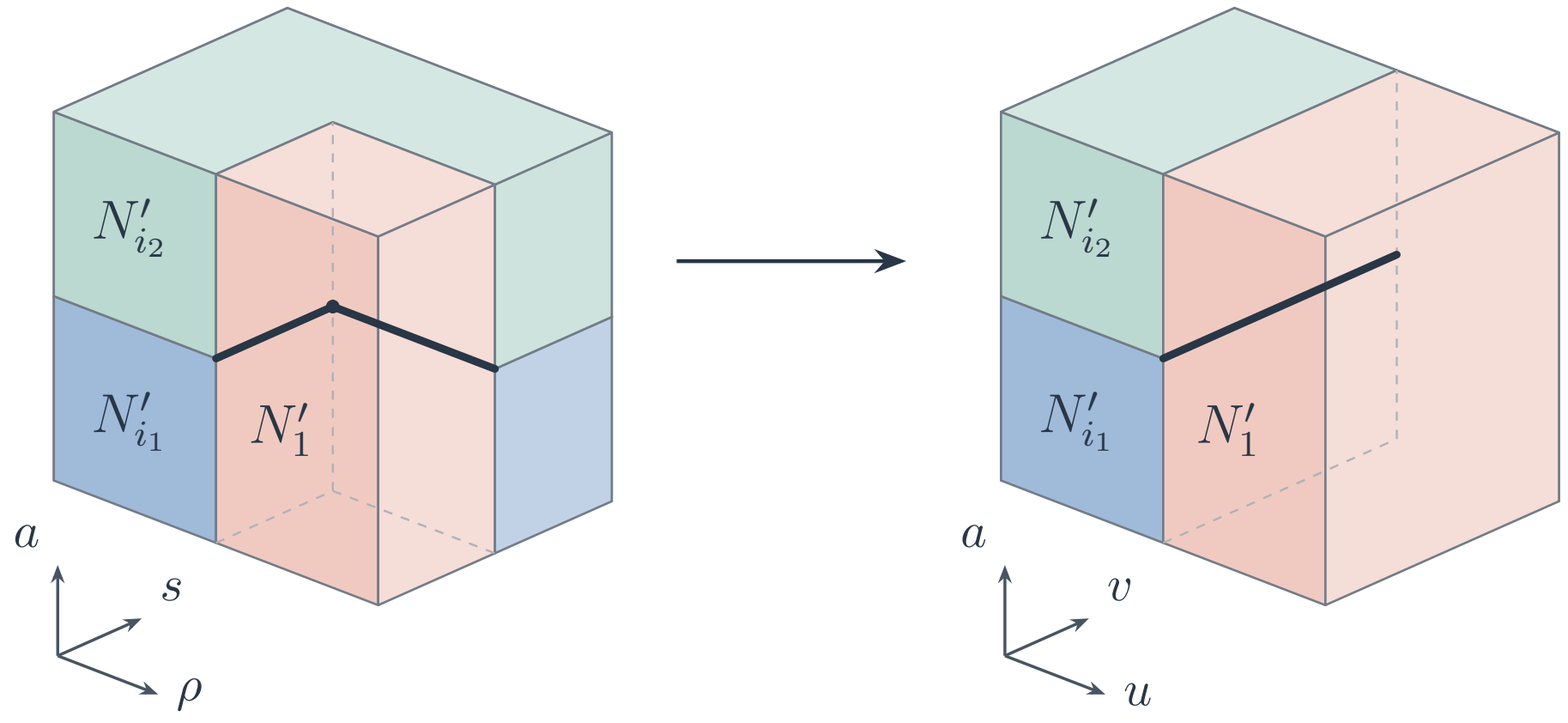
\begin{figure}
    \centering
    \definecolor{donor}{HTML}{D56C54}
\definecolor{lowerreceiver}{HTML}{497AB8}
\definecolor{upperreceiver}{HTML}{41947F}
\definecolor{ink}{HTML}{283646}

\tikzset{
  edge/.style={
    draw=ink!65,
    line width=.5pt,
    line join=round
  },
  faint/.style={
    draw=ink!38,
    line width=.45pt,
    dash pattern=on 2pt off 2pt
  },
  triple/.style={
    draw=ink,
    line width=1.7pt,
    line join=round,
    line cap=round
  },
  axes/.style={
    draw=ink!85,
    line width=.65pt,
    -{Stealth[length=4pt]}
  },
  flabel/.style={
    font=\large,
    inner sep=1.5pt,
    text=ink
  }
}
    \begin{tikzpicture}[font=\normalsize,text=ink,
  x={(1.25cm,-.48cm)},y={(.90cm,.40cm)},z={(0cm,1.42cm)}]

\begin{scope}
  \filldraw[edge,fill=lowerreceiver!13] (-1,1,-1)--(1,1,-1)--(1,1,0)--(-1,1,0)--cycle;
  \filldraw[edge,fill=upperreceiver!13] (-1,1,0)--(1,1,0)--(1,1,1)--(-1,1,1)--cycle;
  \filldraw[edge,fill=lowerreceiver!13] (-1,-1,-1)--(-1,1,-1)--(-1,1,0)--(-1,-1,0)--cycle;
  \filldraw[edge,fill=upperreceiver!13] (-1,-1,0)--(-1,1,0)--(-1,1,1)--(-1,-1,1)--cycle;

  \fill[lowerreceiver!12] (-1,-1,0)--(0,-1,0)--(0,0,0)--(1,0,0)--(1,1,0)--(-1,1,0)--cycle;
  \draw[faint] (-1,-1,0)--(-1,1,0)--(1,1,0)--(1,0,0);

  \filldraw[edge,fill=lowerreceiver!52] (-1,-1,-1)--(0,-1,-1)--(0,-1,0)--(-1,-1,0)--cycle;
  \filldraw[edge,fill=upperreceiver!36] (-1,-1,0)--(0,-1,0)--(0,-1,1)--(-1,-1,1)--cycle;
  \filldraw[edge,fill=donor!36] (0,-1,-1)--(1,-1,-1)--(1,-1,1)--(0,-1,1)--cycle;

  \filldraw[edge,fill=donor!23] (1,-1,-1)--(1,0,-1)--(1,0,1)--(1,-1,1)--cycle;
  \filldraw[edge,fill=lowerreceiver!34] (1,0,-1)--(1,1,-1)--(1,1,0)--(1,0,0)--cycle;
  \filldraw[edge,fill=upperreceiver!26] (1,0,0)--(1,1,0)--(1,1,1)--(1,0,1)--cycle;

  \filldraw[edge,fill=upperreceiver!22] (-1,-1,1)--(0,-1,1)--(0,0,1)--(1,0,1)--(1,1,1)--(-1,1,1)--cycle;
  \filldraw[edge,fill=donor!22] (0,-1,1)--(1,-1,1)--(1,0,1)--(0,0,1)--cycle;

  \draw[faint] (0,0,-1)--(0,0,1);
  \draw[faint] (0,-1,-1)--(0,0,-1)--(1,0,-1);
  \draw[triple] (0,-1,0)--(0,0,0)--(1,0,0);
  \fill[ink] (0,0,0) circle[radius=1.5pt];

  \node[flabel] at (-.53,-1,-.52) {$N'_{i_1}$};
  \node[flabel] at (-.53,-1,.53) {$N'_{i_2}$};
  \node[flabel] at (.40,-1,-.25) {$N'_1$};

\end{scope}

\begin{scope}[shift={(7.3cm,0cm)}]
  \filldraw[edge,fill=lowerreceiver!13] (-1,1,-1)--(0,1,-1)--(0,1,0)--(-1,1,0)--cycle;
  \filldraw[edge,fill=upperreceiver!13] (-1,1,0)--(0,1,0)--(0,1,1)--(-1,1,1)--cycle;
  \filldraw[edge,fill=donor!12] (0,1,-1)--(1,1,-1)--(1,1,1)--(0,1,1)--cycle;
  \filldraw[edge,fill=lowerreceiver!13] (-1,-1,-1)--(-1,1,-1)--(-1,1,0)--(-1,-1,0)--cycle;
  \filldraw[edge,fill=upperreceiver!13] (-1,-1,0)--(-1,1,0)--(-1,1,1)--(-1,-1,1)--cycle;

  \fill[lowerreceiver!12] (-1,-1,0)--(0,-1,0)--(0,1,0)--(-1,1,0)--cycle;
  \draw[faint] (-1,-1,0)--(-1,1,0)--(0,1,0);

  \filldraw[edge,fill=lowerreceiver!52] (-1,-1,-1)--(0,-1,-1)--(0,-1,0)--(-1,-1,0)--cycle;
  \filldraw[edge,fill=upperreceiver!36] (-1,-1,0)--(0,-1,0)--(0,-1,1)--(-1,-1,1)--cycle;
  \filldraw[edge,fill=donor!36] (0,-1,-1)--(1,-1,-1)--(1,-1,1)--(0,-1,1)--cycle;
  \filldraw[edge,fill=donor!23] (1,-1,-1)--(1,1,-1)--(1,1,1)--(1,-1,1)--cycle;
  \filldraw[edge,fill=upperreceiver!22] (-1,-1,1)--(0,-1,1)--(0,1,1)--(-1,1,1)--cycle;
  \filldraw[edge,fill=donor!22] (0,-1,1)--(1,-1,1)--(1,1,1)--(0,1,1)--cycle;

  \draw[faint] (0,1,-1)--(0,1,1);
  \draw[faint] (0,-1,-1)--(0,1,-1);
  \draw[triple] (0,-1,0)--(0,1,0);

  \node[flabel] at (-.53,-1,-.52) {$N'_{i_1}$};
  \node[flabel] at (-.53,-1,.53) {$N'_{i_2}$};
  \node[flabel] at (.40,-1,-.25) {$N'_1$};
\end{scope}

\begin{scope}[x={(1cm,0cm)},y={(0cm,1cm)}]
  \draw[axes] (-2.12,-2.69)--(-1.345,-2.988) node[right] {$\rho$};
  \draw[axes] (-2.12,-2.69)--(-1.472,-2.402) node[above right] {$s$};
  \draw[axes] (-2.12,-2.69)--(-2.12,-1.99) node[above left] {$a$};
  \draw[axes] (5.18,-2.69)--(5.955,-2.988) node[right] {$u$};
  \draw[axes] (5.18,-2.69)--(5.828,-2.402) node[above right] {$v$};
  \draw[axes] (5.18,-2.69)--(5.18,-1.99) node[above left] {$a$};
  \draw[-{Stealth[length=7pt]},line width=.85pt,ink] (2.65,.35)--(4.42,.35);
\end{scope}
     \end{tikzpicture}
   \caption{The local filling near $\partial D^p\times \partial D^h$ before (left) and after (right) straightening the corner of $N'_1$. The black line represents $N'_1\cap N'_{i_1}\cap N'_{i_2}$.}
    \label{fig:11}
\end{figure}

\medskip
\noindent\emph{Step 3. Points at which some collar coordinates are zero.}
Let $K_0$ denote the indices $\beta\in \{2, \ldots, j\}$ for which $t_\beta(z)=0$.  
Take a neighborhood of the point under consideration in which
$t_\beta<1$ for $\beta\in K_0$ and all remaining collar coordinates
are positive. First restrict to the slice $S$ on which $t_\beta=0$ for 
every $\beta\in K_0$, and consider the induced filling of $S$ by the 
sectors other than those indexed by $\beta\in K_0$.
The argument in Step 2 (or Step 1 if no collar direction remains) 
gives the local model on $S$ in the remaining coordinates.

Restore the omitted coordinates in the order used in constructing
the compatible neighborhood. For each restored coordinate $t_\beta$, extend
each piece already present by taking its product with the half-line $t_\beta\ge 0$, 
and assign the half space $t_\beta\le 0$ to the unchanged donor sector of label $\beta$. 
Each step is the construction described in Lemma~\ref{l:prodcuct_filling}. After all the omitted coordinates have been added, we recover the required local decomposition.

If $p=0$, every $E_i$ is empty and the point $D^0$ belongs to exactly
one nonempty $B_i$. The cases involving $\partial D^p$ are absent.
If $h=0$, the cases involving $\partial D^h$ are absent.
The same proof therefore covers both possibilities.

Finally, unchanged donor sectors remain nonempty and receiving sectors
only gain pieces. If $j=1$, then $N'_1=A$, which is nonempty by
hypothesis. If $j>1$, 
$N_1$ has points outside $U$ beyond an upper face, so $N'_1$ is nonempty. Thus the modified sets form a filling.
All modifications were made away from the boundary collar of $\partial M$, so that collar is fixed. 
\end{proof}

\begin{proposition}\label{prop:two-receiver-trade}
If $j=1$, assume that $A$ is non-empty. If $M$ has boundary, assume that the image of $\tilde\Phi$ is disjoint from a collar of $\partial M$. Then the filling $M=N_1'\cup\cdots\cup N_\ell'$ has the following properties: 
\begin{enumerate}[series=trade]
    \item The chosen terminal handle $H$ is removed from $N_J=A\cup H$ so that $N_J'\cong A$. 
    \item Every intersection $N_K$ of a nonempty proper subcollection of donor sectors is preserved up to a PL homeomorphism, i.e., $N'_K\cong N_K$ for $\emptyset\neq K\subsetneq J$. 
\item Suppose $K\subsetneq J$, including $K=\emptyset$, and $\emptyset\neq T\subseteq R$. Then, up to PL homeomorphism, $N'_{K\cup T}$ is obtained from $N_{K\cup T}$ by attaching handles with the same indices as those in a relative handle decomposition of $(B_T, E_T)$. Consequently, if $N_{K\sqcup T}\neq\emptyset$, then 
\[
\operatorname{gdim}N'_{K\cup T}
\leq
\max\!\left\{
   \operatorname{gdim}N_{K\cup T},
   \operatorname{gdim}(B_T,E_T)
\right\}.
\]
If $N_{K\cup T}=\emptyset$, then $E_T=\emptyset$ and $N'_{K\cup T}\cong B_T\times D^{h+j-1-|K|}$. In particular, a previously empty intersection may become nonempty, and when it does, $\operatorname{gdim}N'_{K\cup T}\le\operatorname{gdim}B_T$.
\item Suppose $\emptyset\neq T\subseteq R$. Then $N'_{J\cup T}$ is obtained from $N_{J\cup T}$ by replacing $E_T\times D^h$ with $B_T\times \partial D^h$ along their common attaching region $E_T\times \partial D^h$. 
\end{enumerate}
The construction fixes the prescribed collar of $\partial M$. In particular, the above estimates also hold with $\operatorname{gdim}(N_I,G_I)$ and $\operatorname{gdim}(N'_I,G_I)$ in place of the corresponding absolute geometric dimensions, where $G_I=N_I\cap\partial M=N'_I\cap\partial M$. 
\end{proposition}

\begin{remark}
\begin{itemize}
\item In Proposition~\ref{prop:two-receiver-trade}(1), we do not claim that $\operatorname{gdim} N'_J\le \operatorname{gdim} N_J$. For example, attaching a terminal $2$-handle to the inner boundary of an annulus produces a disk. Removing that handle raises geometric dimension from $0$ to $1$. What decreases in this case is the number of handles of top index in the handle decomposition of $N_J$. 
\item The replacement in Proposition~\ref{prop:two-receiver-trade}(4) may change the geometric dimension and may create a previously empty intersection. If $N_{J\cup T}=\emptyset$, then $E_T=\emptyset$ and
$N'_{J\cup T}\cong B_T\times\partial D^h$. No estimate analogous to (3) is asserted for these intersections.
\item In Proposition~\ref{prop:two-receiver-trade}(4), only intersections of higher depths are changed. 
\end{itemize}
\end{remark}

\begin{proof} The proof consists of several steps.

\vspace{5mm}
\noindent\emph{Step 1. The donor intersections.} Donor intersections indexed by nonempty subsets of $\{2, \ldots, j\}$ are unchanged as sets. Hence $N_K'=N_K$ if $\emptyset\neq K\subseteq J$ does not contain $1$. Suppose therefore that $K\subseteq J$ is a subset containing $1$. We will describe the modification of $N_K$. The boundary of $C$ consists of \emph{upper faces} $\{t_\beta=1\}$, whose union we denote by $C_+$, and the \emph{lower faces}, $\{t_\beta=0\}$ whose union is denoted by $C_-$. Here $2\le \beta\le j$, and we set $C_+=C_-=\emptyset$ when $j=1$.  
Similarly, the boundary of $D^h\times C$ consists of the upper faces $D^h\times C_+$, lower faces $D^h\times C_-$ and  \emph{side faces} $\partial D^h\times C$. Recall that $Q$ consists of the upper and side faces of the boundary of $D^h\times C$.  Define $D_K$ to be the subset of $D^h\times C$ with $t_\beta=0$ for every $\beta\in K\setminus \{1\}$. Then 
$D_K\cong D^h\times [0,1]^u$, where $u=j-|K|$, and the coordinates in $[0,1]^u$ correspond to the coordinates $t_\beta$ with $\beta\in J\setminus K$. 
Define
$Q_K=Q\cap D_K$. By the sector formulas \eqref{eq:N_i}, $N_i\cap U=\{t_i=0\}\cap U$ for $2\le i\le j$. Since $N_1\cap U=U$, we have $N_K\cap U = D^p\times D_K$. On the other hand, $N_1'\cap U=D^p\times Q$, and $N_i'=N_i$ for $2\le i\le j$. We have 
\[
\begin{aligned}
N'_K\cap U
&={(N'_1\cap U)\cap(N_K\cap U)}\\
&={(D^p\times Q)\cap(D^p\times D_K)}\\
&=D^p\times Q_K.
\end{aligned}
\]

\begin{lemma}\label{l:product_collar}
    If $1\in K\subsetneq J$, then $u\ge1$, and the pair $(D_K,Q_K)$ is a PL ball with a standard boundary ball.
\end{lemma}
\begin{proof} Let $\partial_-[0,1]^u=\bigcup_{\beta\in J\setminus K}\{t_\beta=0\}$ denote the union of the facets containing the vertex $(0, \ldots, 0)$. It is a standard ball embedded into the boundary of $[0, 1]^u$. Put $
 P_K=D^h\times\partial_-[0,1]^u$. Then $P_K$ is a standard boundary ball of $D_K$. Since $
Q_K=\overline{\partial D_K\setminus P_K}$, it is the complementary standard boundary ball. Hence $(D_K, Q_K)$ is a PL ball with a standard boundary ball. 
\end{proof}

Since $N_K\cap U=D^p\times D_K$, $N_K'\cap U=D^p\times Q_K$, and the two sets $N_K$ and $N_K'$ agree outside $U$,  we deduce that $N_K$ is obtained from $N_K'$ by attaching the ball $D^p\times D_K$ along $D^p\times Q_K$, i.e.,  
\begin{equation}\label{eq:trade-donors}
 N_K=N_K'\bigcup_{D^p\times Q_K}(D^p\times D_K).
\end{equation}
If $K\subsetneq J$, the trade does not change the PL homeomorphism type of $N_K$ by Lemma~\ref{l:product_collar}. 
For $K=J$, we have $D_J=D^h$ and $Q_J=\partial D^h$. Therefore, \eqref{eq:trade-donors} implies that $N_J$ is obtained from $N_J'$ by attaching exactly the terminal handle. 

\vspace{5mm}
\noindent\emph{Step 2. Intersections involving receiving sectors.}
Fix a nonempty $T\subseteq R$, and let $K\subseteq J$. Let $D_K=D^h\times [0,1]^u$, where $u=j-1-|K\setminus \{1\}|$ and the coordinates in $[0,1]^u$ correspond to the coordinates $t_\beta$ with $\beta\in J\setminus (K\cup \{1\})$. This extends the earlier definition of $D_K$ where we assumed that $1\in K$. For this calculation, we set $N_\emptyset=M$.
The equations \eqref{eq:N_i} immediately imply $N_K\cap U=D^p\times D_K$ and $N_T\cap U=E_T\times D^h\times C$. Consequently, $N_{K\cup T}\cap U=E_T\times D_K$. 

First suppose that $1\notin K$. The donor sectors with labels in $K$ are unchanged, and $N_{K\cup T}'\cap U=B_T\times D_K$. Since $N_{K\cup T}'$ agrees with $N_{K\cup T}$ outside $U$,
\begin{equation}\label{eq:trade-without-one}
 N_{K\cup T}'=N_{K\cup T}\bigcup_{E_T\times D_K}(B_T\times D_K).
\end{equation}
 This includes the case $K=\emptyset$, in which case $D_\emptyset =D^h\times C$. Thus, $N_{K\cup T}'$ is obtained from $N_{K\cup T}$ by attaching the relative pair $(B_T\times D_K, E_T\times D_K)$. This attachment can be carried out using handles with the same indices as those in the relative handle decomposition of $(B_T, E_T)$.   

Now suppose that $1\in K$. We have $N_{K\cup T}\cap U=E_T\times D_K$. Put $W=\overline{N_{K\cup T}\setminus U}$. Since \(W\subseteq N_{K\cup T}\cap N'_1\),
\[
{W\cap U}
{\subseteq}
(N_{K\cup T}\cap U)\cap(N'_1\cap U)
=
(E_T\times D_K)\cap(D^p\times Q)
=
E_T\times Q_K.
\]
Conversely, the product extensions beyond the side and upper faces (see Proposition~\ref{prop:compatible-handle-neighborhood}) show that every point of \(E_T\times Q_K\) is a limit of points of \(N_{K\cup T}\setminus U\). Therefore $W\cap U=E_T\times Q_K$. The definitions also give $N'_{K\cup T}\cap U=B_T\times Q_K$, with $E_T\subseteq B_T$, and the old and new intersections $N_{K\cup T}$ and $N'_{K\cup T}$ agree outside $U$.
Thus, $W$ intersects each of $N_{K\cup T}\cap U$ and $N'_{K\cup T}\cap U$ along $E_T\times Q_K$, and  
\begin{align}
N_{K\cup T}&=W\cup_{E_T\times Q_K}(E_T\times D_K),
       \label{eq:trade-old-mixed}\\
 N_{K\cup T}'&=W\cup_{E_T\times Q_K}(B_T\times Q_K).
       \label{eq:trade-new-mixed}
\end{align}
If $K\subsetneq J$, equation \eqref{eq:trade-old-mixed} and Lemma~\ref{l:product_collar} imply that $N_{K\cup T}$ is obtained from $W$ by attaching a product collar along $E_T\times Q_K$, so $W\cong N_{K\cup T}$. Equation \eqref{eq:trade-new-mixed} implies that $N_{K\cup T}'$ is obtained from $W$ by attaching the relative pair $(B_T,E_T)\times Q_K$. 
Since $Q_K$ is a disk of dimension $h+j-|K|-1$, this is a disk thickening of the prescribed relative pair. 

If $K=J$, the collar argument is not applicable: $Q_J=\partial D^h$, not a boundary ball of $D^h$. Instead, \eqref{eq:trade-old-mixed} and \eqref{eq:trade-new-mixed} read
\[
 N_{K\cup T}=W\cup_{E_T\times\partial D^h}(E_T\times D^h),
 \qquad
 N_{K\cup T}'=W\cup_{E_T\times\partial D^h}(B_T\times\partial D^h).
\]

All attachments and collar identifications fix the prescribed boundary collar, so the same argument proves the relative estimate.
\end{proof}

\begin{remark}
The disk factor for $K\subsetneq J$ has dimension $h+j-1-|K|$. For $K=J$, put $W_T=\overline{N_{J\cup T}\setminus(E_T\times D^h)}$. Suppose that $h\ge 1$, and choose a handle decomposition of $B_T$ relative to $E_T$ with all handle indices at most $c_T$, for some $c_T\ge 0$. 
By taking its product with $\partial D^h$, a $q$-handle of $(B_T,E_T)$ contributes two handles of indices $q$ and $q+h-1$ to $N'_{J\cup T}$ relative to $W_T$. Hence $\operatorname{gdim}N_{J\cup T}'\le\max\{\operatorname{gdim}W_T,c_T+h-1\}$. This does not give the same estimate with $N_{J\cup T}$ in place of $W_T$. If $h=0$, then $N'_{J\cup T}=W_T$ since $\partial D^0=\emptyset$.
\end{remark}

\section{Sufficiency of the partition inequalities}
\label{sec:pareto-realization}

Let $m\ge0$ and $\ell\ge1$. Recall that $\mathcal A_{m,\ell}$ consists of the nondecreasing tuples
$a=(a_1,\ldots,a_\ell)$ of nonnegative integers such that
\[
a_{\lambda_1}+\cdots+a_{\lambda_q}\ge m-\ell+1
\qquad
\text{whenever }\quad \lambda_1+\cdots+\lambda_q=\ell,\quad \lambda_i\ge1.
\]
We show that every tuple in $\mathcal A_{m,\ell}$ is a universal relative upper bound in dimension $m$. 
In particular, every tuple on the cohomological Pareto frontier is an optimal universal upper bound. 

We prove the absolute and relative statements simultaneously. Suppose that
$\partial M=G_1\cup\cdots\cup G_\ell$ is a prescribed filling, and put
$G_I=\bigcap_{i\in I}G_i$. An extension is required to agree with the product
filling $(G_i\times [0,1])$ on a collar of $\partial M$. We will use the following elementary observation to start the construction.

\begin{lemma}\label{lem:initial-relative-filling}
Every filling of the boundary of a compact PL manifold $M$ extends to a relative filling of $M$. 
\end{lemma}

\begin{proof}
Choose a collar $\partial M\times[0,2]$ in $M$, with
$\partial M=\partial M\times\{0\}$. For $i\ge2$ take
$N_i=G_i\times[0,1]$, and assign the rest of $M$ to $N_1$.
More precisely, $N_1$ is the union of $G_1\times[0,1]$ and
the complement of $\partial M\times[0,1)$. The local filling models at the inner boundary of the collar are verified exactly as in the second case of Singhof’s proof of Lemma (4.1) \cite[pp. 403–404]{Si}. Thus the pieces $N_1, N_2, \ldots$ form a relative filling extending the prescribed boundary filling.
\end{proof}

\begin{lemma}\label{lem:receiving-prefix}
Suppose that $a\in\mathcal A_{m,\ell}$, that $1\le j<\ell$, and that
$a_j<h\le m-j+1$. Put $p=m-j+1-h$, and $s=\ell-j$. 
Then $(a_1,\ldots,a_s)\in\mathcal A_{p,s}$.
\end{lemma}

\begin{proof}
The tuple $(a_1,\ldots,a_s)$ is nonnegative and nondecreasing. Let
$\lambda_1+\cdots+\lambda_q=s$ be a partition.
Then $\lambda_1+\cdots+\lambda_{q+1}$, where $\lambda_{q+1}=j$, is a partition of $\ell$, so
$a_{\lambda_1}+\cdots+a_{\lambda_q}\ge m-\ell+1-a_j$. 
Since $h\ge a_j+1$, we have
\[
p-s+1
=(m-j+1-h)-(\ell-j)+1
=m-\ell+2-h
\le m-\ell+1-a_j.
\]
Combining the two inequalities proves the assertion.
\end{proof}

\begin{theorem}\label{thm:pareto-realization}
Let $M$ be a compact PL $m$-manifold, and let
$a\in\mathcal A_{m,\ell}$. Every prescribed filling
$\partial M=G_1\cup\cdots\cup G_\ell$ extends to a relative filling
$M=N_1\cup\cdots\cup N_\ell$ that agrees with the product filling on a collar of $\partial M$ and 
such that $\operatorname{gdim}(N_I,G_I)\le a_{|I|}$
for every nonempty intersection $N_I$. In particular, if 
$M$ is closed, then
$\operatorname{gdim}N_I\le a_{|I|}$. 
\end{theorem}

The proof uses the compatible handle neighborhood and trading results
of the preceding sections.

\begin{proof}
We use induction on $\ell$, proving the relative statement in
all dimensions simultaneously. For $\ell=1$, condition~\eqref{eq:intro-partitions} gives $a_1\ge m$.
Take $N_1=M$. A compact PL $m$-manifold has a handle decomposition
relative to its boundary with all handles of index at most $m$, which
proves the assertion for $\ell=1$.

Suppose now that the assertion has been proved for every number of
labels less than $\ell$, in every dimension. Start with an extension
given by Lemma~\ref{lem:initial-relative-filling}, and fix a smaller
collar of $\partial M$ on which the filling will remain unchanged.
We simplify the intersections in increasing order of depth.

\medskip
\noindent\emph{The induction at depth $j$.}
Let $1\le j<\ell$, and suppose that
\begin{equation}\label{eq:already-controlled}
\operatorname{gdim}(N_I,G_I)\le a_{|I|}
\qquad\text{for }|I|<j.
\end{equation}
Choose finite relative handle decompositions of the pairs $(N_I, G_I)$ with $|I|=j$, with handles ordered by index. We will remove every handle whose index is greater than $a_j$. If one of these decompositions has such a handle, choose a terminal handle $H$ of maximal index $h>a_j$.
Let $J$ be its label set, so $|J|=j$, and write
\[
N_J=A\cup H,\qquad H=D^p\times D^h,\qquad p=m-j+1-h.
\]
The relative handle decomposition and the compatible neighborhood
construction can be chosen away from a collar of $\partial M$.
Thus the trade will leave the prescribed boundary filling unchanged.

\medskip
\noindent\emph{Filling the cocore disk.}
Let $R$ be the complement of $J$ and put $s=|R|=\ell-j$.
The receiving sectors induce a prescribed filling
$\partial D^p=\bigcup_{i\in R}E_i$.
By Lemma~\ref{lem:receiving-prefix}, the tuple
$(a_1,\ldots,a_s)$ belongs to $\mathcal A_{p,s}$.
Since $s<\ell$, the relative induction hypothesis applies to this
particular boundary filling of $D^p$. It gives an extension
\[
D^p=\bigcup_{i\in R}B_i,\qquad
\operatorname{gdim}(B_T,E_T)\le a_{|T|},
\quad \emptyset\ne T\subseteq R.
\]

Apply the trade \eqref{eq:N_i'} described in \S\ref{sec:trading}. Then the chosen
terminal handle disappears from $N_J$.

\medskip
\noindent\emph{Preserving the bounds already obtained.}
Consider any other intersection of depth at most $j$.
If all its labels belong to $J$, it is a proper donor intersection,
so its relative geometric dimension is unchanged.
Otherwise write its label set uniquely as $K\cup T$, where
$K\subseteq J$ and $\emptyset\ne T\subseteq R$.
Since $|K|+|T|\le j$, we have $K\subsetneq J$.
The new handles in this intersection have indices at most $a_{|T|}$.
The tuple $a$ is nondecreasing, and hence $a_{|T|}\le a_{|K|+|T|}$.
Consequently, every bound in~\eqref{eq:already-controlled} is preserved,
and no handle of index greater than $a_j$ is added to any other
intersection of depth $j$.

The intersections containing all labels in $J$ and at least one
receiving label have depth greater than $j$. They will be considered when their
depth is reached.

\medskip
\noindent\emph{Termination at depth $j$.}
The total number of handles of index greater than $a_j$, counted
over all intersections of depth $j$, decreases at least by one at each trade.
It is a finite nonnegative integer. Hence, after finitely many trades
we obtain
$\operatorname{gdim}(N_I,G_I)\le a_j$ for $|I|=j$, 
while all bounds at smaller depths remain valid. Carry out this procedure by induction for $j=1,\ldots,\ell-1$.

\medskip
\noindent\emph{The deepest intersection.}
The only remaining intersection is $N_{\{1,\ldots,\ell\}}$.
If it is nonempty, it has dimension $m-\ell+1$.
Every relative handle index is at most this dimension.
The one part partition $\ell$ in~\eqref{eq:intro-partitions}
gives $a_\ell\ge m-\ell+1$, and therefore
\[
\operatorname{gdim}
\bigl(N_{\{1,\ldots,\ell\}},G_{\{1,\ldots,\ell\}}\bigr)
\le m-\ell+1\le a_\ell.
\]

All changes were made away from the fixed boundary collar. This proves
the relative assertion, and the closed assertion follows by taking
$\partial M=\emptyset$. 
\end{proof}

\begin{proof}[Proof of Theorem~\ref{th:intro-main}]
Theorem~\ref{thm:pareto-realization} gives $\textup{(i)}\Rightarrow\textup{(iii)}\Rightarrow\textup{(ii)}$. Applying $\textup{(ii)}$ to $\mathbb{RP}^m$ gives $\textup{(ii)}\Rightarrow\textup{(iv)}$. Finally, Proposition~\ref{prop:rpm} gives
$\textup{(iv)}\Rightarrow\textup{(i)}$. Thus all four conditions are equivalent.
\end{proof}

\section{Lusternik-Schnirelmann category estimate}\label{s:ls_estimate}

Our approach draws on ideas from Lusternik-Schnirelmann category theory. Recall that a subset $A$ of a topological space $X$ is \emph{categorical} in $X$ if the inclusion $A\hookrightarrow X$ is null-homotopic. We use the convention that $\operatorname{cat} X$ is the least number of open categorical subsets covering $X$. 
In particular, $\operatorname{cat}(\mathrm{point})=1$. The following estimate is a consequence of the standard dimension bound and subadditivity properties of LS category~\cite{CLOT}. 

\begin{proposition}\label{prop:ls} Let $M$ be a compact connected PL manifold with a filling $M=N_1\cup\cdots\cup N_\ell$. If $\operatorname{gdim}N_i\le n_i$, where each $n_i$ is a nonnegative integer, then $\operatorname{cat}(M)\le\sum_{i=1}^{\ell}(n_i+1)$. In particular, if every sector has geometric dimension at most $n$, then $\ell(n+1)\ge\operatorname{cat}(M)$.
\end{proposition}
\begin{proof} Each connected component of $N_i$ has the homotopy type of a CW complex of dimension at most $n_i$, and hence admits a cover by at most $n_i+1$ open sets categorical in that component. Since $M$ is connected, the contraction points can be joined by paths to a common point in $M$. Thus $N_i$ is covered by at most $n_i+1$ open sets categorical in $M$. Choose an open neighborhood of $N_i$ that deformation retracts onto $N_i$. Pulling back the cover under the retraction gives at most $n_i+1$ open subsets of $M$, categorical in $M$, covering this neighborhood. Taking these covers for all nonempty sectors gives a categorical open cover of $M$ with at most $\sum_i(n_i+1)$ members. \end{proof}

\begin{corollary}\label{c:ls} Suppose that every closed PL manifold of dimension $m$ admits a filling into $\ell$ sectors of geometric dimension at most $n$. Then $\ell(n+1)\ge m+1$. 
\end{corollary}
\begin{proof} Apply Proposition~\ref{prop:ls} to a filling of $\mathbb{RP}^m$, using $\operatorname{cat}(\mathbb{RP}^m)=m+1$. 
\end{proof}

When $\operatorname{cat}(M)<\dim M+1$, the obstruction in Proposition~\ref{prop:ls} to the existence of fillings with small sector dimensions is weaker than the corresponding universal obstruction. Hence, such a manifold may admit fillings with smaller geometric-dimension bounds than those required for all manifolds of the same dimension.

\section{Appendix}

 A \emph{unit upper triangular PL transformation} is a map $\R^n\to \R^n$ given by $(t_1, \ldots t_n)\mapsto (t_1 + f_1(t_2,\ldots, t_n), \ldots, t_{n-1}+f_{n-1}(t_n), t_n+c)$, where $c\in \R$ and each $f_i$ is a PL function. Such a transformation is a PL homeomorphism.

\begin{proof}[Proof of Proposition~\ref{prop:filling}]
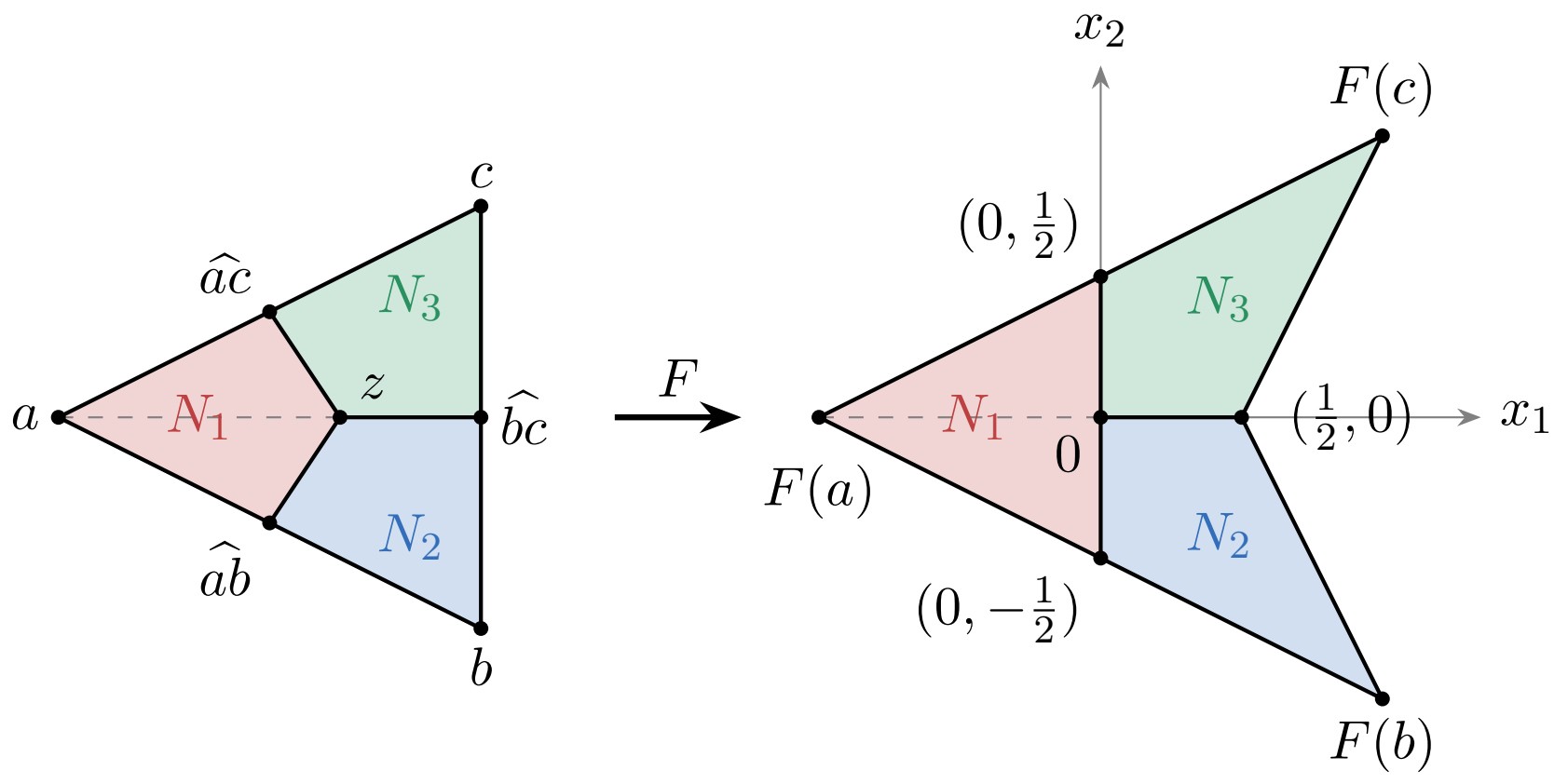
\begin{figure}[h]
\centering
\begin{tikzpicture}[>=Stealth]
\coordinate (a) at (0,1.5);\coordinate (b) at (3,0);\coordinate (c) at (3,3);
\coordinate (ab) at (1.5,.75);\coordinate (ac) at (1.5,2.25);
\coordinate (bc) at (3,1.5);\coordinate (z) at (2,1.5);
\fill[sectorone!22] (a)--(ab)--(z)--(ac)--cycle;
\fill[sectortwo!22] (b)--(bc)--(z)--(ab)--cycle;
\fill[sectorthree!22] (c)--(ac)--(z)--(bc)--cycle;
\draw[gray,dashed] (a)--(z);
\draw[thick] (a)--(b)--(c)--cycle;
\draw[thick] (ab)--(z)--(ac) (z)--(bc);
\foreach \p in {a,b,c,ab,ac,bc,z}{\fill (\p) circle (1.5pt);}
\node[left] at (a) {$a$};\node[below] at (b) {$b$};\node[above] at (c) {$c$};
\node[below left] at (ab) {$\widehat{ab}$};\node[above left] at (ac) {$\widehat{ac}$};
\node[right] at (bc) {$\widehat{bc}$};\node[above right] at (z) {$z$};
\node[text=sectorone] at (1,1.5) {$N_1$};
\node[text=sectortwo] at (2.5,.65) {$N_2$};
\node[text=sectorthree] at (2.5,2.35) {$N_3$};
\draw[->,very thick] (3.95,1.5)--(4.85,1.5) node[midway,above] {$F$};
\begin{scope}[shift={(7.4,1.5)},x=2cm,y=2cm]
\coordinate (A) at (-1,0);\coordinate (B) at (1,-1);\coordinate (C) at (1,1);
\coordinate (AB) at (0,-.5);\coordinate (AC) at (0,.5);\coordinate (BC) at (.5,0);\coordinate (O) at (0,0);
\fill[sectorone!22] (A)--(AB)--(O)--(AC)--cycle;
\fill[sectortwo!22] (AB)--(B)--(BC)--(O)--cycle;
\fill[sectorthree!22] (AC)--(O)--(BC)--(C)--cycle;
\draw[gray,dashed] (A)--(O);
\draw[thick] (A)--(B)--(BC)--(C)--cycle;
\draw[thick] (AB)--(AC) (O)--(BC);
\draw[->,gray] (BC)--(1.35,0) node[right,black] {$x_1$};
\draw[->,gray] (AC)--(0,1.25) node[above,black] {$x_2$};
\foreach \p in {A,B,C,AB,AC,BC,O}{\fill (\p) circle (1.5pt);}
\node[below=2mm] at (A) {$F(a)$};\node[below] at (B) {$F(b)$};\node[above] at (C) {$F(c)$};
\node[below left] at (AB) {$(0,-\tfrac12)$};
\node[above left] at (AC) {$(0,\tfrac12)$};
\node[right=2mm] at (BC) {$(\tfrac12,0)$};
\node[below left] at (O) {$0$};
\node[text=sectorone] at (-.45,0) {$N_1$};
\node[text=sectortwo] at (.42,-.42) {$N_2$};
\node[text=sectorthree] at (.42,.42) {$N_3$};
\end{scope}
\end{tikzpicture}
\caption{The transformation $F=(x_1, x_2)$ of sectors in a triangle $\tau$ whose vertices $(a, b, c)$ have colors $1, 2, 3$, respectively.}
\label{fig:transformation}
\end{figure}

By Lemma~\ref{prop:canonical}, the sectors cover $M$ and have pairwise disjoint interiors. Thus, it remains to show that for each point $z\in M$ there is a model neighborhood. 

For a vertex $v$ of $K'$, we define the \emph{barycentric distance to $v$} by $\delta_v(x)=1-\lambda_v(x)$, where $\lambda_v\co M\to [0,1]$ denotes the barycentric coordinate function of $v$ extended by zero over simplices not containing $v$.  Let $\delta_i(x)$ denote the barycentric distance from $x$ to the set of vertices of color $i$, i.e., $\delta_i(x)=\min \delta_v(x)$, where $v$ ranges over vertices of color $i$.   

Fix $z\in M$, and let $\tau$ be the simplex of $K'$ containing $z$ in its interior $\tau\setminus \partial\tau$. Let $A=A(z)$ consist of the vertices $v\in K'$ that minimize $\delta_v(z)$. These vertices $v$ belong to $\tau$. We may assume that the colors of vertices in $A$ are $1,\ldots,\nu$ for some $\nu\ge 1$. After choosing a sufficiently small neighborhood $U\subset M$ of $z$, no vertex outside $A$ attains the minimum anywhere in $U$. 
Specifically, put 
\begin{equation} \label{eq:U}
U=\left\{x\in M:
\max_{v\in A}\delta_v(x)
<
\min_{w\notin A}\delta_w(x)
\right\}.
\end{equation}

By continuity, the set $U$ is indeed an open neighborhood of $z$. At every point of $U$, the minimum among all vertex functions $\delta_v$ is assumed in $A$. Thus only $N_1,\ldots,N_\nu$ meet $U$. Let $m_j=\min\{\delta_{j}, \ldots, \delta_{\nu}\}$ denote the barycentric distance to the set of vertices of colors $j, \ldots, \nu$. For $j<\nu$, put
$x_j=\delta_j-m_{j+1}$.

\begin{lemma}\label{l:N}
    On $U$, the sectors are given by $N_{j}\cap U=\{x_1\ge 0, \ldots, x_{j-1}\ge 0, x_j\le 0\}$ for $j<\nu$, and $N_{\nu}\cap U=\{x_1\ge 0, \ldots, x_{\nu-1}\ge 0\}$.
\end{lemma}   
\begin{proof} For each $j$ the sector $N_j\cap U$ consists of the points $x\in U$ at which $\delta_j(x)=\min_{1\le i\le \nu} \delta_i(x)$, or equivalently $\delta_j(x)=m_1(x)$. We will show that this condition is equivalent to the stated inequalities.

Consider first a point in $N_j\cap U$. For every $i<j$, the equality $m_{i+1}=\min\{\delta_{i+1}, \ldots, \delta_j, \ldots, \delta_\nu\}$ implies that $m_{i+1}=\delta_j$, since $\delta_j$ is the minimum among all $\delta_i$. Hence $x_i=\delta_i-\delta_j\ge 0$ for all $i<j$. 
If $j<\nu$, then we also have $\delta_j\le m_{j+1}$, which implies $x_j\le 0$, see Fig.~\ref{fig:transformation}.

Conversely, suppose that the stated inequalities hold. For $i<j$, each of the inequalities $x_i\ge 0$ is equivalent to $\delta_i\ge m_{i+1}$. Since $m_{i}=\min\{\delta_{i}, m_{i+1}\}$, we deduce that $m_1=\cdots=m_j$. If $j<\nu$, the remaining condition $x_j\le 0$ is equivalent to $\delta_j\le m_{j+1}$. Therefore, $m_j=\min\{\delta_j, m_{j+1}\}=\delta_j$. If $j=\nu$, then $m_j=\delta_j$ holds by definition. Thus, $\delta_j=m_1$, and the point belongs to $N_j$.  
\end{proof}

For each color $i=1,\ldots, \nu$, let $A_i$ denote the set of vertices $v\in A$ of color $i$, and fix a representative vertex $v_i\in A_i$.

\begin{lemma}\label{l:PL_coordinates}
For $\nu>1$, the functions $x_1,\ldots, x_{\nu-1}$ extend to a PL coordinate system near $z$. 
\end{lemma}
\begin{proof}
Put $p=\dim\tau$, where $\tau$ is the simplex that contains $z$ in its interior. Then the restrictions of $(\delta_v)_{v\in I}$ to $\tau$ are affine coordinates near $z$ on $\tau$, where $I=\operatorname{Vert}(\tau)\setminus \{v_\nu\}$. For each $v\in I$, define a function on $U$ by
\[
y_v=
\begin{cases}
\delta_i - \min\{\delta_{i+1}, \ldots, \delta_{\nu}\}, &v=v_i, \quad \mathrm{for}\quad i=1,\ldots, \nu-1,\\
\delta_v-\delta_{v_i},&v\in A_i\setminus\{v_i\},\\
\delta_v,&v\in B, 
\end{cases}
\]
where $B=\operatorname{Vert}(\tau)\setminus A$. Since $y_{v_i}=x_i$ for $i=1,\ldots,\nu-1$, the resulting family of functions contains the prescribed functions $x_1,\ldots,x_{\nu-1}$.
Throughout this part of the proof, we restrict all functions to $U\cap \operatorname{Int}(\tau)$, and write $f|$ for the restriction of $f$ to $U\cap \Int(\tau)$. 
We will show that the transformation  $(\delta_v|)_{v\in I}\mapsto (y_v|)_{v\in I}$ is a PL coordinate change on $\tau$ near $z$, which will prove Lemma~\ref{l:PL_coordinates} in the case $p=m$.
We write the transformation  $(\delta_v|)_{v\in I}\mapsto (y_v|)_{v\in I}$ as a composition of three elementary PL coordinate changes. Specifically, 
\[
 (\delta_v|)_{v\in I} \mapsto (r_1, \ldots, r_{\nu-1}, (y_v|)_{v\in I'}) \mapsto (h_1, \ldots h_{\nu-1}, (y_v|)_{v\in I'}) \mapsto (x_1|, \ldots, x_{\nu-1}|, (y_v|)_{v\in I'}), 
\]
where $I'=I\setminus \{v_1, \ldots, v_{\nu-1}\}$, and the functions $r_i$ and $h_i$ are defined below. 

For $1\le i\le \nu-1$, the functions $r_i:=\delta_{v_i}|-\delta_{v_\nu}|$ measure the differences between the chosen representatives. Put $r_\nu=0$. The first transformation is an injective affine map. Indeed, if two points $x$ and $x'$ in $U\cap \Int(\tau)$ are transformed to the same point, then $\delta_v(x)=\delta_v(x')$ for $v\in B$. 
For $v\in A_i\setminus \{v_i\}$, we have $\delta_v|=y_v|+r_i+ \delta_{v_\nu}|$, and therefore $\delta_v(x)-\delta_v(x')=\delta_{v_\nu}(x)-\delta_{v_\nu}(x')$. The sum of barycentric coordinates of a point in a simplex is $1$. 
Hence $\sum \delta_v(x)=\sum(1-\lambda_{v}(x))=p=\sum \delta_v(x')$, where the sums range over all vertices of $\tau$. Consequently,  
\[
 0 =\sum_{v\in \operatorname{Vert}(\tau)} (\delta_v(x)-\delta_v(x'))=\sum_{v\in A} (\delta_v(x)-\delta_v(x')) = |A|\cdot (\delta_{v_{\nu}}(x)-\delta_{v_\nu}(x')).
\]
Hence $\delta_{v_\nu}(x)=\delta_{v_\nu}(x')$,  and therefore $\delta_v(x)=\delta_v(x')$ for all $v\in A\setminus \{v_i\}$. The equality $r_i(x)=r_i(x')$ together with $\delta_{v_\nu}(x)=\delta_{v_\nu}(x')$
implies that $\delta_{v_i}(x)=\delta_{v_i}(x')$. Thus, the first transformation is invertible.  

Next, for each $i=1,\ldots, \nu$, we introduce functions $a_i$ such that $\delta_i=\delta_{v_i}+a_i$. Explicitly, 
\[
  a_i=\min \bigl(\{0\}\cup \{y_v: v\in A_i\setminus \{v_i\}\}\bigr). 
\]
To show that these functions satisfy $\delta_i=\delta_{v_i}+a_i$, we note that 
\[
\delta_i
=\min\bigl(\{\delta_{v_i}\}\cup
\{\delta_{v_i}+y_v:v\in A_i\setminus\{v_i\}\}\bigr)\\
=\delta_{v_i}
+\min\bigl(\{0\}\cup\{y_v:v\in A_i\setminus\{v_i\}\}\bigr) = \delta_{v_i}+a_i.
\]

For the next transformation, we set $h_i=\delta_i|-\delta_\nu|$ for $i=1, \ldots, \nu-1$. This transformation is a PL homeomorphism since $h_i=(\delta_{v_i}|+a_i|)-(\delta_{v_\nu}|+a_\nu|)=r_i+a_i|-a_\nu|$, and the functions $a_i|$ and $a_\nu|$ depend only on the coordinates $(y_v|)_{v\in I'}$ which remain unchanged. Finally, since
\[
  x_j|=\delta_j|-\min_{j<i\le \nu} \delta_i| = (\delta_j|-\delta_\nu|) - \min_{j<i\le \nu} \bigl(\delta_i|-\delta_\nu|\bigr) = h_j - \min\bigl(\{0\}\cup \{h_i: j<i<\nu\}\bigr). 
\]
the last transformation is unit upper triangular, and hence, a PL homeomorphism. Thus, the restrictions of $(y_v)_{v\in I}$ give $p$ PL coordinates on $\tau$ near $z$. 

If $p=m$, then the proof is complete. Suppose now that $p<m$. Put $L=\mathop\mathrm{Lk}_{K'}(\tau)$. Let $o*L$ be the cone on $L$ with cone point $o$. We use $\lambda_v$ for barycentric coordinates both on $M$ and $o*L$, and similarly for $\delta_v$. Choose an open neighborhood $U_\tau$ of $z$ in $\tau$ whose closure is contained in $U\cap \Int(\tau)$. Choose $\varepsilon\in (0,1)$ sufficiently small so that  $\lambda_v(u)>\varepsilon/(p+1)$ for every $u\in \overline{U}_\tau$, and every vertex $v$ of $\tau$. 
Let  
$W=\{\xi\in o*L:\delta_o(\xi)<\varepsilon\}$ be the $\varepsilon$-neighborhood of $o$ in the cone. Define a PL map
$F\co U_\tau \times W\longrightarrow M$ in terms of barycentric coordinates by 
\[
\begin{aligned}
\lambda_v(F(u,\xi))
   &=\lambda_v(u)-\frac{\delta_o(\xi)}{p+1},
   &&v\in\operatorname{Vert}(\tau),\\
\lambda_w(F(u,\xi))
   &=\lambda_w(\xi),
   &&w\in\operatorname{Vert}(L),
\end{aligned}
\]
with all remaining coordinates zero. We have $\sum \lambda_v(u)=1$ for each $u\in U_\tau$, where the sum ranges over vertices in $\operatorname{Vert}(\tau)$. Also $\delta_o(\xi)=1-\lambda_o(\xi)=\sum \lambda_w(\xi)$, where $w$ ranges over vertices in $L$. Therefore the coordinates $\lambda_v(F(u, \xi))$ and $\lambda_w(F(u, \xi))$ indeed sum to one. Choose a simplex $\sigma\in L$ such that $\xi\in o*\sigma$. Then  $\lambda_w(\xi)=0$ for every vertex $w$ outside $\sigma$ so only vertices $w$ of $\sigma$ can have positive coordinates $\lambda_w$. Therefore these coordinates determine a point of the simplex $\tau*\sigma$ of $K'$. The inverse of $F$ can be described explicitly in terms of barycentric coordinates. For a point $q$ in the open star of $\tau$, define $\xi=\xi(q)\in o*L$ by 
\[
\lambda_w(\xi)=\lambda_w(q),\qquad
\lambda_o(\xi)=1-\sum_{w\in\operatorname{Vert}(L)}\lambda_w(q), 
\]
and define $u=u(q)\in\operatorname{Int}(\tau)$ by
\[
\lambda_v(u)=\lambda_v(q)+\frac{\delta_o(\xi)}{p+1},
\qquad v\in\operatorname{Vert}(\tau).
\]
These formulas restricted to the image of $F$ define a PL inverse to $F$. Thus $F$ is a PL homeomorphism onto an open neighborhood $V$ of $z$. Decreasing $\varepsilon$ if necessary, we may assume that $V\subset U$. For every vertex $v$ in $\tau$, we have 
\[
\delta_v(F(u,\xi))
=\delta_v(u)+\frac{\delta_o(\xi)}{p+1}.
\]
Consequently, 
\[
x_j(F(u,\xi))
=\delta_j(F(u,\xi))-\min_{j<i\le \nu} \delta_i(F(u,\xi))
=\delta_j(u)-\min_{j<i\le \nu} \delta_i(u) 
=x_j(u),
\]
which means that the functions $x_j$, and therefore sector membership, are constant along the normal fibers $F(\{u\}\times W)$. 

The link $L$ is a PL sphere if $z$ is in the interior of $M$, and it is a ball if $z$ is on the boundary of $M$. Therefore $W$ admits $m-p$ PL coordinates. Under the map $F(u, \xi)$, the coordinates $y_v(u)$ together with normal coordinates on $W$ produce a coordinate system about $z$ containing the original functions $x_1, \ldots, x_{\nu-1}$. 
\end{proof}

For $\nu>1$, by Lemmas~\ref{l:N} and \ref{l:PL_coordinates}, the sectors $\{N_j\}$ satisfy the local model condition. If $\nu=1$, then all vertices in $A$ are of the same color $i=1$. In this case $\delta_1(z)<\delta_i(z)$ for every $i\ne 1$, and these inequalities still hold on a sufficiently small neighborhood $U'$ of $z$. Hence $U'\subset N_1$, and therefore a sufficiently small coordinate neighborhood of $z$ gives the required local model. The local models show that the sectors are PL submanifolds of codimension $0$ with boundary. Lemma~\ref{prop:canonical} implies that $M=N_1\cup\cdots\cup N_r$, and the interiors of the sectors $N_i$ are pairwise disjoint. Finally, the number of sectors containing a given point $z$ is $\nu\le |A|\le m+1$. Thus, the ordered collection $(N_j)$ is a PL filling when $M$ is closed. The case where $\partial M$ is non-empty is similar.  
\end{proof}

\section*{Tool use}

During the preparation of this paper, the author used ChatGPT (OpenAI) to assist with revising the exposition, developing and reviewing mathematical arguments, locating references, and preparing figures. The author takes full responsibility for the content of the paper.

\end{document}